
  \magnification 1200
  \input amssym  


  \newcount\fontset
  \fontset=1
  \def \dualfont#1#2#3{\font#1=\ifnum\fontset=1 #2\else#3\fi}

  \dualfont\bbfive{bbm5}{cmbx5}
  \dualfont\bbseven{bbm7}{cmbx7}
  \dualfont\bbten{bbm10}{cmbx10}

  \font \eightbf = cmbx8
  \font \eighti = cmmi8 \skewchar \eighti = '177
  \font \eightit = cmti8
  \font \eightrm = cmr8
  \font \eightsl = cmsl8
  \font \eightsy = cmsy8 \skewchar \eightsy = '60
  \font \eighttt = cmtt8 \hyphenchar\eighttt = -1

  \font \sixi = cmmi6 \skewchar \sixi = '177
  \font \sixrm = cmr6
  \font \sixsy = cmsy6 \skewchar \sixsy = '60
  \font \tensc = cmcsc10

  \scriptfont \bffam = \bbseven
  \scriptscriptfont \bffam = \bbfive
  \textfont \bffam = \bbten

  \newskip \ttglue

  \def \eightpoint {\def \rm {\fam0 \eightrm }%
  \textfont0 = \eightrm
  \scriptfont0 = \sixrm \scriptscriptfont0 = \fiverm
  \textfont1 = \eighti
  \scriptfont1 = \sixi \scriptscriptfont1 = \fivei
  \textfont2 = \eightsy
  \scriptfont2 = \sixsy \scriptscriptfont2 = \fivesy
  \textfont3 = \tenex
  \scriptfont3 = \tenex \scriptscriptfont3 = \tenex
  \def \it {\fam \itfam \eightit }%
  \textfont \itfam = \eightit
  \def \sl {\fam \slfam \eightsl }%
  \textfont \slfam = \eightsl
  \def \bf {\fam \bffam \eightbf }%
  \textfont \bffam = \bbseven
  \scriptfont \bffam = \bbfive
  \scriptscriptfont \bffam = \bbfive
  \def \tt {\fam \ttfam \eighttt }%
  \textfont \ttfam = \eighttt
  \tt \ttglue = .5em plus.25em minus.15em
  \normalbaselineskip = 9pt
  \def \MF {{\manual opqr}\-{\manual stuq}}%
  \let \sc = \sixrm
  \let \big = \eightbig
  \setbox \strutbox = \hbox {\vrule height7pt depth2pt width0pt}%
  \normalbaselines \rm }



  \newcount \secno \secno = 0
  \newcount \stno \stno =0
  \newcount \eqcntr \eqcntr=0

  \def \ifn #1{\expandafter \ifx \csname #1\endcsname \relax }

  \def \track #1#2#3{\ifn{#1}\else {\tt\ [#2 \string #3] }\fi}

  \def \advseqnumbering {\global \advance \stno by 1 \global \eqcntr=0}

  \def \current {\number \secno \ifnum \number \stno = 0 \else
    .\number \stno \fi }

  \def \laberr#1#2{\message{*** RELABEL CHECKED FALSE for #1 ***}
      RELABEL CHECKED FALSE FOR #1, EXITING.
      \end}

  \def \syslabel#1#2{%
    \ifn {#1}%
      \global \expandafter 
      \edef \csname #1\endcsname {#2}%
    \else
      \edef\aux{\expandafter\csname #1\endcsname}%
      \edef\bux{#2}%
      \ifx \aux \bux \else \laberr{#1=(\aux)=(\bux)} \fi
      \fi
    \track{showlabel}{*}{#1}}

  \def \subeqmark #1 {\global \advance\eqcntr by 1
    \edef\aux{\current.\number\eqcntr}
    \eqno {(\aux)}
    \syslabel{#1}{\aux}}

  \def \eqmark #1 {\advseqnumbering
    \eqno {(\current)}\syslabel{#1}{\current}}

  \def \fcite#1#2{\syslabel{#1}{#2}\lcite{#2}}

  \def \label #1 {\syslabel{#1}{\current}}

  \def \lcite #1{(#1\track{showcit}{$\bullet$}{#1})}

  \def \cite #1{[{\bf #1}\track{showref}{\#}{#1}]}

  \def \scite #1#2{{\rm [\bf #1\track{showref}{\#}{#1}{\rm \hskip 0.7pt:\hskip 2pt #2}\rm]}}


 \def \Headlines #1#2{\nopagenumbers
    \advance \voffset by 2\baselineskip
    \advance \vsize by -\voffset
    \headline {\ifnum \pageno = 1 \hfil
    \else \ifodd \pageno \tensc \hfil \lcase {#1} \hfil \folio
    \else \tensc \folio \hfil \lcase {#2} \hfil
    \fi \fi }}

  \def \Date #1 {\footnote {}{\eightit Date: #1.}}


  \def \lcase #1{\edef \auxvar {\lowercase {#1}}\auxvar }

  \def \goodbreak {\vskip0pt plus.1\vsize \penalty -250 \vskip0pt
plus-.1\vsize }

  \def \section #1{\global\def \SectionName{#1}\stno = 0 \global
\advance \secno by 1 \bigskip \bigskip \goodbreak \noindent {\bf
\number \secno .\enspace #1.}\medskip \noindent \ignorespaces}

  \long \def \sysstate #1#2#3{%
    \advseqnumbering
    \medbreak \noindent 
    {\bf \current.\enspace #1.\enspace }{#2#3\vskip 0pt}\medbreak }
  \def \state #1 #2\par {\sysstate {#1}{\sl }{#2}}
  \def \definition #1\par {\sysstate {Definition}{\rm }{#1}}
  \def \remark #1\par {\sysstate {Remark}{\rm }{#1}}


  \def \proof {\medbreak \noindent {\it Proof.\enspace }}
  \def \proofend {\ifmmode \eqno \square \else \hfill \square
\looseness = -1 \medbreak \fi }

  \def \$#1{#1 $$$$ #1}
  \def \=#1{\buildrel \hbox{\sixrm #1} \over =}

  \def \pilar #1{\vrule height #1 width 0pt}

  \def \Item #1{\smallskip \item {{\rm #1}}}
  \newcount \zitemno \zitemno = 0

  \def \izitem {\zitemno = 0}
  \def \zitemplus {\global \advance \zitemno by 1\relax}
  \def \rzitem{\romannumeral \zitemno}
  \def \rzitemplus {\zitemplus \rzitem}
  \def \zitem {\Item {{\rm(\rzitemplus)}}}

  \newcount \nitemno \nitemno = 0
  
  \def \nitem {\global \advance \nitemno by 1 \Item {{\rm(\number\nitemno)}}}

  \newcount \aitemno \aitemno = 0
  \def\boxlet#1{\hbox to 6.5pt{\hfill #1\hfill}}
  \def \iaitem {\aitemno = 0}
  \def \aitem {\Item {(\ifcase \aitemno \boxlet a\or \boxlet b\or
\boxlet c\or \boxlet d\or \boxlet e\or \boxlet f\or \boxlet g\or
\boxlet h\or \boxlet i\or \boxlet j\or \boxlet k\or \boxlet l\or
\boxlet m\or \boxlet n\or \boxlet o\or \boxlet p\or \boxlet q\or
\boxlet r\or \boxlet s\or \boxlet t\or \boxlet u\or \boxlet v\or
\boxlet w\or \boxlet x\or \boxlet y\or \boxlet z\else zzz\fi)} \global
\advance \aitemno by 1}

  \newcount \footno \footno = 1
  \newcount \halffootno \footno = 1
  \def \footcntr {\global \advance \footno by 1
  \halffootno =\footno
  \divide \halffootno by 2
  $^{\number\halffootno}$}
  \def \fn#1{\footnote{\footcntr}{\eightpoint#1\par}}

  \begingroup
  \catcode `\@=11
  \global\def\eqmatrix#1{\null\,\vcenter{\normalbaselines\m@th%
      \ialign{\hfil$##$\hfil&&\kern 5pt \hfil$##$\hfil\crcr%
	\mathstrut\crcr\noalign{\kern-\baselineskip}%
	#1\crcr\mathstrut\crcr\noalign{\kern-\baselineskip}}}\,}
  \endgroup


  \font\mf=cmex10
  \def\union {\mathop{\raise 9pt \hbox{\mf S}}\limits}
  \def\inters{\mathop{\raise 9pt \hbox{\mf T}}\limits}
  
  \def\bigcap{\inters}

  \def \N {{\bf N}}
  
  \def \<{\left \langle \vrule width 0pt depth 0pt height 8pt }
  \def \>{\right \rangle }  
  
  \def \ds{\displaystyle}
  \def \and {\hbox {,\quad and \quad }}
  \def \calcat #1{\,{\vrule height8pt depth4pt}_{\,#1}}
  \def \labelarrow#1{\ {\buildrel #1 \over \longrightarrow}\ }
  \def \imply {\mathrel{\Rightarrow}}
  \def \for #1{,\quad \forall\,#1}
  \def \square {\hbox {$\sqcap \!\!\!\!\sqcup $}}
  
  \def \stress #1{{\it #1}\/}
  \def \inv {^{-1}}
  \def \*{\otimes}

  \newcount \bibno \bibno =0
  \def \newbib #1{\global\advance\bibno by 1 \edef #1{\number\bibno}}
  \def \bibitem #1#2#3#4{\smallskip \item {[#1]} #2, ``#3'', #4.}
  \def \references {
    \begingroup
    \bigskip \bigskip \goodbreak
    \eightpoint
    \centerline {\tensc References}
    \nobreak \medskip \frenchspacing }

  \input pictex

  \newcount\ax
  \newcount\ay
  \newcount\bx
  \newcount\by
  \newcount\dx
  \newcount\dy
  \newcount\vecNorm
  \newcount\pouquinho \pouquinho = 200

  \def\beginmypicture{ \begingroup \noindent \hfill \beginpicture}
  \def\endmypicture{\endpicture \hfill\null \endgroup}
  \def\arwlabel#1#2#3{\put{${\scriptstyle #1}$} at #2 #3}

  \def\myarrow#1#2#3#4{\arrow <0.15cm> [0.25,0.75] from #1 #2 to #3 #4 }%
  \def\morph#1#2#3#4{%
    \ax = #1
    \ay = #2
    \bx = #3
    \by = #4
    \dx = \bx \advance \dx by -\ax
    \dy = \by \advance \dy by -\ay
    \vecNorm = \dx 
    \ifnum\vecNorm<0 \vecNorm=-\vecNorm \fi
    \advance \vecNorm by \ifnum\dy>0 \dy \else -\dy \fi
    \multiply \dx by \pouquinho \divide \dx by \vecNorm
    \multiply \dy by \pouquinho \divide \dy by \vecNorm
    \advance \ax by \dx
    \advance \bx by -\dx
    \advance \ay by \dy
    \advance \by by -\dy
    \myarrow{\number\ax}{\number\ay}{\number\bx}{\number\by}}

  \font \tensc = cmcsc10 
  \font\rs=rsfs10
  \font\rssmall=rsfs8
  \def\mycalfont{\ifdim\normalbaselineskip = 12pt \rs \else \rssmall \fi}

  \def\sysGexp#1#2{#1^{(#2)}}
  \def\Gexp#1{\sysGexp{\G}{#1}}
  \def\tn#1{\def\aux{\kern-1.4pt|}|\aux\aux#1|\aux\aux}
  \def\d{\delta}
  \def\lft{\delta}
  \def\rgt{\delta}
  
  \def\a{\alpha}
  \def\iota{i}

  \def\c{\subseteq}
  \def\leq{\leqslant}
  \def\geq{\geqslant}
  \def\emptyset{\varnothing}
  \def\smallcalbox#1{{\hbox{\rssmall #1}\,}}
  \def\calbox#1{{\hbox{\mycalfont #1}\,}}
  \def\cryout#1{\removelastskip \medskip {\bf #1} \medskip}


  \def\S{{\cal S}}
  \def\G{{\cal G}}   \def\H{{\cal H}}
  \def\E{E} 
  \def\Ehat{\hat E}
  \def\spec{\Ehat}
  \def\specz{\Ehat_0}
  \def\speci{\Ehat_\infty}
  
  \def\tightbox{{\hbox{\eightsl tight}}}
  \def\spect{\Ehat_\tightbox}
  \def\Gt{\G_\tightbox}
  \def\B{\calbox{B}}
  \def\I{{\cal I}}
  \def\P{\calbox{P}}
  \def\Psmall{\smallcalbox{P}}
  \def\Pregerm{\Omega}

  \def\inf{\wedge}
  \def\nega{\neg\,} 
  \def\p#1{#1#1^*}
  \def\q#1{#1^*#1}
  \def\O{\Theta}

  \def\its{\Cap}

  \def\r{\hbox{\bf r}} \def\r{r}
  \def\s{\hbox{\bf d}}   \def\s{d}

  \def\M{\Lambda}
  \def\Mu{\tilde\M}
  \def\Mt{\M^{(2)}}
  \def\Catego{\calbox{C}}%
  \def\stem{\omega}
  \def\D#1{\M^{#1}}
  \def\bitem{\Item{$\bullet$}}
  \def\dil{\mathrel{|}}
  \def\disj{\perp}
  \def\lcm{{\rm lcm}}
  \def\Q{\calbox{Q}}
  \def\SM{\S(\M)}
  \def\SG{\S(\G)}
  \def\orep{\pi} 
  \def\irep{\sigma} 
  \def\brep{\beta} 
  \def\act{\theta}  
  \def\ini{q}
  \def\iniq{Q}
  \def\fin{p}
  \def\Dom{\calbox{D}}

  \def \newbib #1#2{\global\advance\bibno by 1 \edef #1{#2}}
  \def \newbib #1#2{\global\advance\bibno by 1 \edef
#1{\number\bibno}}

  \newbib\AbadieGroupoid{A}
  \newbib\BHRS{BHRS}  
  \newbib\BPRS{BPRS}  
  \newbib\Connes{Co}
  \newbib\Cuntz{Cu}
  \newbib\CK{CK}
  \newbib\Watatani{EW}  
  \newbib\newpim{E1}
  \newbib\tpa{E2}
  \newbib\ExSemiGpdAlg{E3}
  \newbib\ExelAlgebra{E4}
  \newbib\infinoa{EL}
  \newbib\Muhly{FMY}  
  \newbib\FLR{FLR}  
  \newbib\KPActions{KP1}  
  \newbib\KP{KP2}  
  \newbib\KPZActions{KP3}  
  \newbib\KPR{KPR}  
  \newbib\KPPR{KPPR}  
  \newbib\Lawson{L}
  \newbib\McClanahan{McC}
  \newbib\PQR{PQR}  
  \newbib\PRRS{PRRS}  
  \newbib\PatBook{P1}
  \newbib\PatGraph{P2} 
  \newbib\QuigSieb{QS}
  \newbib\RaeBook{Ra}  
  \newbib\RSY{RSY}  
  \newbib\RaeSzy{RS}  
  \newbib\PTW{PTW}  
  \newbib\RenaultThesis{Re1}
  \newbib\RenaultInfinoa{Re2}
  \newbib\Sieben{Si}
  \newbib\Szendrei{Sz}
  \newbib\Tomforde{T} 


  \def\titletextOne{INVERSE SEMIGROUPS AND COMBINATORIAL}
  \def\titletextTwo{C*-ALGEBRAS}

  \Headlines {\titletextOne \ \titletextTwo} {R.~Exel}

  \null\vskip -1cm
  \centerline{\bf \titletextOne} \smallskip
  \centerline{\bf \titletextTwo}
  \footnote{\null}
  {\eightrm 2000 \eightsl Mathematics Subject Classification:
  \eightrm 
  Primary 46L05; 
  secondary 18B40. 
  }

  \bigskip
  \centerline{\tensc 
    R.~Exel\footnote{*}{\eightpoint Partially supported by
CNPq.}}

  \bigskip
  \Date{29 Feb 2008}

  \null\footnote{\null}
  {\eightrm Keywords: C*-algebras,
  Cuntz-Krieger algebras,
  graphs,
  higher-rank graphs,
  groupoids,
  inverse semigroups,
  semilattices,
  ultra-filters,
  boolean algebras,
  tight Hilbert space representations,
  crossed products,
  germs,
  semigroupoids,
  categories.
  }

  \midinsert 
  \narrower \narrower
  \eightpoint \noindent We describe a special class of representations
of an inverse semigroup $\S$ on Hilbert's space which we term
\stress{tight}.  These representations are {supported} on a subset of
the spectrum of the idempotent semilattice of $\S$, called the
\stress{tight spectrum}, which is in turn shown to be precisely the
closure of the space of ultra-filters, once filters are identified
with semicharacters in a natural way.  These representations are
moreover shown to correspond to representations of the C*-algebra of
the groupoid of germs for the action of $\S$ on its tight spectrum.
  We then treat the case of certain inverse semigroups constructed
from semigroupoids, generalizing and inspired by inverse semigroups
constructed from ordinary and higher rank graphs.  The tight
representations of this inverse semigroup are in one-to-one
correspondence with representations of the semigroupoid, and
consequently the semigroupoid algebra is given a groupoid model.  The
groupoid which arises from this construction is shown to be the same
as the \stress{boundary path} groupoid of Farthing, Muhly and Yeend,
at least in the singly aligned, sourceless case.
  \endinsert

\section{Introduction}
  By a \stress{combinatorial} C*-algebra we loosely mean any
C*-algebra which is constructed from a combinatorial object.  Among
these we include the Cuntz-Krieger algebras built out of 0--1
matrices, first studied in the finite case in \cite{\CK}, and quickly
recognized to pertain to the realm of Combinatorics by Enomoto and
Watatani \cite{\Watatani}.  Cuntz-Krieger algebras were subsequently
generalized to row-finite matrices in \cite{\KPPR}, and to general
infinite matrices in \cite{\infinoa}.  Another important class of
combinatorial C*-algebras, closely related to the early work of Cuntz
and Krieger, is formed by the graph C*-algebras 
  \cite{%
  \BHRS,
  \BPRS,
  \FLR,
  \KPActions,
  \KPR,
  \PatGraph,
  \RaeSzy,
  \PTW,
  \Tomforde
  }, including the case of higher rank graphs introduced by Kumjian
and Pask in \cite{\KP}, and given its final form by Farthing, Muhly
and Yeend in \cite{\Muhly}. See also
  \cite{%
  \KPZActions,
  \PQR,
  \PRRS,
  \RSY
  }.  
  The monograph \cite{\RaeBook} is an excellent source of well
organized information and references.

Attempting to understand all of these algebras from a single
perspective, I have been interested in the notion of semigroupoid
C*-algebras \cite{\ExSemiGpdAlg}, which includes the Cuntz-Krieger
algebras and the higher rank graph C*-algebras in the general infinite
case, provided some technical complications are not present including,
but not limited to, \stress{sources}.

The most efficient strategy to study combinatorial C*-algebras has
been the construction of a dynamical system which intermediates
between Combinatorics and Algebra.  In the case of \cite{\infinoa},
the dynamical system took the form of a partial action of a free group
on a topological space, but more often it is represented by an
\'etale, or $r$-discrete groupoid. In fact, even in the case of
\cite{\infinoa}, the partial action may be encoded by a groupoid
\cite{\AbadieGroupoid}, \cite{\RenaultInfinoa}.
  It therefore seemed natural to me that semigroupoid C*-algebras
could also be given groupoid models.  But, unfortunately, the
similarity between the terms \stress{semigroupoid} and
\stress{groupoid} has not made the task any easier.

The vast majority of combinatorial C*-algebras may be defined
following a standard pattern: the combinatorial object chosen is used
to suggest a list of relations, written in the language of
C*-algebras, and then one considers the universal C*-algebra generated
by partial isometries satisfying such relations.

Partial isometries can behave quite badly from an algebraic point of
view, and in particular the product of two such elements needs not be
a partial isometry.  Should the most general and wild partial
isometries be involved in combinatorial algebras, the study of the
latter would probably be impossible.  Fortunately, though, the partial
isometries one usually faces are, without exception, of a tamer nature
in the sense that they always generate a *-semigroup consisting
exclusively of partial isometries or, equivalently, an inverse
semigroup.

In two recent works, namely \cite{\PatGraph} and \cite{\Muhly}, this
inverse semigroup has been used in an essential way, bridging the
combinatorial input object and the groupoid.

  \medskip\beginmypicture
  \setcoordinatesystem units <0.0010truecm, -0.001truecm> point at 0 0
  \putrectangle corners at 0000 0000 and 2500 800
  \put {\eightrm Combinatorial} at 1250 200
  \put {\eightrm Object} at 1250 550
  \arrow <0.15cm> [0.25,0.75] from 2600 400 to 3400 1400
  \putrectangle corners at 3500 1000 and 6000 1800
  \put {\eightrm Inverse} at 4750 1200
  \put {\eightrm Semigroup} at 4750 1550
  \arrow <0.15cm> [0.25,0.75] from 6100 1400 to 6900 400
  \putrectangle corners at 7000 0000 and 9500 800
  \put {\eightrm Groupoid} at 8250 400
  \arrow <0.15cm> [0.25,0.75] from 9600 400 to 10400 400
  \putrectangle corners at 10500 0000 and 13000 800
  \put {\eightrm Combinatorial} at 11750 200
  \put {\eightrm C*-algebra} at 11750 550
  \endmypicture

  \vskip 0.5cm
  \def\Strategy{1.1}
  \centerline{\eightrm Diagram \Strategy}
  \vskip 0.3cm

\noindent
  In both \cite{\PatGraph} and in \cite{\Muhly} the relevant inverse
semigroup is made to act on a topological space by means of partial
homeomorphisms.  The groupoid of germs for this action then turns out
to be the appropriate groupoid.
  However, the above diagram does not describe this strategy quite correctly
because the topological space where the inverse semigroup acts is a space of
\stress{paths} whose description requires that one looks back at the
combinatorial object.

Attempting to adopt this strategy, I stumbled on the fact that it is
very difficult to guess the appropriate path space in the case of a
semigroupoid.  Moreover, earlier experience with partial actions of
groups suggested that the path space should be intrinsic to the
inverse semigroup.  

Searching the literature one indeed finds intrinsic dynamical systems
associated to a given inverse semigroup $\S$, such as the natural
action of $\S$ on the semicharacter space of its idempotent
semilattice \scite{\PatBook}{Proposition 4.3.2}.  But, unfortunately,
the groupoid of germs for this action turns out not to be the correct
one.  For example, if one starts with the most basic of all
combinatorial algebras, namely the Cuntz algebra ${\cal O}_n$, the
appropriate combinatorial object is an $n\times n$ matrix of zeros and
ones, which in this case consists only of ones, and the inverse
semigroup is the Cuntz inverse semigroup, as defined by Renault in
\scite{\RenaultThesis}{III.2.2}.  But the groupoid of germs
constructed from the above intrinsic action is not the Cuntz groupoid
because its C*-algebra is the Toeplitz extension of ${\cal O}_n$
\scite{\Cuntz}{Proposition 3.1}, rather than ${\cal O}_n$ itself.  See
also \scite{\RenaultThesis}{III.2.8.i}.

If $E=E(\S)$ is the idempotent semilattice of an inverse semigroup
$\S$, one says that a nonzero map
  $$
  \phi:E\to\{0,1\}
  $$
  is a \stress{semicharacter} if $\phi(x y) = \phi(x) \phi(y)$, for
all $x$ and $y$ in $E$.  The set of all semicharacters, denoted
$\spec$, is called the \stress{semicharacter space} of $E$, and it is
a locally compact topological space under the topology of pointwise
convergence.  The intrinsic action we referred to above is a certain
very natural action of $\S$ on $\spec$.

If $\S$ contains a zero element 0, a quite common situation which can
otherwise be easily arranged, then 0 is in $E$ but the above popular
definition of semicharacter strangely does not require that
$\phi(0)=0$.  In fact the space of all semicharacters is too big, and
this is partly the reason why $\spec$ does not yield the correct
groupoid in most cases.  This is also clearly indicated by the need to
reduce the universal groupoid in \cite{\PatGraph}.

One of the main points of this work is that this reduction can be
performed in a way that is entirely intrinsic to $\S$, and does not
require any more information from the combinatorial object which gave
rise to $\S$.  In other words, the diagram above can be made to work
exactly as indicated.

  If $\phi$ is a semicharacter of a semilattice $E$ then the set
  $$
  \xi = \xi_\phi =\{e\in E: \phi(e)=1\},
  $$
  which incidentally characterizes $\phi$, is a \stress{filter} in the
sense that it contains $ef$, whenever $e$ and $f$ are in $\xi$, and
moreover $f\geq e\in\xi$ implies that $f\in\xi$.  In case $\S$
contains 0, and one chooses to add to the definition of semicharacters
the sensible requirement that $\phi(0)$ should be equal to zero then,
in addition to the above properties of $\xi$, one gets $0\notin\xi$.
We then take these simple properties as the definition of a filter.

With the exception of Kellendonk's topological groupoid
\scite{\Lawson}{9.2}, most authors have not paid too much attention to
the fact that ultra-filters form an important class of filters, and
that these are present in abundance, thanks to Zorn's Lemma.
Kellendonk's treatment is however not precisely what we need, perhaps
because of the reliance on sequences with only countable many terms.

It then makes sense to pay attention to the set $\speci$ formed by all
semicharacters $\phi$ for which $\xi_\phi$ is an ultra-filter. Our
apology of ultra-filters notwithstanding, $\speci$ is not always
tractable by the methods of Topology not least because it may fail to
be closed in $\spec$.  But in what follows we will try to convince the
reader that the closure of $\speci$ within $\spec$, which we denote by
$\spect$, is the right space to look at.

This explains several
instances in the literature where \stress{finite paths} shared the
stage with \stress{infinite paths}.  Not attempting to compile a
comprehensive list, we may cite as examples, in chronological order:

  \bitem The description of the spectrum of the Cuntz--Krieger
relations for arbitrary 0--1 matrices given at the end of 
\scite{\infinoa}{Section 5}.  See also \scite{\infinoa}{7.3}.

  \bitem Paterson's description of the unit space of the path
groupoid of a graph \scite{\PatGraph}{Proposition 3}.
See also \scite{\PatGraph}{Proposition 4}.

  \bitem The closed invariant space $\partial\Lambda$ within the
space of all finite and infinite paths in a higher rank graph $\Lambda$,
constructed by Farthing, Muhly and Yeend in \scite{\Muhly}{Definition
5.10}.  See also \scite{\Muhly}{Theorem 6.3}.

\bigskip
  To fully explain the connection between the groupoid of germs for
the natural action of $\S$ on $\spect$ and the above works would make
this paper even longer than it already is, so we have opted instead to
restrict attention to semigroupoid C*-algebras.  On the one hand these
include most of the combinatorial algebras mentioned so far, but on
the other hand we have made significant restrictions in order to fend
off well known technical complications.

While we do not compromise on \stress{infiniteness}, we assume that our
semigroupoid has no \stress{springs}, and \stress{admits least common
multiples}.  These hypotheses correspond, in the case of a higher rank
graph $\M$, to the absence of \stress{sources}, and to the fact that
$\M$ is \stress{singly aligned}. Besides allowing for technical
simplifications, the existence of least common multiples evokes
important connections to Arithmetics, and has an important geometrical
interpretation in higher rank graphs.

When $\M$ is a semigroupoid satisfying all of the these favourable
hypotheses, we construct an inverse semigroup $\SM$, and then prove in
Theorem \fcite{CstarAlgOfGpd}{18.4} that the semigroupoid C*-algebra
is isomorphic to the C*-algebra for the groupoid of germs for the
natural action of $\SM$ on $\spect$, exactly following the strategy
outlined in Diagram \lcite{\Strategy}.

Although we have not invested all of the necessary energy to study the
inverse semigroup constructed from a general higher rank graph, as in
\cite{\Muhly}, we conjecture that the groupoid there denoted by
$\G_\Lambda|_{\partial\Lambda}$ is the same as the groupoid $\Gt$ of
Theorem \fcite{EquivTightRepISGandGPD}{13.3} below, or at least our
findings seem to give strong indications that this is so.  Should this
be confirmed, the assertion made in the introduction of \cite{\Muhly}
that their groupoid is \stress{fairly far removed from the universal
groupoid of $\S_\Lambda$} might need rectification.

The first part of this work, comprising Sections
  \fcite{GroupoidSection}{3}--\fcite{ActOnSpecSect}{10}
  is based on Renault's Thesis \cite{\RenaultThesis} and Paterson's
book \cite{\PatBook}, and should be considered as a survey of the
technical methods we use in the subsequent sections, beginning with a
careful study of non-Hausdorff \'etale groupoids and their
C*-algebras.  We also discuss actions of inverse semigroups on
topological spaces and describe the associated groupoid of germs in
detail.  Sieben's theory of crossed products by inverse semigroups
\cite{\Sieben} is included.

We have made a special effort to assume as few hypotheses as possible,
and this was of course facilitated by our focus on \'etale groupoids.
We hope this can be used as a guide to the beginner who
is primarily interested in the \'etale case and hence needs not spend
much energy on Haar systems.

As a result of our economy of assumptions we have found
generalizations of some known results, most notably Theorem
\fcite{GenQSP}{9.9} below, which shows that the C*-algebra of an
\'etale groupoid is a crossed product in Sieben's sense, even in the
non-Hausdorff case, with much less stringent hypotheses than the
\stress{additivity} assumption of \scite{\PatBook}{Theorem 3.3.1} or
the \stress{fullness} condition of \scite{\QuigSieb}{8.1}.
  We also present a minor improvement on \scite{\PatBook}{Proposition
3.3.3}, by removing the need for condition (ii) of
\scite{\PatBook}{Definition 3.3.1}.  This is presented in Proposition
\fcite{CoressponReps}{9.7} below.

Even though we do most of our work based on non-Hausdorff groupoids,
we have found an interesting sufficient condition for the groupoid of
germs to be Hausdorff, related to the order structure of inverse
semigroups.
  We show in Theorem \fcite{HausdorffGPG}{6.2} that if the inverse
semigroup $\S$ is a semilattice with respect to its natural order
  $$
  s\leq t \iff s = ts^*s,
  $$
  then every action of $\S$ on a locally compact Hausdorff space, for
which the domains of the corresponding partial homeomorphisms are
clopen, one has that the associated groupoid of germs is Hausdorff.

A special class of semigroups possessing the above mentioned property
  (see \fcite{ISGLattice}{6.4}) is formed by the $E^*$-unitary inverse
semigroups, sometimes also called $0$-$E$-unitary, which was
defined by Szendrei \cite{\Szendrei} and has been intensely studied in
the semigroup literature.  See, for example, \scite{\Lawson}{Section
9}.
  Kellendonk's topological groupoid is Hausdorff when $\S$ is
$E^*$-unitary \scite{\Lawson}{9.2.6}, and the related class of
$E$-unitary inverse semigroups have also been shown to provide
Hausdorff groupoids \scite{\PatBook}{Corollary 4.3.2}.


It is with section
  \fcite{RepSemilatSec}{11}
  that our original work takes off, where we develop the crucial
notion of \stress{tight} representations of a semilattice in a Boolean
algebra \fcite{DefLatTightRep}{11.6}.  Strangely enough, it is in the
realm of these very elementary mathematical constructs that we have
found the most important ingredient of this paper.  The concept of
tight representations may be considered a refinement of an idea which
has been dormant in the literature for many years, namely condition
(1.3) in \cite{\infinoa}.

In the following section we study representations of a semilattice
into the Boolean algebra $\{0,1\}$, and its relation to filters and
ultra-filters.  The central result, Theorem
\fcite{ClosureOfUltrafilters}{12.9}, is that the space of tight
characters is precisely the closure of the set of characters
associated to ultra-filters.  We also show in
\fcite{TightInvariance}{12.11} that tight characters on the idempotent
semilattice of an inverse semigroup $\S$ are preserved under the
natural action of $\S$, thus giving rise to the action of $\S$ on
$\spect$, the dynamical system which occupies our central stage.

In the short section 
  \fcite{TightRepISGSect}{13}
  we consider tight Hilbert space representations of a given inverse
semigroup $\S$ and show in \fcite{EquivTightRepISGandGPD}{13.3} that
they are in one-to-one correspondence to the representations of
$C^*(\Gt)$, where $\Gt$ is the groupoid of germs associated to the
natural action of $\S$ on the tight part of the spectrum of its
idempotent semilattice.  Perhaps this result could be interpreted as
saying that the C*-algebra generated by the range of a universal tight
representation of $\S$, which is isomorphic to $C^*(\Gt)$ by the
result mentioned above, is an important alternative to the classical
C*-algebra of an inverse semigroup studied, e.g.~in
\scite{\PatBook}{2.1}.  We also believe this addresses the concern
expressed by Renault in \scite{\RenaultThesis}{III.2.8.i}.

From section
  \fcite{ISGFromGPDSec}{14}
  onwards we start our study of the C*-algebra of a semigroupoid $\M$
and, as in \cite{\PatGraph} and \cite{\Muhly}, the first step is to
construct an inverse semigroup, which we denote by $\SM$.  This could
be thought of as traversing the leftmost arrow of Diagram
\lcite{\Strategy}.

In sections 
  \fcite{RepSGPDSect}{15}--\fcite{ExtRepSec}{17}
  we show that tight representations of $\SM$ correspond to certain
representations of $\M$, which by abuse of language we also call
\stress{tight}.  This step is crucial because the definition of the
semigroupoid C*-algebra naturally emphasizes the semigroupoid, rather
than its associated inverse semigroup, so one needs to be able to
determine tightness by looking at $\M$ only.

In the following section we essentially piece together the results so
far obtained to arrive at another main result, namely Theorem
  \fcite{CstarAlgOfGpd}{18.4},
  where we show that the semigroupoid C*-algebra is isomorphic to the
groupoid C*-algebra of $\Gt$, where where $\Gt$ is the groupoid of
germs associated to the natural action of $\S$ on the tight part of
the spectrum of its idempotent semilattice.

Our approach has an aesthetical advantage over \cite{\PatGraph} or
\cite{\Muhly} in the sense that our groupoid is constructed based on a
very simple idea which can be conveyed in a single sentence, namely
that one needs to focus on the set of ultra-filters, and necessarily
also on the filters in its boundary.  The disadvantage is that it
leads to a very abstract picture of our groupoid and one may argue
that a more concrete description is desirable.

We believe this concern may be addressed in the most general
situation, and a dynamical system much like the one studied in
\cite{\infinoa} will certainly emerge, although, rather than a partial
action of a group, it will be an action of an inverse semigroup.
Given the widespread interest in combinatorial objects taking the form
of a category, we instead specialize in section
  \fcite{CategoricalSec}{19}
  to \stress{categorical} semigroupoids, as defined in
\fcite{DefineCatSGPD}{19.1}.  This notion captures an essential
property of categories which greatly simplifies the study of $\SM$,
and hence also of $\Gt$.  In Proposition
\fcite{CharacTightInCatego}{19.12} we then give a simple
characterization of tight characters, resembling very much the
description of boundary paths of \cite{\Muhly}.

In the closing section we focus directly on higher rank graphs and
some effort is spent to determine which such objects lead to a
semigroupoid admitting least common multiples.  Not surprisingly only
\stress{singly aligned} higher rank graphs pass the test, in which
case we may apply our machinery, arriving at groupoid model of
its C*-algebra.

The literature is rich in very interesting examples of inverse
semigroups, such as certain inverse semigroups associated to tilings
\scite{\Lawson}{9.5}, and we believe it might be a very fulfilling
task to explore some of these from the point of view of tight
representations.

We would like to thank Aidam Sims for bringing to our attention many
relevant references on the subject of Higher Rank Graphs.



%
%

\section{A quick motivation}
  Let us now briefly discuss the example of the Cuntz inverse semigroup
\cite{\Cuntz}, \scite{\RenaultThesis}{III.2.2}, since it is the one of
the main motivations for this work.  The reader is invited to keep
this example in mind throughout.

To avoid unnecessary complications
we will restrict ourselves to the case $n=2$. 

 Consider the semigroup
$S$ consisting of an identity 1, a zero element 0, and all the words
in four letters, namely $p_1,p_2, q_1, q_2$, subject to the relations
$q_jp_i=\delta_{i, j}$.

It is shown in \scite{\Cuntz}{1.3} that every element of $S$ may be
uniquely written as 
  $$
  p_{i_1}\ldots p_{i_k}  \ q_{j_l}\ldots p_{j_1},
  $$
  where $k, l\geq 0$, and $i_1,\ldots,i_k,j_1,\ldots,j_l\in \{1,2\}$.
It turns out that $S$ is an inverse semigroup with $1^*=1$, $0^*=0$,
and $p_i^* = q_i$.

Given a nondegenerated representation $\irep$ of $S$ on a Hilbert
space $H$, vanishing on $0$, put $S_i=\irep(p_i)$.  It is then
elementary to prove that the $S_i$ satisfy
  $$
  S_i^*S_j = \delta_{i, j},
  $$
  which means that the $S_i$ are isometries on $H$ with pairwise
orthogonal ranges.  Conversely, given any two isometries on $H$ with
pairwise orthogonal ranges one may prove that there is a unique
representation $\irep$ of $S$ such that $\irep(p_i)=S_i$.  In other
words the representations of $S$ are in one-to-one correspondence with the
pairs of isometries having orthogonal ranges.

If the reader is acquainted with the Cuntz algebra ${\cal O}_2$
he or she will
likely wonder under which conditions on $\irep$ does the relation
  $$
  S_1S_1^*+S_2S_2^*=1 
  \eqno{(\dagger)}
  $$
  also holds.  After fiddling a bit whith this question one will
realize that the occurence of the plus sign above is not quite in
accordance with the language of semigroups (in which one only has the
multiplication operation).  In other words it is not immediately clear
how to state $(\dagger)$ in the language of semigroups.

In order to approach this problem first notice that the
idempotent semillatice $E(S)$ consists of 0, 1, and the idempotents
  $$
  e_{i_1, \ldots, i_k}:=
  p_{i_1}\ldots p_{i_k}  \ q_{i_k}\ldots q_{i_1},
  $$
  where  $i_1,\ldots,i_k\in \{1,2\}$.
Clearly 1 is the largest element of $E(S)$, while 0 is the
smallest.  

Next observe that any nonzero idempotent $f$ \stress{intersects} either
$e_1$ or $e_2$, in the sense that either $fe_1$ or $fe_2$ is nonzero.
The set of idempotents $\{e_1,e_2\}$ is therefore what we shall call a
\stress{cover} for $E(S)$.
One could try to indicate this fact by saying that 
  $$
  e_1\vee e_2=1,
  $$
  except that semillatices are only equipped with a meet operation
``$\wedge$", rather than a join operation ``$\vee$" as we seem to be
in need of.   However any (meet)-semilattice of projections on a Hilbert
space is contained in a Boolean algebra of projections, in which a
join operation is fortunately available, namely
  $$
  p\vee q = p+q-pq.
  $$
  Given a representation $\irep$ of $S$ on a Hilbert
space, one might 
therefore impose the condition that 
  $$
  \irep(e_1) \vee \irep(e_2) = 1,
  $$
  which is tantamount to $(\dagger)$.  Representations of $S$ obeying
this conditions will therefore be in one-to-one correspondence with
representations of the Cuntz algebra ${\cal O}_2$.  

The conclusion is
therefore that the missing link between the representation theory of
$S$ and the Cuntz algebras lies in the order structure of the
semillatice $E(S)$ in relation to Boolean algebras of projections on
Hilbert's space.

\section{\'Etale groupoids}
  \label GroupoidSection
  In this section we will review the basic facts about \'etale
groupoids which will be needed in the sequel.  We follow more or less
closely two of the most basic references in the subject, namely
\cite{\RenaultThesis} and \cite{\PatBook}.  We will moreover strive to
assume as few axioms and hypotheses as possible.

We assume the reader is familiar with the notion of groupoids (in the
purely algebraic sense) and in particular with its basic notations: a
groupoid is usually denoted by $\G$, its unit space by $\Gexp0$, and
the set of composable pairs by $\Gexp2$.  Finally the source and range
maps are denoted by $\s$ and $\r$, respectively.

Given our interests, we go straight to the definition of \'etale
groupoids without attempting to first define general locally compact
groupoids.
  We nevertheless begin by recalling from \scite{\RenaultThesis}{I.2.1} that
a \stress{topological groupoid} is a groupoid with a (not necessarily
Hausdorff) topology with respect to which both the multiplication and
the inversion are continuous.

\definition
  \label RenaultDefineRDiscrete
  \scite{\RenaultThesis}{I.2.8} \
  An \stress{\'etale} (sometimes also called \stress{$r$-discrete})
groupoid is a topological groupoid $\G$, whose unit space $\Gexp0$ is
locally compact and Hausdorff in the relative topology, and such that
the range map
  \ $
  \r: \G\to\Gexp0
  $ \
  is a local homeomorphism.
 
\cryout{From now on we will assume that we are given an \'etale
groupoid $\G$.}
  
It is well known that $\s(x) = \r(x\inv)$, for every $x$ in $\G$, and
hence $\s$ is a local homeomorphism as well.  Like any local
homeomorphism, $\s$ and $\r$ are open maps.

\state Proposition
  \label GZeroOpen
  $\Gexp0$ is an open subset of $\G$.
    
\proof
  Let $x_0\in\Gexp0$.  By assumption there is an open subset $A$ of
$\G$ containing $x_0$, and an open subset $B$ of $\Gexp0$ containing
$\r(x_0)$, such that $\r(A)=B$, and $\r|_{A}$ is a homeomorphism onto
$B$.  Set $B'=A\cap B$, and notice that 
  $$
  x_0=\r(x_0)\in A\cap B = B'.
  $$
  Given that $A$ is open in $\G$ we see that $B'$ is open in $B$,
hence $A':=\r\inv(B')\cap A$ is open in $A$, and moreover $\r$ is a
homeomorphism from $A'$ to $B'$.

We next claim that $B'\c A'$.  In order to prove it let $x\in
B'$.  So $x\in B\c\Gexp0$, and hence $x=\r(x)$.  This implies that
$x\in\r\inv(B')$, and we already know that $x\in A$, so $x\in
r\inv(B')\cap A=A'$.

We conclude that $\r$ is a bijective map from $A'$ to $B'$, which
restricts to a surjective map (namely the identity) on the subset
$B'\c A'$.  This implies that $B'=A'$, and since $A'$ is open in $\G$,
so is $B'$.
  The conclusion then follows from the fact that
  $$
  x_0\in B'\c\Gexp0.
  \proofend
  $$

\definition
  An open  subset $U\c\G$ is said to be a 
  \stress{slice}\fn{Slices are sometimes referred to as \stress{open
$\G$-sets.   Some authors use the notation $\G^{op}$ for the set of all slices.}}
  if the restrictions of $\s$ and $\r$ to $U$ are injective.

Since $\s$ and $\r$ are local homeomorphisms, for every slice $U$ one
has that $\s$ and $\r$ are homeomorphisms from $U$ onto 
  $\s(U)$ and $\r(U)$, respectively.
  For the same reason $\s(U)$ and $\r(U)$ are open subsets of $\Gexp0$, and hence
also open in $\G$.  It is obvious that every open subset of a slice is
also a slice.

  \state Proposition
  \label GzeroIsSlice
  $\Gexp0$ is a slice.

\proof By \lcite{\GZeroOpen} $\Gexp0$ is open in $\G$.  Since $\r$ and
$\s$ coincide with the identity on $\Gexp0$, they are injective.
  \proofend  

We next present a crucial property of slices, sometimes used as the
definition of \'etale groupoids \scite{\PatBook}{Definition 2.2.3}:

\state Proposition
  \label HasSlices
  The collection of all slices forms a basis for the topology of $\G$.

\proof
  Let $V$ be an open subset of $\G$ and let $x_0\in V$.  We must prove
that there exists a slice $U$ such that $x_0\in U\c V$.

  Since $\r$ is a local homeomorphism there is an open subset $A_1$ of
$\G$ containing $x_0$, and an open subset $B_1$ of $\Gexp0$ containing
$\r(x_0)$, such that $\r(A_1)=B_1$, and $\r|_{A_1}$ is a homeomorphism
onto $B_1$.  Since $\s$ is also a local homeomorphism, we may choose
an open subset $A_2$ of $\G$ containing $x_0$, and an open subset
$B_2$ of $\Gexp0$ containing $\s(x_0)$, such that $\s(A_2)=B_2$, and
$\s|_{A_2}$ is a homeomorphism onto $B_2$.  Therefore $U:= A_1\cap
A_2\cap V$ is a slice containing $x_0$, and contained in $V$.
  \proofend

\medskip If $U$ is a slice then $\r(U)$ is an open subset of the
locally compact Hausdorff space $\Gexp0$, and hence $\r(U)$ also
possesses these properties.  Since $U$ is homeomorphic to $\r(U)$ we
have:

\state Proposition
  Every slice is a locally compact Hausdorff space in the relative
topology.

If $V$ is a subset of $\G$ which is open and Hausdorff, observe that
the set of all intersections $V\cap U$, where $U$ is a slice, forms a
basis for the relative topology of $V$ by \lcite{\HasSlices}.  Notice
that $V\cap U$ is locally compact because it is an open subset of the
locally compact space $U$.  From this it is easy to see that $V$
itself is locally compact.  This proves:

\state Proposition
  \label OpenLocCpct
  Every open Hausdorff subset of $\G$ is locally compact.

 The following result is proved in Proposition
 (2.2.4) of \cite{\PatBook}, and although we are
 not assuming the exact same set of hypotheses, the
 proof given there works under our conditions:
 
 \state Proposition
   \label ProdSlices
   If\/ $U$ and $V$ are slices then
   \izitem
   \zitem $U\inv = \{u\inv : u\in U\}$ is a slice,
 and
   \zitem $UV = \{uv: u\in U,\ v\in V,\
 (u,v)\in\Gexp2\}$ is a (possibly empty) slice.

The theory of continuous functions on locally compact Hausdorff spaces
has many rich features which one wishes to retain in the study of non
necessarily Hausdorff groupoids.  The definition of $C_c(\G)$, first
used in \cite{\Connes}, and also given in \cite{\PatBook}, takes
advantage of the abundance of Hausdorff subspaces of $\G$.

\definition
  \label DefineCc
  We shall denote by $C^0_c(\G)$ the set of all complex valued
functions $f$ on $\G$ for which there exists a subset $V\c\G$,
such that
  \izitem 
  \zitem $V$ is open and Hausdorff in the relative topology,
  \zitem $f$ vanishes outside $V,$ and
  \zitem the restriction of $f$ to $V$ is continuous and compactly
supported\fn{That is $f|_V\in C_c(V)$, where the latter has the usual
meaning.}.  
  \medskip\noindent We finally define $C_c(\G)$ as the linear span of
$C^0_c(\G)$ within the space of all complex valued functions on $\G$.

We would like to stress that functions in $C_c(\G)$ might not be continuous
relative to the global topology of $\G$.

Suppose that $V$ is an open Hausdorff subset of $\G$ and let $f\in
C_c(V)$.  Considering $f$ as a function on $\G$ by extending it to be
zero outside $V$, it is immediate that $f\in C^0_c(\G)$, and hence
also $f\in C_c(\G)$.  This said we will henceforth view $C_c(V)$ as a
subset of $C_c(\G)$.

\state Proposition
  \label LinComb
  Let $\hbox{\rs C}$  be a covering of $\G$ consisting of slices.  Then
$C_c(\G)$ is linearly spanned by the collection of all subspaces of
the form $C_c(U)$, where $U\in\hbox{\rs C}$.

\proof Given $f\in C^0_c(\G)$ pick $V$ as in \lcite{\DefineCc}.  Observe
that $V$ is locally compact Hausdorff by \lcite{\OpenLocCpct}, and that
$\{U\cap V: U\in\hbox{\rs C}\, \}$ is a covering for $V$.  We may then
use a standard partition of unit argument to prove that $f$ may be
written as a finite sum of functions $f_i\in C_c(V\cap U_i)\c C_c(U_i)$,
where each $U_i$ is a slice in $\hbox{\rs C}$.  This concludes the
proof.
  \proofend

We are now about to introduce the operations that will eventually lead
to the C*-algebra of $\G$.  Normally this is done by first introducing a
Haar system on $\G$.  In \'etale groupoids a Haar system is just a
collection of counting measures, so the whole issue becomes a lot
simpler.  So much so that we can get away without even mentioning Haar
systems.

\state Proposition
  \label IntroduceOperations
  Given $f,g\in C_c(\G)$ define, for every $x\in\G$,
  $$
  (f\star g)(x) = \sum_{
  \buildrel {\scriptstyle (y,z)\in\Gexp2}
    \over {x=yz \pilar {6pt}}
  }f(y)g(z)
  ,\qquad\hbox{and}\qquad
  f^*(x) = \overline{f(x\inv)}.
  $$
  \vskip -20pt \noindent
  Then 
  \izitem 
  \zitem $f\star g\,$ and $\,f^*$ are well defined complex functions
on $\G$ belonging to $C_c(\G)$,
  \zitem if $f\in C_c(U)$ and $g\in C_c(V)$, where $U$ and $V$ are
slices, then $f\star g\in C_c(UV)$.
  \zitem if $f\in C_c(U)$, where $U$ is a slice, then $f^*\in
C_c(U\inv)$.

\proof
  The parts of the statement concerning $f^*$ are trivial, so we leave
them as exercises.
  We begin by addressing the finiteness of the sum above.  For this we
use \lcite{\LinComb} to write $f=\sum_{i=1}^n f_i$, where each $f_i\in
C_c(U_i)$, and $U_i$ is a slice.

If $x=yz$, and $f(y)g(z)$ is nonzero, then $f_i(y)$ is nonzero for some
$i=1,\ldots,n$, and hence $y\in U_i$.  Observing that $\r(y)=\r(x)$, and
that there exists at most one $y\in U_i$ with that property, we see that
there exists at most $n$ pairs $(y,z)$ such that $yz=x$, and $f(y)\neq
0$.  This proves that the above sum is finite, and hence that $f\star g$
is a well defined complex valued function on $\G$.

Since ``$\star$" is clearly a bilinear operation, in order to prove
that $f\star g\in C_c(\G)$, one may again use \lcite{\LinComb} in
order to assume that $f\in C_c(U)$ and $g\in C_c(V)$, where $U$ and
$V$ are slices.  That is, it suffices to prove (ii), which we do next.

So let us be given $f$ and $g$ as in (ii).  If $(f\star g)(x) \neq 0$,
then there exists at least one pair $(y,z)\in\Gexp2$, such that $x=yz$,
$y\in U$, and $z\in V$, but since $\r(y)=\r(x)$, and $\s(z)=\s(x)$, we
necessarily have that
  $$
  y= \r_U\inv\r(x)
  \and
  z = \s_V\inv\s(x),
  $$
  where we are denoting by $\r_U$ the restriction of $\r$ to $U$, and by
$\s_V$ the restriction of $\s$ to $V$.  It is then easy to see that
  $$
  (f\star g)(x) = 
  \left\{\matrix{
  f\big(\r_U\inv\r(x)\big) \ g\big(\s_V\inv\s(x)\big), & \hbox{ if }
x\in UV,\cr
  \pilar {17pt} \hfill 0\hfill , & \hbox{ otherwise.}
  }
  \right.
  \subeqmark FormulaForProd
  $$
  In addition the above formula for $f\star g$ proves that it is
continuous on $UV$, so we must only show that $f\star g$ is compactly
supported on $UV$.  If $A\c U$ and $B\c V$ are the compact supports of
$f$ and $g$ in $U$ and $V$, respectively, we claim that
  $AB$ is compact.
  In fact, since $\Gexp0$ is Hausdorff, we have that 
  $$
  \Gexp2 = \{(x,y)\in\G\times\G: \s(x)=\r(y)\} 
  $$
  is closed in $\G\times \G$.  
  So \ $(A\times B) \cap \Gexp2$ \ is closed in $A\times B$, and hence
compact.  Since $AB$ is the image of $(A\times B) \cap \Gexp2$ under
the continuous multiplication operator, we conclude that $AB$ is
compact, as claimed.
  Observing that $f\star g$ vanishes outside $AB$, we deduce that
$f\star g\in C_c(UV)$.
  \proofend

It is now routine to show that $C_c(\G)$ is an associative complex *-algebra
with the operations defined above. 

We have already commented on the fact that $C_c(V)\c C_c(\G)$, for
every open Hausdorff subset $V$ of $\G$, and hence $C_c(\Gexp0)\c
C_c(\G)$.  A quick glance at the definitions of the operations will
convince the reader that $C_c(\Gexp0)$ is also a *-subalgebra of
$C_c(\G)$, the induced multiplication and adjoint operations
corresponding to the usual pointwise operations on $C_c(\Gexp0)$.

\state  Proposition
  \label BundleUnitFiberProp
  Let $U$ be a slice and let $f\in C_c(U)$. Then $f\star f^*$ lies in
$C_c(\Gexp0)$.

\proof
  Given that $U$ is a slice, we have that $UU\inv=\r(U)$.  Moreover,
by \lcite{\IntroduceOperations.iii} we have that $f^*\in C_c(U\inv)$,
and hence by \lcite{\IntroduceOperations.ii},
  $$
  f\star f^*\in C_c(UU\inv)\c C_c\big(\r(U)\big) \c C_c(\Gexp0).
  \proofend
  $$ 

We would now like to discuss representations of $C_c(\G)$.  So let $H$
be a Hilbert space and let
  $$
  \pi:C_c(\G)\to B(H)
  $$
  be a *-representation.  Obviously the restriction of $\pi$ to
$C_c(\Gexp0)$ is a *-representation of the latter.  Since
$C_c(\Gexp0)$ is the union of 
  C*-algebras\fn{Namely the subalgebras of $C_c(\Gexp0)$ formed by all
continuous functions that vanish outside a fixed compact subset
$K\c\Gexp0$.}, $\pi$ is necessarily contractive with respect to the
norm $\|\cdot\|_\infty$ defined by
  $$
  \|f\|_\infty=\sup_{x\in\Gexp0}|f(x)|
  \for f\in C_c(\Gexp0).
  $$
  Therefore, if $U$ is a slice and $f\in C_c(U)$, we have
  $$
  \|\pi(f)\|^2 = 
  \|\pi(f)\pi(f)^*\| =
  \|\pi(f\star f^*)\| \leq
  \|f\star f^*\|_\infty,
  \eqmark PreContractiveOnSlice
  $$
  because $f\star f^*\in C_c(\Gexp0)$, by
\lcite{\BundleUnitFiberProp}.  Notice that for every $x\in \Gexp0$ we
have
  $$
  (f\star f^*)(x) = \sum_{x=yz} f(y) \overline{f(z\inv)},
  $$
  where any nonzero summand must correspond to a pair $(y,z)$ such that
$\s(y) = \r(z) = \s(z\inv)$, with both $y,z\inv\in U$.  Since $U$ is a
slice this implies that $y=z\inv$, but since $\r(y)=\r(x)$, we have that
$y=\r_U\inv\big(\r(x)\big)$, so, provided $(f\star f^*)(x)$ is nonzero one
has that
  $
  (f\star f^*)(x) = |f(\r_U\inv\big(\r(x)\big)|^2,
  $
  therefore
  $$
  \|f\star f^*\|_\infty =
  \sup_{x\in\Gexp0}|(f\star f^*)(x)| =
  \sup_{u\in U} |f(u)|^2 = \|f\|_\infty^2,
  $$
  where we are also denoting by $\|\cdot\|_\infty$ the sup norm on
$C_c(U)$.
  Combining this with \lcite{\PreContractiveOnSlice} we have proven:

\state Proposition 
  \label ContractToInftyNorm
  If $\pi$ is any *-representation of $C_c(\G)$ on a Hilbert space $H$
then for every slice $U$ and for every $f\in C_c(U)$ one has that
  $
  \|\pi(f)\| \leq \|f\|_\infty.
  $

\medskip By \lcite{\LinComb} any $f\in C_c(\G)$ may be written a
finite fum
  $
  f = \sum_{i=1}^n f_i,
  $
  where $f_i\in C_c(U_i)$, and $U_i$ is a slice.  So for every
representation $\pi$ of $C_c(\G)$ we have
  $$
  \|\pi(f)\| \leq
  \sum_{i=1}^n \|\pi(f_i)\| \leq
  \sum_{i=1}^n \|f_i\|_\infty,
  \eqmark ContractiveOnAll
  $$
  by \lcite{\ContractToInftyNorm}.  Regardless of its exact
significance,%
 \fn{It is related to the so called $I$-norm of $f$
\scite{\RenaultThesis}{II.1.3}.}
  the right-hand side of \lcite{\ContractiveOnAll} depends only on $f$
and not on $\pi$.  This means that
  $$
  \tn f := \sup_\pi\|\pi(f)\| <\infty,
  \eqmark DefineTripleNorm
  $$
  for all $f\in C_c(\G)$.  It is then easy to see that $\tn\cdot$ is a
C*-seminorm on $C_c(\G)$ and hence its Hausdorff completion is a
C*-algebra.

\definition
  The C*-algebra of $\G$, denoted $C^*(\G)$, is defined to be the
completion of $C_c(\G)$ under the norm $\tn{\cdot}$ defined above.  We
will moreover denote by
  $$
  \iota: C_c(\G)\to C^*(\G)
  \subeqmark NaturalMap
  $$
  the natural inclusion given by the completion process, which is
injective by \scite{\RenaultThesis}{4.2.i}.

Let us now study approximate units in $C^*(\G)$.  For this recall that
$C_c(\Gexp0)$ is a subalgebra of $C_c(\G)$.

\state Proposition 
  \label ResultOnApproxUnits
  Let $\{u_i\}_{i\in I}$ be a bounded selfadjoint approximate unit
for $C_c(\Gexp0)$ relative to the norm $\|\cdot\|_\infty$.  Then
$\{\iota(u_i)\}_{i\in I}$ is an approximate unit for $C^*(\G)$.

  \proof
  It is obviously enough to prove that $\{u_i\}_{i\in I}$ is a bounded
approximate unit for $C_c(\G)$ relative to $\tn\cdot$.  In view of
\lcite{\LinComb} it in fact suffices to verify that for every slice
$U$ and for every $f\in C_c(U)$ one has that
  $$
  \lim_i \tn{f\star u_i - f} = 0.
  $$
  By \lcite{\GzeroIsSlice} we have that $\Gexp0$ is a slice, and it is
easy to see that $U\Gexp0 = U$, so by \lcite{\IntroduceOperations} we
have that $f\star u_i\in C_c(U)$.  Moreover, for every $x\in U$ one has
that
  $$
  (f\star u_i)(x)=f(x)u_i(\s(x)).
  $$
  By \lcite{\ContractToInftyNorm} we conclude that
  $$
  \tn{f\star u_i - f} \leq
  \sup_{x\in U} \big|f(x)u_i(\s(x)) - f(x)\big|,
  $$
  which converges to zero because $u_i$ converges uniformly to 1 on
every compact subset of $\Gexp0$, such as the image under $\s$ of the
compact support of $f$.
  \proofend


\section{Inverse semigroup actions}
  Recall that a semigroup $\S$ is said to be an \stress{inverse
semigroup} if for every $s\in\S$, there exists a unique $s^*\in\S$
such that
  $$
  ss^*s=s
  \and
  s^*ss^*=s^*.
  \eqmark ISGAxioms
  $$
  It is well known that the correspondence $s\mapsto s^*$ is then an
involutive anti-homomorphism.
  One usually denotes by $E(\S)$ the set of all idempotent elements of
$\S$, such as $s^*s$, for every $s\in\S$.  For a thorough treatment of
this subject the reader is referred to \cite{\Lawson}, and
\cite{\PatBook}.

  We next recall the definition of one of the most important examples
of inverse semigroups:

\definition
  \label DefineIX
  If $X$ is any set we denote by $\I(X)$ the inverse semigroup formed by
all bijections between subsets of $X$, under the operation given by
composition of functions in the largest domain in which the composition
may be defined.

The following is a crucial concept to be studied throughout the
remaining of this work.

\definition
  \label DefineAction
  Let $\S$ be an inverse semigroup and let $X$ be a locally compact
Hausdorff topological space.
  An \stress{action} of $\S$ on $X$ is a semigroup homomorphism
  $$
  \act: \S \to \I(X) 
  $$
  such that,
  \izitem
  \zitem for every $s\in\S$ one has that $\act_s$ is continuous and
its domain is open in $X$,
  \zitem the union of the domains of all the $\act_s$ coincides with
$X$.

\cryout{Fix for the duration of this section an action $\act$ of $\S$
on $X$.}

Observe that if $s\in\S$ then from \lcite{\ISGAxioms} we get
  $$
  \act_s\act_{s^*}\act_s=\act_s \and
  \act_{s^*}\act_s\act_{s^*}=\act_{s^*},
  $$
  which implies that $\act_{s^*}=\act_s\inv$.

Notice that the range of each $\act_s$ coincides with the
domain of $\act_s\inv=\act_{s^*}$, and hence it is open as well.
This also says that $\act_s\inv$ is continuous, so $\act_s$ is
necessarily a homeomorphism onto its range.

In the absence of the property expressed in the last sentence of the
above definition one may replace $X$ by the open subspace $X_0$ formed
by the union of the domains of all the $\act_s$.  It is then apparent
that $\act$ gives an action of $\S$ on $X_0$ with all of the desired
properties.  In other words, the restriction imposed by that requirement
is not so severe.

It is well known that if $e$ is an idempotent, that is, if $e^2=e$, then
$\act_e$ is the identity map on its domain.

\def\RemarkOnDe{Some authors adopt the notation $D_s$, even if $s$ is
not idempotent, to mean the \stress{range} of $\act_s$, which
therefore coincides with our $D_{ss^*}$.  We shall however not do so in
order to avoid introducing an unnecessary convention: one could
alternatively choose to use $D_s$ to denote the \stress{domain} of
$\act_s$.  Once one is accustomed to the idea that the source and
range projections of a partial isometry $u$ are $\q u$ and $\p u$,
respectively, the notations $D_{\q s}$ and $D_{\p s}$ require no
convention to convey the idea of domain and range.}

\sysstate{Notation}{\rm}{For every idempotent $e\in E(\S)$ we will
denote\fn{\RemarkOnDe} by $D_e$ the domain (and range) of $\act_e$.}

  \halffootno =\footno
  \divide \halffootno by 2
  \edef\FootnoteOnDe{\number\halffootno}

  It is easy to see that $\act_s$ and $\act_{\q s}$ share
domains, and hence the domain of $\act_s$ is $D_{s^*s}$.  Likewise
the range of $\act_s$ is given by $D_{ss^*}$.  Thus $\act_s$ is a
homeomorphism between the open sets
  $$
  \act_s: D_{s^*s} \to D_{ss^*}.
  $$

If $e$ and $f$ are idempotents it is easy to conclude from the identity
$\act_e\act_f = \act_{ef}$, that $D_e\cap D_f = D_{ef}$.  The next
result appears in \scite{\Sieben}{4.2}.

\state Proposition
  \label MovingDe
  For each $s\in \S$ and $e\in E(\S)$
one has that
  $$
  \act_s(D_e\cap D_{\q s}) = D_{ses^*}.
  $$

  \proof By the observation above we have 
  $$
  \act_s(D_e\cap D_{\q s}) =   \act_s(D_{e\q s}),
  $$
  which coincides with the range of 
  $
  \act_s\act_{e\q s} = 
  \act_{se\q s}.
  $
  The conclusion then follows from the
following calculation:
  $$
  se\q s (se\q s)^* =
  se\q s \ \q s es^* = ses^*.
  \proofend
  $$

Our next short term goal is to construct a \stress{groupoid of germs}
from $\act$.  However, given the examples of inverse semigroups that
we have in mind, we would rather not assume that $\act$ is given in
terms of a \stress{localization}, as in \scite{\PatBook}{Theorem
3.3.2}.  Nor do we want to assume that $\S$ is \stress{additive}, as in
\scite{\PatBook}{Corollary 3.3.2}.  We also want to avoid using the
condition of \stress{fullness} \scite{\QuigSieb}{5.2}, which is used to
prove a result \scite{\QuigSieb}{8.1} similar to what we are looking for
in the Hausdorff case.

\definition
  \label DefineGermEquiv
  \scite{\PatBook}{page 140} \
  We will denote by $\Pregerm$ the subset of $\S\times X$ given by
  $$
  \Pregerm = \{(s,x)\in \S\times X: x\in D_{s^*s}\},
  $$
  and for every $(s,x)$ and $(t,y)$ in $\Pregerm$ we will say that 
  $
  (s,x)\sim(t,y),
  $
  if $x=y$, and there exists an idempotent $e$ in $E(\S)$ such that
$x\in D_e$, and  $se=te$.  The 
equivalence class of $(s,x)$ will be called the \stress{germ of $s$ at $x$},
and will be denoted by
$[s,x]$.

Given $(s,x)$ and $(t,y)$ in $\Pregerm$ such that $(s,x)\sim(t,y)$, and
letting $e$ be the idempotent mentioned in \lcite{\DefineGermEquiv},
observe that
  $$
  x\in D_e\cap D_{\q s}\cap D_{\q t} = D_{e\q s\q t}.
  $$
  If we set $e_0= e\q s\q t$, it then follows that 
  $
  se_0 = 
  te_0.
  $
  So, upon replacing $e$ by $e_0$, we may always assume that the
idempotent $e$ in \lcite{\DefineGermEquiv} satisfies $e\leq \q s,\q
t$.

\state Proposition
  Given $(s,x)$ and $(t,y)$ in $\Pregerm$ such that $x=\act_t(y)$, one
has that
  \izitem
  \zitem $(st,y)\in\Pregerm$, and
  \zitem the germ $[st,y]$ depends only on the germs $[s,x]$ and
$[t,y]$.

  \proof Initially observe that
  $$
  y = \act_{t^*}(x) \in 
  \act_{t^*}(D_{s^*s}\cap D_{tt^*}) \={(\MovingDe)} 
  D_{t^*s^*st} =   D_{(st)^*st},
  $$
  so $(st,y)$ indeed belongs to $\Pregerm$.  Next let $(s',x)$ and
$(t',y)$ be elements of $\Pregerm$ such that $(s',x)\sim(s,x)$ and
$(t',y)\sim(t,y)$.  Therefore there are idempotents $e$ and $f$ such
that $x\in D_e$, $y\in D_f$, $se=s'e$, and $tf=t'f$.  We then have that
  $$
  \act_{t'}(y) =   
  \act_{t'}\big(\act_{f}(y)\big) = 
  \act_{t'f}(y) =   
  \act_{tf}(y) =   
  \act_{t}\big(\act_{f}(y)\big) = 
  \act_{t}(y) = x.
  $$
  In other words, the fact that $\act_t(y)=x$ does not depend on
representatives.  By (i) it then follows that $(s't',y)\in\Pregerm$ and
we will be finished once we prove that
  $(st,y)\sim (s't',y)$.  For this let $d$ be the idempotent given by
$d=ft^*et$, and we claim that $y\in D_d$.  To see this notice that
since $x\in D_{e}\cap
D_{tt^*}$, it follows that
  $$
  y=\act_{t^*}(x) \in
  \act_{t^*}(D_{e}\cap D_{tt^*}) = 
  D_{t^*et},
  $$
  and since $y\in D_f$ by assumption, we deduce that
  $$
  y\in
  D_f \cap D_{t^*et}  = 
  D_{ft^*et}  = 
  D_d.
  $$
  This proves our claim.  In addition we have 
  $$
  s't'd = 
  s't' ft^*et =
  s't ft^*et =
  s'et ft^*t =
  set ft^*t =
  st ft^*et =
  std,
  $$
  proving that $(st,y)\sim(s't',y)$, as required.
  \proofend

  Let 
  $$
  \G =\Pregerm/{\sim}
  $$
  be the set of all germs, and put
  $$
  \Gexp2 = 
  \left\{\Big([s,x],[t,y]\Big)\in\G\times\G: x = \act_t(y)\right\}.
  \eqmark DefineGTwoGpd
  $$
  For $\Big([s,x],[t,y]\Big)\in\Gexp2$ define 
  $$
  [s,x] \cdot [t,y] = [st,y],
  \eqmark DefineProductGpd
  $$
  and
  $$
  [s,x]\inv = [s^*,\act_s(x)].
  \eqmark DefineInverseGpd
  $$

We leave it for the reader to prove:

  \state Proposition 
  \label GIsGroupoid
  $\G$ is a groupoid with the operations defined above, and the unit
space $\Gexp 0$ of\/ $\G$ naturally identifies with $X$ under the
correspondence
  $$
  [e,x] \in \Gexp 0 \mapsto x \in X, 
  $$
  where $e$ is any idempotent such that $x\in D_e$.

Although we are not providing a proof of the result above, we observe
that the last part of the statement depends upon the assumption made in
the last sentence of Definition \lcite{\DefineAction}.

The \stress{source}\fn{Given the several uses of the letter ``s" in
the setting of semigroups, we have decided to allow the idea of
``domain" to determine the letter to denote the source map,
a convention that is not rare in the literature.}  map of $\G$ is
clearly given for every $[t,x]\in\G$ by
  $$
  \s[t,x] =
  [t,x]\inv[t,x] = 
  [t^*,\act_t(x)]\ [t,x] = 
  [t^*t,x].
  $$
  Enforcing the identification referred to in 
\lcite{\GIsGroupoid} we will write
  $$
  \s[t,x] =
  x.
  $$  
  With respect to the \stress{range} map, a similar reasoning gives 
  $$
  \r[t,x] = \act_t(x).
  $$
  
We would now like to give $\G$ a topology.  For this, given any $s\in
S$, and any open subset $U\c D_{\q s}$, let
  $$
  \O(s,U) = \{[s,x]\in\G: x\in U\}.
  \eqmark DefineOsuBasicSlice
  $$

\state Proposition
  Let $s$ and $t$ be elements of $\S$ and let $U$ and $V$ be open sets
with $U\c D_{\q s}$, and $V\c D_{\q t}$.
  If \ $[r,z]\in \O(s,U)\cap \O(t,V)$ then there exists an idempotent
$e$ and an open set $W\c D_{(re)^*re}$ such that
  $$
  [r,z]\in \O(re,W) \c \O(s,U)\cap \O(t,V).
  $$
  
\proof  By assumption $[r,z] = [s,x] = [t,y]$, for some $x\in U$ and
$y\in V$.  But this implies that $z=x=y$, so $z\in U\cap V$.  In
  addition there are idempotents $e$ and $f$ such that $z\in D_e$,
$z\in D_f$, $re=se$, and $rf=tf$.  Replacing $e$ and $f$ by $ef$, we
  may assume without loss of generality that $e=f$,  hence 
$re=se=te$.  
  Set $W  = U \cap V\cap D_{(re)^*re}$.
  Since 
  $
  z\in   D_{r^*r} \cap D_{e} =
  D_{r^*re} = 
  D_{(re)^*re},
  $ 
  we see that $z\in W$, and hence
  $$
  [r,z] =   [re,z]\in \O(re,W).
  $$
  In order   to prove that $\O(re,W) \c \O(s,U)\cap \O(t,V)$, let 
$[re,x]$ be a generic element of $\O(re,W)$, so that $x\in W$.
Noticing that  $x\in U$, and  that 
  $$
  [re,x] = [se,x] = [s,x],
  $$
  we see that   $[re,x] \in \O(s,U)$, and a similar reasoning gives 
$[re,x] \in \O(t,V)$.
  \proofend

By the result above we see that the collection of all
$\O(s,U)$ forms the basis of a topology on $\G$.
  From now on $\G$ will be considered to be equipped with this
topology, and hence $\G$ is a topological space.

\state Proposition
  With the above topology $\G$ is a topological groupoid.

\proof
  Our task is to prove that the multiplication and inversion
operations on $\G$ are continuous.  For this let $[s,x]$ and $[t,y]$
be elements of $\G$ such that $\big([s,x],[t,y]\big)\in\Gexp2$.
Moreover suppose that the product of these elements lie in a given
open set $W\c\G$.  Therefore, there exists some $r\in\S$ and an open
set $V\c D_{\q r}$, such that
  $$
  [s,x][t,y] = [st,y] \in \O(r,V)\c W.
  $$
  This implies that $y\in V$ and that there exists some idempotent $e$
such that $y\in D_e$, and $ste=re$.

Setting $U=V\cap D_e\cap D_{\q t}$, we will prove that the product of
any pair of elements 
  $$
  \big([s,x'],[t,y']\big)\ \in\ 
  \big(\O(s,D_{\q s}) \times \O(t,U)\big) \ \cap \ \Gexp2
  \subeqmark RelTopOnGTwo
  $$
  lies in W.  The product referred to is clearly given by $[st,y']$, 
  and since $y'\in U\c D_e$, we have
  $$
  [st,y'] = [r,y'] \in \O(r,V)\c W.
  $$
  Observing that $x\in D_{\q s}$, and $y\in V\c U$, we see that the set
appearing in
  \lcite{\RelTopOnGTwo}
  is a neighborhood of $\big([s,x],[t,y]\big)$ in the relative topology
of $\Gexp2$.  This proves that multiplication is continuous.

   With respect to inversion let $s\in\S$ and let $U\c D_{\q s}$
be an open set.  From the definition of the inversion in
\lcite{\DefineInverseGpd} it is clear that 
  $$
  \O(s,U)^* = \O(s^*,\act_s(U)),
  $$
  from which the continuity of the inversion follows immediately.
  \proofend

We shall now begin to work towards proving that $\G$ is an \'etale
groupoid.

\state Proposition
  \label PhiHomeo
  Given $s\in\S$, let $U\c D_{\q s}$ be an open set.  Then the map
  $$
  \phi: x\in U \mapsto [s,x]\in\O(s,U)
  $$
  is a homeomorphism, where $\O(s,U)$ of course carries the topology
induced from $\G$.

  \proof By the definition of the equivalence relation in
\lcite{\DefineGermEquiv} it is obvious that $\phi$ is a bijective map.
 Let $V\c U$ be an open subset.  Then clearly
$\phi(V) = \O(s,V)$, which is open in $\G$, proving that $\phi$ is an
open mapping.  To prove that $\phi$ is continuous at any given $x\in
U$, let $W$ be a neighborhood of $\phi(x)$ in $\O(s,U)$, so there
exists some $t\in\S$ and an open set $V\c D_{t^*t}$, such that
  $$
  [s,x] = \phi(x)\in\O(t,V)\c W \c \O(s,U).
  $$
  Clearly this implies that $x\in V\c U\c D_{\q s}$. 
  In addition there exists some idempotent $e$ such that $x\in D_e$,
and $se=te$. 
  %
  For every $y\in D_e\cap V$ observe that 
  $$
  \phi(y) = 
  [s,y] = 
  [t,y] \in
  \O(t,V) \c W,
  $$
  which means that $\phi(D_e\cap V)\c W$.  Since $D_e\cap V$ is a
neighborhood of $x$ in $U$, we see that $\phi$ is continuous.
  \proofend

We have already seen that $\Gexp 0$, the unit space of $\G$,
corresponds to $X$.  The result above helps to complete that picture
by showing that the correspondence is topological:

\state Corollary 
  \label GZeroIsX
  The identification of\/ $\Gexp 0$ with $X$ given by
\lcite{\GIsGroupoid} is a homeomorphism.

\proof
  Given $[e,x]\in\Gexp 0$ we have that $D_e$ is an open subset of $X$
containing $x$ and $\O(e,D_e)$ is an open subset of $\Gexp 0$
containing $[e,x]$.  The result then follows from the fact that
  $$
  \phi: y\in D_e \mapsto [e,y]\in\O(e,D_e)
  $$
  is a homeomorphism by \lcite{\PhiHomeo}.
  \proofend

Having assumed that $X$ is locally compact and Hausdorff, it follows
from the above result that $\Gexp0$ shares these properties.  The
source map on every basic open set $\O(s,U)$ is a homeomorphism onto
$U$ because it is the inverse of the map $\phi$ of \lcite{\PhiHomeo}.
This implies that $\s$ is a local homeomorphism, and hence so is $\r$.
This implies that:

\state Proposition 
  \label TheMainGroupoid
  The groupoid $\G= \G(\act,\S,X)$ constructed above, henceforth 
  called the \stress{groupoid of germs} of the system $(\act,\S,X)$,
is an \'etale groupoid.

The following identifies important slices in $\G$.

\state Proposition 
  \label OSsIsSliceLegal
  For every $s\in\S$ and every open subset $U\c D_{\q
s}$, one has that $\O(s,U)$ is a slice.

\proof By definition of the topology on $\G$ we have that $\O(s,U)$ is
open in $\S$.  Recall that source map is given by
  $
  \s: [s,x]\mapsto x,
  $
  whose restriction to $\O(s,U)$ is the inverse of the map $\phi$ of
\lcite{\PhiHomeo}, so it is injective.  With respect to the restriction
of the range map $\r$ on $\O(s,U)$, notice that $\r = \act_s\circ\s$,
which is injective.
  \proofend

\section{Example: Action of the inverse semigroup of slices}
  The main goal of this section is to present an example of inverse
semigroup actions which is intrinsic to every \'etale groupoid.  We
therefore fix an \'etale groupoid $\G$ from now on.  Denote%
  \fn{As already observed some authors denote this set by $\G^{op}$.}
   by $\SG$ the set of all slices in $\G$.  It is well known
\scite{\PatBook}{Proposition 2.2.4} that $\SG$ is an inverse semigroup
under the operations
  $$
  UV = \{uv: u\in U,\ v\in V,\ (u,v)\in\Gexp2\}
  \and
  U^* = \{u\inv: u\in U\},
  $$
  for all slices $U$ and $V$ in $\SG$.
The idempotent semilattice of $\SG$ is easily seen to consist
precisely of the open subsets of $\Gexp0$.

Henceforth denoting  by 
  $$
  X := \Gexp0,
  $$
  we wish to define an action $\act$ of $\SG$ on $X$.  Given a slice
$U$ we have already mentioned that $\s(U)$ and $\r(U)$ are open
subsets of $X$, and moreover that the maps
  $$
  \s_U : U \to \s(U)
  \and
  \r_U : U \to \s(U),
  $$
  obtained by restricting $\s$ and $\r$, respectively, are
homeomorphisms.  Given $x\in\s(U)$ we let
  $$
  \act_U(x) = \r_U(\s_U\inv(x)).
  \eqmark DefineActU
  $$
  Clearly $\act_U$ is a homeomorphism from $\s(U)$ to $\r(U)$.
It is interesting to observe that $\act_U(x)=y$, if and only if there
exists some $u\in U$ such that $\s(u)=x$ and $\r(u)=y$.  Thus, if we
view $\act_U$ as a set of ordered pairs, according to the technical
definition of functions, we have
  $$ 
  \act_U = \big\{\big(\s(u),\r(u)\big) : u\in U\big\}.
  \eqmark TechDefActU
  $$
  We would like to show that $\act_U\act_V = \act_{UV}$, for all
$U,V\in\SG$.  Assuming that $\act_U\act_V(x) =z$, or equivalently that
$(x,z)\in\act_U\act_V$, there exists $y\in X$, such that
$(y,z)\in\act_U$ and $(x,y)\in\act_V$, so we may pick $u\in U$ and
$v\in V$ such that $\s(v)=x$, $\r(v)=y=\s(u)$, and $\r(u)=z$.

  \bigskip
  \medskip\beginmypicture
  \setcoordinatesystem units <0.0017truecm, 0.001truecm> point at 0 0
  \setquadratic
  \put {$\bullet$} at 0000 0000 \put {$x$} at 0000 -300 
  \put {$\bullet$} at 1000 0000 \put {$y$} at 1000 -300 
  \put {$\bullet$} at 2000 0000 \put {$z$} at 2000 -300 
  \plot 0000 0000 0500 0400  1000 0000 / \put {$v$} at 0500 700 
  \plot 1000 0000 1500 0400  2000 0000 / \put {$u$} at 1500 700 
  \arrow <0.15cm> [0.25,0.75] from 0480 0400 to 0540 0400
  \arrow <0.15cm> [0.25,0.75] from 1480 0400 to 1540 0400
  \endmypicture

  \bigskip \noindent 
  Therefore we have that $uv\in UV$, and since
  $$
  (x,z)= \big(\s(v),\r(u)\big) = 
  \big(\s(uv),\r(uv)\big) \in \act_{UV},
  $$
  we see that $\act_{UV}(x)=z$.  
  Conversely, if we are given that $\act_{UV}(x)=z$, there exists some
$w\in UV$ such that $\s(w)=x$, and $\r(w)=z$.  Writing $w=uv$, with
$u\in U$ and $v\in V$, set
  $y=\r(v)=\s(u)$.  Then
  $$
  (x,y)=
  \big(\s(w),\r(v)\big) = 
  \big(\s(v),\r(v)\big) \in\act_V,
  $$
  and similarly 
  $(y,z)=\big(\s(u),\r(u)\big)\in\act_U$, and we see that $(x,z)\in
\act_U\act_V$, thus proving  that $\act_U\act_V = \act_{UV}$.  

The last condition to be checked in order to prove that $\act$ is an
action is \lcite{\DefineAction.ii}, but this is obvious because
$X=\Gexp0$ is a slice by \lcite{\GzeroIsSlice}, and $\act_{X}$ is
clearly the identity map defined on the whole of $X$.
  With this we have proven:

\state  Proposition
  The correspondence \ $U\mapsto\act_U$, \ defined by
\lcite{\DefineActU}, gives an action of $\SG$ on the unit space of\/
$\G$.

Given any 
  *-subsemigroup\fn{A subsemigroup of an inverse semigroup is said to
be a *-subsemigroup if it is closed under the * operation, in which
case it is clearly an inverse semigroup in itself.}
   $\S\c\SG$, one may restrict $\act$ to $\S$, thus obtaining a
semigroup homomorphism
  $$
  \act|_\S: \S \to \I(X)
  $$
  which is an action of $\S$ on $X$, provided \lcite{\DefineAction.ii}
may be verified.
  The next result gives sufficient conditions for the groupoid of
germs for such an action to be equal to $\G$.

\state  Proposition
  \label TwoGroupoids
  Let $\G$ be an \'etale groupoid and let $\S$ be a *-subsemigroup of
$\S(\G)$ such that
  \izitem
  \zitem $\G = \union_{U\in\S} U$, and
  \zitem for every $U,V\in\S$, and every $u\in U\cap V$, there exists
$W\in\S$, such that $u\in W\c U\cap V$.
  \medskip\noindent
  Then $\act|_\S$ is an action of $\S$ on $X=\Gexp0$,
  and the groupoid of germs for $\act|_\S$ is isomorphic to $\G$.

\proof
  Given $x\in X$, there exists some $U\in\S$ such that $x\in U$, by
(i), and so $(x,x) = \big(\s(x),\r(x)\big) \in\act_U$, and in
particular $x$ is in the domain of $U$.  This proves
\lcite{\DefineAction.ii} and hence $\act|_\S$ is indeed an action of
$\S$ on $X$.

  Let us temporarily denote the groupoid of germs for $\act|_\S$ by
$\H$.
  Observe that the domain of $\act_U$ is $\s(U)$, so $\H$  is given by
  $$
  \H = \big\{[U,x]: U\in\S,\ x\in\s(U)\big\}.
  $$
  Given a germ $[U,x]\in\H$ we therefore have that there exists a
unique $u_0\in U$ such that $\s(u_0)=x$, because $\s|_U$ is injective.

We claim that $u_0$ depends only on the germ $[U,x]$.  For this
suppose that $[U,x]=[V,x]$, for some $V\in\S$, which means that there
is an idempotent $E\in\S$ such that $x\in\s(E)$ and $UE=VE$.  As
observed earlier, $E$ is necessarily a subspace of $X$ and hence
$E=\s(E)$.
  Applying the definition of the product one
gets
  $$
  UE = \{u\in U: \s(u)\in E\},
  $$
  and since $\s(u_0) = x\in\s(E)=E$, we conclude that $u_0\in UE$.
Therefore also $u_0\in VE$, and in particular $u_0\in V$.  This is to
say that the unique element $v\in V$, with $\s(v)=x$, is $u_0$, so the
claim is proved.  We may then set
  $
  \phi([U,x]) = u,
  $
  thus obtaining a well defined map
  $$
  \phi:\H\to\G.
  $$
  Employing the homeomorphisms $\s_U =\s|_U : U \to \s(U)$, for every
slice $U$, one may concretely describe $\phi$ by
  $$
  \phi([U,x]) = \s_U\inv(x).
  \subeqmark ConcreteDefOfPhi
  $$
  Another interesting characterization of $\phi$ is 
  $$
  \phi([U,x]) = u \iff u\in U, \hbox{ and } \s(u)=x.
  \subeqmark CleverDefOfPhi
  $$
  for every $U\in\S$, and every $x\in\s(U)$.
  To see that $\phi$ is surjective let $u\in\G$.  We may invoke (i) to
find some $U\in\S$ such that $u\in U$, and hence $[U,\s(u)]$ is in
$\H$ and 
  $$
  \phi\big([U,\s(u)]\big)=u.
  \subeqmark PhiSurjectice
  $$

We will next prove that $\phi$ is injective, and for this we let
$[U_1,x_1]$ and $[U_2,x_2]$ be germs in $\H$ such that
  $$
  \phi([U_1,x_1])=\phi([U_2,x_2]).
  $$
  Denoting by $w$ the common value of the terms above we have by
\lcite{\CleverDefOfPhi} that 
  $$
  w\in U_i \and \s(w)=x_i \for i=1,2.
  $$
  In particular $w\in U_1\cap U_2$, so (ii) applies providing some
$W\in\S$ such that $w\in W\c U_1\cap U_2$.
  The fact that $W\c U_i$ may be described in terms of the semigroup
structure of $\SG$ by saying that $W=U_iW^*W$, (compare
\fcite{PartialOrderISG}{6.1}), which in particular implies that
  $$
  U_1W^*W= U_2W^*W.
  $$
  Moreover
  $$
  x_1 = x_2 = \s(w) \in \s(W) = \s(W^*W),
  $$
  thus proving that $[U_1,x_1]=[U_2,x_2]$.

  Let us now prove that $\phi$ is a homeomorphism.  For this pick a
germ $[U,x]\in\H$ and recall from \lcite{\OSsIsSliceLegal} that
  $$
  \O_U := 
  \O\big(U,\s(U)\big) =
  \big\{[U,y]: y\in \s(U)\big\}
  $$
  is a slice in $\H$, which clearly contains $[U,x]$.  The image of
  $\O_U$ under $\phi$ is obviously $U$, and the restriction of $\phi$
to $\O_U$ is certainly continuous on $\O_U$ by
\lcite{\ConcreteDefOfPhi}.  On the other hand, for each $u\in U$, we
have that
  $$
  \phi\inv(u)=[U,\s(u)],
  $$
  by \lcite{\PhiSurjectice}.  Since $\phi\inv$ sends $U$ into the
slice $\O_U$, in order to prove that $\phi\inv$ is continuous on $U$,
it is enough to prove that $\delta\circ\phi\inv$ is continuous, where
we are denoting by $\delta$ the source map for the groupoid $\H$.
That composition is clearly given by
  $$
  \delta\circ\phi\inv(u) = \delta\big([U,\s(u)]\big) = \s(u)
  \for u\in U,
  $$
  which is well known to be continuous.  

It remains to prove that $\phi$ is an isomorphism of groupoids.  For
this let $[U,x]$ and $[V,y]$ be germs is $\H$ and let $u=\phi([U,x])$,
and $v=\phi([V,y])$, so that $u\in U$, $\s(u)=x$, $v\in V$, and
$\s(v)=y$, according to \lcite{\CleverDefOfPhi}.
  It is useful to remark that 
  $$
  \big(y,\r(v)\big)= \big(\s(v),\r(v)\big) \in \act_V,
  $$
  by \lcite{\TechDefActU}, so $\act_V(y) = \r(v)$.
  By \lcite{\DefineGTwoGpd} we have that $\big([U,x][V,y]\big)\in
\sysGexp{\H}{2}$ if and only $x=\act_V(y)$, which is equivalent to
saying that $\s(u)=\r(v)$, or that 
  $$
  \Big(\phi\big([U,x]\big),\phi\big([V,y]\big)\Big) = (u,v)\in\Gexp2.
  $$
  This says that two elements in $\H$ may be multiplied if and only if
their images under $\phi$ in $\G$ may be multiplied.  In this case we
have by \lcite{\DefineProductGpd} that
  $$
  [U,x][V,y] = [UV,y].
  $$
  On the other hand, notice that $uv\in UV$, and that
$\s(uv)=\s(v)=y$, so
  $$
  \phi\big([UV,y]\big) = uv = \phi\big([U,x]\big) \phi\big([V,y]\big),
  $$
  thus proving that $\phi$ is a homomorphism of groupoids.
  \proofend

Conditions \lcite{\TwoGroupoids.i-ii} look very much like the
definition of a topological base for $\G$.  Therefore if $\S$ is a
\stress{full} *-subsemigroup of $\SG$, in the sense of
\scite{\QuigSieb}{5.2}, then $\S$ clearly satisfies
\lcite{\TwoGroupoids.i-ii}.  However the latter conditions are clearly
much weaker than to require that $\S$ be a base for the topology of
$\G$.  For example, if $\G$ is a groupoid consisting only of units,
that is, if $\G$ is a topological space, then $\G$ itself is a slice
and the singleton $\{\G\}$ is a *-subsemigroup of $\SG$ which is not
full, but satisfies \lcite{\TwoGroupoids.i-ii}.

  \section{The Hausdorff property for the groupoid of germs}
  Quoting Paterson \cite{\PatBook}, the theory of non Hausdorff
groupoids presented in section \lcite{\GroupoidSection}, and employed
throughout this paper, already has enough of the Hausdorff property to
allow for the efficient use of standard topological methods.  However
should a groupoid be Hausdorff in the true sense of the word it is
definitely good to be aware of it.

It is not easy to determine conditions on an inverse semigroup $\S$ to
ensure that the groupoid $\G(\act,\S,X)$ of \lcite{\TheMainGroupoid}
be Hausdorff for any action $\act$ of $\S$ on any space $X$,
especially because even groups may present difficulties.  However the
actions we are interested in have a special property which may be
exploited in order to obtain such a characterization.  In what follows
we would like to describe this result.

Recall e.g.~from \scite{\Lawson}{1.4.6} that an inverse semigroup $\S$
is naturally equipped with a partial order defined by
  $$
  s\leq t \iff s=ts^*s
  \for s,t\in\S.
  \eqmark PartialOrderISG
  $$

\state Proposition
  \label HausdorffGPG
  Suppose that $\S$ is an inverse semigroup which is a semilattice%
  \fn{Not to be confused with the semilattice of idempotents of $\S$,
this means  that for every $s,t\in\S$, there is a maximum among
the elements of $\S$ which are smaller than both $s$ and $t$.  
  Tradition suggests that this element be denoted by $s\wedge t$. 
  It is convenient to observe that if $e$ and $f$ are idempotents in
$\S$ then the product $ef$ coincides with $e\inf f$.  However if $s$
and $t$ are not idempotents then the product $st$ is not always the
same as $s\inf t$.  }
  with respect to its natural order.  Let $\act$ be an action of $\S$
on a locally compact Hausdorff space $X$, such that for each $s\in\S$,
the domain $D_{s^*s}$ of\/ $\act_s$ is closed (besides being open).
Then $\G(\act,\S,X)$ is Hausdorff.

\proof
  Let $[s,x]$ and $[t,y]$ be two distinct elements of $\G(\act,\S,X)$.
We need to find disjoint open subsets $U$ and $V$ of $\G(\act,\S,X)$,
such that $[s,x]\in U$, and $[t,y]\in V$.  If $x\neq y$ this is quite
easy:  separate $x$ and $y$ within $X$ using disjoint open sets $A,B\c
X$, and take $U=\O(s,A\cap D_{\q s})$ and $V=\O(t,B\cap D_{\q t})$.

Let us then treat the less immediate case in which $x=y$.  For this
let $u=s\wedge t$ and notice that
  $$
  s\q u = u = t\q u,
  $$
  and hence $x\notin D_{\q u}$, or else $[s,x]=[t,x]$, by
\lcite{\DefineGermEquiv}.  As we are assuming that $D_{\q u}$ is
closed we deduce that $V=X\setminus D_{\q u}$ is an open neighborhood
of $x$ in $X$.  Setting $W = V\cap D _{\q s}\cap D _{\q t}$, it is clear that 
  $$
  [s,x] \in \O(s,W)
  \and 
  [t,x] \in \O(t,W).
  $$
  It therefore suffices to prove that $\O(s,W)$ and $\O(t,W)$ are
disjoint sets.  Arguing by contradiction suppose that $[r,z]\in
\O(s,W)\cap\O(t,W)$.  It follows that $[r,z]=[s,z]=[t,z]$, and hence 
there are idempotents $e$ and $f$ such that $z$ lies in  $D_e$ and in
$D_f$, and moreover such that 
  $re=se$, and $rf =tf$.  By replacing $e$ and $f$ with $ef$, we may
assume that $e=f$, in which case $re = se = te$.
  Then 
  $$
  s(re)^*(re) = s e r^*re = r e r^*re = r r^*re e = re,
  $$
  so $re\leq s$, and similarly $re\leq t$, 
so $re\leq u$.  This implies that 
  $
  re = re u^*u,
  $
  whence 
  $
  \q r e = r^*re u^*u \leq \q u,
  $
  and therefore 
  $$
  z\in D_{\q r} \cap D _e = D_{\q re} \c D_{\q u},
  $$
  which contradicts the fact that $z\in W$.
  \proofend


In view of this result it is interesting to find examples of inverse
semigroups which are semilattices.
Recall that a \stress{zero} in an inverse
semigroup $\S$ is an element $0\in\S$ such that
  $$
  0s=s0=0
  \for s\in\S.
  $$
  An inverse semigroup $\S$ with zero is said to be $E^*$-unitary if
for every $e,s\in\S$, one has that
  $$
  0\neq e^2=e\leq s \ \Longrightarrow \ s^2=s.
  $$
  In other words, if an element dominates a nonzero idempotent then
that element itself is an idempotent. The $E^*$-unitary inverse
semigroups have been intensely studied in the semigroup literature.
See, for example, \scite{\Lawson}{Section 9}.

The following result resembles the fact that two analytic functions on
a common connected domain, and agreeing on an open subset, must be
equal.

\state Lemma
  \label Analiticity
  Let $\S$ be an $E^*$-unitary inverse semigroup and let $s,t\in\S$ be
such that $s^*s=t^*t$, and $se=te$, for some nonzero idempotent $e\leq
s^*s$.  Then $s=t$.

\proof
  Notice that the idempotent $f=ses^*$ is nonzero because
  $
  e= s^*s e s^* s.
  $
  We have that 
  $$
  ts^*f =
  ts^*ses^* = 
  tt^*tes^* = 
  tes^* = 
  ses^* = f,
  $$
  so $f\leq ts^*$, which implies that $ts^*$ is idempotent.  In
particular it follows that $ts^* = (ts^*)^* = st^*$, so $st^*$ is
idempotent as well.  We next claim that $ss^*=tt^*$.  In fact
  $$
  tt^* = tt^*tt^* = ts^*st^* = st^*ts^* = ss^*ss^* =  ss^*.
  $$
  Setting $u = ts^*t$, we have that 
  $$
  u^*u = t^*st^* ts^*t = t^*ss^* ss^*t = t^*ss^*t =  t^*tt^*t = t^*t.
  $$
  Therefore also   $u^*u = s^* s$, while
  $$
  t=tt^*t = tu^*u
  \and
  s=ss^*s = su^*u,
  $$
  so it is enough to prove that $tu^*=su^*$.  We have
  $$
  su^* = st^*st^* = st^* = ts^* = tt^*ts^* = tt^*st^* = tu^*.
  \proofend
  $$

The following result is probably well known to semigroup theorists:

\state Proposition
  \label ISGLattice
  If $\S$ is an $E^*$-unitary inverse semigroup with zero, then $\S$
is a semilattice with respect to its usual order.

  \proof
  We must prove that $s\inf t$ exists for every $s,t\in\S$.  If there
exists no nonzero $u\in\S$, such that 
  $
  u\leq s,t,
  $
  it is clear that $s\inf t = 0$.  So suppose the contrary and fix any
such nonzero $u$.
  We then claim that
  $$
  st^*t = ts^*s = tt^*s = ss^*t.
  \subeqmark OldClaimOne
  $$
  Let $f = s^*s t^*t$.  Since $u^*u \leq s^*s$, and $u^*u\leq t^*t$,
we have that $u^*u \leq f$.
  Setting \def\alt{\tilde}
  $$
  \alt s=sf \and \alt t=tf, 
  $$
  notice that $\alt s$ and $\alt t$ share
initial projections because
  $$
  \alt s^*\alt s = fs^*sf = f = ft^*tf = \alt t^*\alt t.
  $$
  Also notice that 
  $$
  \alt su^*u = sfu^*u = su^*u = u = tu^*u = tfu^*u = \alt tu^*u.
  $$
  Employing \lcite{\Analiticity} we then deduce that $\alt s=\alt t$.  So
  $$
  st^*t = ss^*st^*t = sf = \alt s = \alt t = tf = t s^*s.
  $$
  This shows the equality between the first and second terms in \lcite{\OldClaimOne}.
  Since $0\neq u^*\leq s^*,t^*$, we may apply the above argument to
$s^*, t^*, u^*$ in order  to prove that 
  $
  s^*tt^* = t^* ss^*,
  $
  which implies that $tt^*s = ss^*t$, so the third and fourth terms in
\lcite{\OldClaimOne} agree.

The fact that $u\leq s,t$ \ implies that 
  $
  su^*u = u = tu^*u.
  $
  Left multiplying this by $t^*$ we have that 
  $$
  t^*su^*u = t^*tu^*u = u^*u,
  $$
  so $t^*s$ is idempotent by the fact that $\S$ is $E^*$-unitary.
Applying the same reasoning to $s^*$, $t^*$ and $u^*$, we have that
$ts^*$ is idempotent as well.  Thus both $t^*s$ and $ts^*$ are
selfadjoint, and hence
  $$
  st^*t =   ts^*t =  tt^*s,
  $$
  proving the equality between the first and third terms in
\lcite{\OldClaimOne}, hence concluding the proof of our claim.  We
shall next prove that the element $m(s,t):= st^*t$, satisfies
  $$
  u\leq m(s,t)\leq s,t. 
  $$
  It is obvious that $m(s,t)\leq s,t$.  Recalling that $u^*u\leq f$,
notice that 
  $$
  u = su^*u = sfu^*u = ss^*st^*tu^*u = st^*tu^*u = m(s,t)u^*u,
  $$
  so $u\leq m(s,t)$.   Therefore $m(s,t)$ is the infimum of $s$ and
$t$.
  \proofend 

\section{Pre-grading structure of $C^*(\G)$}
  In this section we return to our earlier standing hypotheses, namely
that $\act$ is an action of the inverse semigroup $\S$ on the locally
compact Hausdorff space $X$.  We will again be dealing with the
groupoid of germs of the system $(\act,\S,X)$, denoted here simply by
$\G$.

  Our aim is to show that $C^*(\G)$ admits a \stress{pre-grading} over
$\S$, as explained below:

\definition
  \label DefineGrading
  Let $A$ be any C*-algebra and let $\S$ be an inverse semigroup.  A
\stress{pre-grading}\fn{We use the term pre-grading to suggest that we
are not requiring any sort of linear independence of the subspaces
$A_s$, as is usually required for gradings over groups.}  of $A$ over
$\S$ is a family of closed linear subspaces $\{A_s\}_{s\in\S}$ of $A$,
such that for every $s,t\in\S$ on has that
  \izitem
  \zitem $A_sA_t\c A_{st}$,
  \zitem $A_s^*= A_{s^*}$,
  \zitem if $s\leq t$ \lcite{see \PartialOrderISG}, then $A_s\c A_t$,
  \zitem $A$ is the closed linear span of the
union of all $A_s$.
  \medskip \noindent
  The pre-grading is said to be \stress{full} if in addition $A_sA_t$
is dense in $A_{st}$.

We begin by introducing some terminology:

\sysstate{Notations}{\rm}{\label BunchOfNotations
  \izitem
  \zitem For each $s\in\S$ and each $f\in C_0(D_{\q s})$ we will
denote by $\a_s(f)$ the element of $C_0(D_{\p s})$ given by
  $$
  \a_s(f)\calcat x = f\big(\act_{s^*}(x)\big)
  \for x\in D_{\p s}.
  $$
  \zitem Given $s\in\S$ we will use the shorthand notation $\O_s$ for
the slice $\O(s,D_{\q s})$.
  \zitem The restriction of the source and range maps to $\O_s$ will
be denoted by $\s_s$ and $\r_s$, respectively.
  \zitem If $f$ is any complex valued function on $D_{\q s}$ we will
denote the composition $f\circ\s_s$ by $\lft_sf$. This is by definition
a function on $\O_s$ which we shall also view as a function on $\G$
by extending it to be zero outside $\O_s$.
  \zitem If $f$ is any complex valued function on $D_{\p s}$ we
will denote the composition $f\circ\r_s$ by $f\rgt_s$, with the same
convention making $f\rgt_s$ a function supported on $\O_s$.
  }

\bigskip Since $\O_s$ is a slice we have that $\s_s$ is a
homeomorphism, with domain $\O_s$, onto $\s(\O_s)=D_{\q s}$.  The
inverse of $\s_s$ is then given by
  $$
  \s_s\inv: x\in D_{\q s} \mapsto [s,x]\in\O_s.
  $$
  Compare \lcite{\PhiHomeo}

It is important not to mistake $\lft_sf$ by $f\circ \s$, since the
latter does not necessarily vanish outside $\O_s$.  Also notice that
because $\s_s$ is a homeomorphism one has that $\lft_sf\in
C_c\big(\O_s\big)$ if and only if $f\in C_c(D_{\q s})$.  In this case
we obviously have that $\lft_sf\in C_c(\G)$.  Similar observations
apply to $f\rgt_s$.

The reader will be able to tell between the notations of
\lcite{\BunchOfNotations.iv} and \lcite{\BunchOfNotations.v} by taking
note of which side of $f$ does $\delta_s$ appear.
The following is intended to conciliate these points of view.

\state Proposition
  \label ChangeConv
  Given $f\in C_c(D_{\q s})$ one has that $\lft_sf = \a_s(f)\rgt_s$.

\proof
  Clearly both $\lft_sf$ and $\a_s(f)\rgt_s$ are functions
supported on $\O_s$.  Thus, given any $x\in D_{\q s}$ we have 
  $$
  (\a_s(f)\rgt_s)([s,x]) =
  \a_s(f)\big(\r([s,x])\big) =
  \a_s(f)\big(\act_s(x)\big) =
  f(\act_{s^*}(\act_s(x))) \$=
  f(x) = 
  f(\s\big([s,x]\big) = 
  (\lft_sf)([s,x]).
  \proofend
  $$

The two notations are therefore completely interchangeable.  We shall
however prefer to use $\lft_s f$, perhaps because our notation for
$[s,x]$ already favours sources over ranges.  After all when one
speaks of the ``germ of a function $f$ at a point $x$", the emphasis
is on the point $x$ in the domain of $f$, rather that the point $f(x)$
in the range of $f$.

\state Proposition
  \label MultOsOt
  If $s,t\in\S$ then 
  \izitem
  \zitem $\O_s\O_t=\O_{st}$,
  \zitem $\O_s\inv = \O_{s^*}$.

\proof
  Given $[st,y]\in\O_{st}$ we have that $y\in D_{(st)^*st}$.  By
\lcite{\MovingDe} if follows that
  $$
  D_{(st)^*st} = 
  D_{t^*s^*st} =
  \act_{t^*}(D_{s^*s}\cap D_{tt^*}).
  $$
  Therefore, there exists $x\in D_{s^*s}\cap D_{tt^*}$ such that
$y=\act_{t^*}(x)$ and hence $x=\act_t(y)$.
  After verifying that $y\in D_{\q t}$ we then conclude that
  $$
  \big([s,x],[t,y]\big)\in (\O_s\times\O_t) \cap \Gexp2,
  $$
  and hence
  $[st,y]=[s,x]\,[t,y] \in \O_s\O_t$.  This proves that
$\O_{st}\c\O_s\O_t$.  The converse inclusion is trivial, so (i) is
proved.  The proof of (ii) follows by inspection.
  \proofend

\state Proposition
  \label PartialProducts
  Given $s,t\in\S$, let $f\in C_c(D_{\q s})$, and $g\in C_c(D_{\q t})$.
Then
  \izitem
  \zitem $(\lft_sf)\star(\lft_t g) = \lft_{st} h$, where 
  $h=\a_{t^*}(f\a_t(g))$,
  \zitem $(\lft_sf)^* =\lft_{s^*}\,\a_s(\bar f).$
  \medskip
  \noindent On the other hand, if $f\in C_c(D_{\p s})$, and $g\in
C_c(D_{\p t})$, then
  \medskip
  \zitem $(f\rgt_s)\star(g\rgt_t) = h\rgt_{st} $, where
  $h=\a_{s}(\a_{s^*}(f)g)$,
  \zitem $(f\rgt_s)^* =\,\a_{s^*}(\bar f)\rgt_{s^*}.$

\proof
  Since $\lft_sf\in C_c(\O_s)$ and $\lft_tg\in C_c(\O_t)$ we have by
\lcite{\IntroduceOperations.ii} that 
  $$
  (\lft_sf)\star(\lft_t g)\in 
  C_c(\O_s\O_t)
  \={(\MultOsOt)}
  C_c(\O_{st}).
  $$
  Given $[st,y]\in\O_{st}$ recall from the proof of \lcite{\MultOsOt}
that $[st,y]=[s,x]\,[t,y]$, where $x=\act_t(y)$.
Therefore 
  $$
  (\lft_sf)\star(\lft_t g)([st,y]) =
  (\lft_sf)([s,x])\ (\lft_tg)([t,y]) = 
  f(x) g(y) =
  f\big(\act_t(y)\big) g(y) = h(y),
  $$
  where the last equality is to be taken as the definition of $h$.
  It is tempting to write $h=\a_{t^*}(f)g$, except that we are
reserving the expression $\a_{t^*}(f)$, defined in
\lcite{\BunchOfNotations.i}, for functions $f\in C_c(D_{tt^*})$, and
all we know about $f$ is that it lies in $C_c(D_{\q s})$.  The reader
will find that the expression given in the statement is an alternative
way to describe $h$ which respects the domains of $\a_t$ and
$\a_{t^*}$, the fundamental point being that $\a_t(g)$ is in
$C_c(D_{\p t})$, and the latter is an ideal in the space of continuous
functions.
  This proves (i).

With respect to (ii), for every $\gamma\in\G$ we have that
  $$
  (\lft_s f)^*(\gamma) = \overline{(\lft_s f)(\gamma\inv)},
  $$
  so the support of $(\lft_s f)^*$ is contained in $\O_s\inv=\O_{s^*}$.
  Given $[s^*,x]\in\O_{s^*}$ we than compute
  $$
  (\lft_s f)^*([s^*,x]) = 
  \overline{(\lft_s f)([s^*,x]\inv)} =
  \overline{(\lft_s f)([s,\act_{s^*}(x)])} =
  \overline{f(\act_{s^*}(x))} \$=
  \a_s(\bar f)(x) =
  \big(\lft_{s^*}\,\a_s(\bar f)\big) ([s^*,x]).
  $$
  Points (iii) and (iv) follow respectively from (i) and (ii), with
the aid of \lcite{\ChangeConv}. 
  \proofend

The expression for $h$ in \lcite{\PartialProducts.iii} is a
fundamental formula underlying the algebrization of partially defined
maps.  It's first appearance in the literature dates back at least to
\scite{\newpim}{3.4}, and may be found also in \scite{\Sieben}{Section
5} and in \scite{\tpa}{Section 2}, the latter being its twisted
version.

\state Proposition
  \label OrderInBundle
  If $s,t\in\S$ are such that $s\leq t$, then 
  \izitem 
  \zitem $D_{\q s}\c D_{\q t}$,
  \zitem for every $f\in C_c(D_{\q s})$ one has that
$\lft_sf=\lft_tf$,
  \zitem for every $f\in C_c(D_{\p s})$ one has that
$f\rgt_s=f\rgt_t$,
  \zitem $\O_s\c\O_t$,
  \zitem $C_c(\O_s)\c C_c(\O_t)$.

\proof
  Given that $s = t\q s$ we have
  $$
  D_{\q s} = 
  D_{(ts^*s)^*ts^*s} =
  D_{\q s\q t} = 
  D_{\q s} \cap D_{\q t},
  $$
  from where (i) follows.
  To prove (iv) let $[s,x]\in\O_s$, so we have that $x\in D_{\q s}\c
D_{\q t}$ and hence $[t,x]$ belongs to $\O_t$.  In addition, setting
$e=\q s$, the fact that $x\in D_e$, and $se=te$ implies that
  $$
  [s,x]=[t,x]\in \O_t,
  \subeqmark SxEqualTx
  $$
  proving (iv), and consequently also proving (v).

To prove (ii) let $f\in C_c(D_{\q s})$, so also $f\in C_c(D_{\q t})$
by (i).  Using (v) we may view both $\lft_sf$ and $\lft_tf$ as
elements of $C_c(\O_t)$.  Given $[t,x]\in\O_t$, where $x\in D_{\q t}$,
we either have that $x\notin D_{\q s}$, in which case $[t,x]\notin
\O_s$ and hence, recalling that $\lft_sf$ is supported in $\O_s$, we
have
  $$
  \lft_s f([t,x])=
  0 = 
  f(x) = 
  \lft_t f([t,x]).
  $$
  On the other hand, if $x\in D_{\q s}$, we have
that
  $$
  \lft_t f([t,x]) =
  f(x) = 
  \lft_s f([s,x]) \={(\SxEqualTx)}
  \lft_s f([t,x]).
  $$
  This concludes the proof of (ii), and (iii) follows as well in view
of \lcite{\ChangeConv}.
  \proofend

In what follows  we give the result promised at the beginning of this
section:

\state Proposition
  Let $\iota:C_c(\G)\to C^*(\G)$ be the natural map defined in
\lcite{\NaturalMap}.
  For each $s\in\S$, let $A_s$ denote the closure of
$\iota\big(C_c(\O_s)\big)$ within $C^*(\G)$.  Then the collection
$\{A_s\}_{s\in\S}$ is a full pre-grading of $C^*(\G)$.

\proof
  It is obvious that $\{\O_s\}_{s\in\S}$ is a covering of $\G$, so
\lcite{\DefineGrading.iv} follows immediately from \lcite{\LinComb}
and the fact that $\iota(C_c(\G))$ is dense in $C^*(\G)$.

Given that $\s_s$ is a homeomorphism is is clear that $C_c(\O_s)$
consists precisely of the elements of the form $\lft_sf$, where $f$ runs
in $C_c(D_{\q s})$.  Therefore \lcite{\DefineGrading.i--ii} follow
respectively from \lcite{\PartialProducts.i--ii}.  The third axiom of
pre-gradings is an obvious consequence of \lcite{\OrderInBundle.v}, so
we are left with proving that our pre-grading is full.  For this let
$s,t\in\S$ and pick any element in $C_c(\O_{st})$, which is necessarily
of the form $\lft_{st}h$, where $h\in C_c(D_{(st)^*st})$.
  Recall e.g.~from the proof of \lcite{\MultOsOt.i} that
  $
  D_{(st)^*st} = \act_{t^*}(D_{s^*s}\cap D_{tt^*}),
  $
  so 
  $$
  \a_t(h) = h\circ \act_{t^*}\in C_c(D_{s^*s}\cap D_{tt^*}).
  $$
  We may then write $\a_t(h) =fk$, where both $f$ and $k$ are in
$C_c(D_{s^*s}\cap D_{tt^*})$.  Observing that 
  $$
  k\in C_c(D_{s^*s}\cap D_{tt^*}) \c
  C_c(D_{tt^*}) =
  C_c\big(\act_t(D_{t^*t})\big),
  $$
  the function 
  \ $
  g=k\circ\act_t = \a_{t^*}(k)
  $ \
  lies in $C_c(D_{\q t})$, and hence $\lft_tg\in C_c(\O_t)$.  In
addition we have that $\lft_sf\in\O_s$, so
  $$
  C_c(\O_s)  C_c(\O_t) \ni
  (\lft_sf)\star(\lft_tg) =
  \lft_{st}\Big(\a_{t^*}\big(f\a_t(g)\big)\Big) =
  \lft_{st}\big(\a_{t^*}(fk)\big) =
  \lft_{st}h.
  $$
  This shows that 
  $C_c(\O_{st}) \c C_c(\O_s) C_c(\O_t)$, from where one sees that our
pre-grading is in fact full.
  \proofend


\section{Universal property of $C^*(\G)$}
   As before we fix an action $\act$ of an inverse semigroup $\S$ on
a locally compact Hausdorff topological space $X$.  We will assume in
addition that $\S$ is countable and that $X$ is second countable,%
  \fn{In case of absolute necessity one may perhaps dispense with the
second countability assumption at the expense of working with the
$\sigma$-algebra of Baire (instead of Borel) measurable sets, assuming
in addition that every $D_e$ is Baire measurable.}
  due to the use of measure theory methods.
  We shall retain the notation $\G$ for the groupoid of germs of the
system $(\act,\S,X)$.

Recall from \lcite{\BunchOfNotations.i} that for $s\in\S$ we denote by
$\a_s$ the isomorphism from $C_c(D_{\q s})$ to $C_c(D_{\p s})$ given
by $\a_s(f)=f\circ\act_{s^*}$.

\definition 
  \label DefineCovarRep
  A \stress{covariant representation} of the system $(\act,\S,X)$ on
a Hilbert space $H$ is a pair 
  $(\pi,\irep)$, where $\pi$ is a nondegenerate *-representation of
$C_0(X)$ on $H$, and
  \ $
  \irep:\S\to B(H)
  $ \
  satisfies
  \izitem
  \zitem $\irep_{st}=\irep_s\irep_t$,
  \zitem $\irep_{s^*}=\irep_s^*$,
  \zitem $\pi(\a_s(f)) = \irep_s\pi(f)\irep_{s^*}$,
  \zitem $\overline{\pi\big(C_0(D_e)\big)H} = \irep_e(H)$,
  \medskip\noindent
  for every
  $s,t\in\S$, 
  $f\in C_0(D_{\q s})$,
  and 
  $e\in E(\S)$.

\cryout{From now on we fix a covariant representation $(\pi,\irep)$ of
$(\act,\S,X)$ on $H$.}

We will write $\tilde\pi$ for the canonical weakly continuous
extension 
  of $\pi$ to the algebra $\B(X)$ of all bounded Borel measurable
functions on $X$.  It is well known that for each open subset $U\c X$,
one has that the range of $\tilde\pi(1_U)$ coincides with
$\overline{\pi\big(C_0(U)\big)H}$, where $1_U\in\B(X)$ denotes the
characteristic function of $U$.
  Therefore \lcite{\DefineCovarRep.iv} may be expressed by saying that
  $$
  \irep_e = \tilde\pi(1_{D_e})
  \for e\in E(S).
  \eqmark RepePiOneE
  $$
  In particular it follows that 
  $$
  \irep_e \pi(f) = \pi(f)\irep_e
  \for f\in C_0(X).
  \eqmark ComutingPiAndRep
  $$

Our next main goal will be to show that there exists a
*-representation
  $
  \irep\times\pi:C_c(\G)\to B(H),
  $
  such that for every $s\in\S$, and $f\in C_c(D_{\q s})$, one has that
  $
  (\irep\times\pi)(\d_sf) = \irep_s\pi(f).
  $

  \state Lemma
  \label KeyWellDefLemma
  Let $J$ be a finite subset of $\S$ and suppose that for each $s\in
J$ we are given $f_s\in C_c(D_{\q s}\big)$ such that
  $
  \sum_{s\in J}\d_sf_s=0,
  $
  in $C_c(\G)$.
  Then 
  $
  \sum_{s\in J}\irep_s\pi(f_s)=0,
  $
  in $B(H)$.

\proof Fix, for the time being, two elements
$\xi,\eta\in H$.
  For each $s\in \S$, let $\mu_s = \mu_{s,\xi,\eta}$ be the finite Borel
measure on $\O_s$ given by
  $$
  \mu_s(A) = \<\irep_s\tilde\pi(1_{\s(A)})\xi,\eta\>,
  $$
  for every Borel measurable $A\c\O_s$, where $1_{\s(A)}$ stands for
the characteristic function on $\s(A)$.  

Since $\s$ is a homeomorphism
from $\O_s$ to $D_{\q s}$, one has that $\s(A)$ is a measurable
subset of $D_{\q s}$, and hence also of $X$.  Therefore $1_{\s(A)}\in
\B(X)$, so that $\tilde\pi(1_{\s(A)})$ is well defined.  That $\mu_s$
is indeed a countably additive measure follows from the corresponding
well known property of $\tilde\pi$.

If $B\c D_{\q s}$ is a measurable set, let $A=\s_s\inv(B)$, so that $A$
is a measurable subset of $\O_s$ and $B=\s(A)$.  Notice that
  $$
  \d_s1_B = 1_B\circ \s_s = 1_A,
  $$
  so 
  $$
  \int_{\O_s} \d_s1_B \, d\mu_s =
  \int_{\O_s} 1_A \, d\mu_s =
  \mu_s(A) =
  \<\irep_s\tilde\pi(1_{\s(A)})\xi,\eta\> =
  \<\irep_s\tilde\pi(1_B))\xi,\eta\>,
  $$
  from where one easily deduces that 
  $$
  \int_{\O_s}\d_sf\ d\mu_s = \<\irep_s\tilde\pi(f)\xi,\eta\>
  \for  f\in\B(X).
  \subeqmark IntegralAndPi
  $$
  We next claim that for every $s,t\in\S$, and every measurable set
$A\c\O_s\cap \O_t$, one has that
  $$
  \mu_s(A) = \mu_t(A).
  \subeqmark MustA
  $$
  In order to prove it observe that $B:=\s(A) = \s_s(A) = \s_t(A)$ is
a Borel subset of
  $D_{\q s}\cap D_{\q t}$ and
  $$
  A=\{[s,x]: x\in B\} =\{[t,x]: x\in B\}.
  $$
  For every $x\in B$ we moreover have that $[s,x]=[t,x]$, so there
exists $e\in E(\S)$ such that $x\in D_e$, and $se=te$.  It therefore
follows that 
  $$ 
  B\c\union_{
  \buildrel \scriptstyle e\in E(\S) \over {se = te}
  }
  D_e.
  $$
  Since we are assuming that $\S$ is countable, so is $\E(S)$ and we
may decompose $B$ as a disjoint union of measurable subsets
$\{B_n\}_{n\in\N}$, such that each $B_n$ is a subset of some
$D_{e_n}$, and $se_n=te_n$.  Obviously $A$ is then the disjoint union
of the sets
  $$
  A_n = \s_s\inv(B_n) = \s_t\inv(B_n).
  $$
  Notice that for each $n\in\N$ we have
  $$
  \irep_s\tilde\pi(1_{\s(A_n)}) =
  \irep_s\tilde\pi(1_{B_n}) =
  \irep_s\tilde\pi(1_{D_{e_n}}1_{B_n}) \$=
  \irep_s\tilde\pi(1_{D_{e_n}})\tilde\pi(1_{B_n}) \={(\RepePiOneE)}
  \irep_s\irep_{e_n}\tilde\pi(1_{B_n}) =
  \irep_{se_n}\tilde\pi(1_{B_n}),
  $$
  and similarly for $t$.  Since $se_n=te_n$ we have that
  $
  \irep_s\tilde\pi(1_{\s(A_n)}) = \irep_t\tilde\pi(1_{\s(A_n)}),
  $
  whence $\mu_s(A_n) = \mu_t(A_n)$.  The countable additivity of $\mu_s$
and $\mu_t$ then take care of \lcite{\MustA}.

Let $M$ be the measurable subset of $\G$ given by $M={\union}_{s\in
J}\O_s$, where $J$ is as in the statement.
  It is an easy exercise in measure theory to prove that there exists
a measure $\mu$ on $M$ such that $\mu(A) = \mu_s(A)$, for every $s\in
J$, and $A\c\O_s$.  We then have that
  $$
  \Big\langle\sum_{s\in J}\irep_s\pi(f_s)\xi,\eta\Big\rangle
  \ \={(\IntegralAndPi)} \ 
  \sum_{s\in J}\ \int_{\O_s}\d_sf_s\ d\mu_s =
  \sum_{s\in J}\ \int_M\d_sf_s\ d\mu =
  \int_M \ \sum_{s\in J}\d_sf_s\ d\mu = 0.
  $$
  Since $\xi$ and $\eta$ are arbitrary we conclude that 
  $\sum\limits_{s\in J}\irep_s\pi(f_s)=0$, as stated.
  \proofend

We thus arrive at the main result of this section.

\state Theorem
  \label UnivPropGpd
   Let $\S$ be a countable inverse semigroup, let $\act$ be an action
of $\S$ on the second countable locally compact Hausdorff space $X$, and
let $\G$ be the corresponding groupoid of germs
\lcite{\TheMainGroupoid}.  Given any covariant representation
$(\pi,\irep)$ of $(\act,\S,X)$ on a Hilbert space $H$ there exists a
unique *-representation $\irep\times\pi$ of $C^*(\G)$ on $H$ such that
  $$
  (\irep\times\pi)\big(\iota(\d_sf)\big) = \irep_s\pi(f)
  \and
  (\irep\times\pi)\big(\iota(g\d_s)\big) = \pi(g)\irep_s,
  $$
  for every $s\in\S$, every $f\in C_c(D_{\q s})$, and every $g\in
C_c(D_{\p s})$, 
  where $\iota:C_c(\G)\to C^*(\G)$ is the canonical map.

\proof
  Given any $f\in C_c(\G)$ use \lcite{\OSsIsSliceLegal} and
\lcite{\LinComb} to write $f = \sum\limits_{k=1}^n \d_{s_k}f_k$, where
$s_1,\ldots,s_n\in\S$, and $f_k\in C_c(D_{s_k^*s_k})$, for all
$k=1,\ldots,n$.  Define
  $$
  (\irep\times\pi)(f) = \sum_{k=1}^n \irep_{s_k}\pi(f_k).
  $$
  That $\irep\times\pi$ is well defined is a consequence of
\lcite{\KeyWellDefLemma}.  It is obviously also linear, and we claim
that it is a *-homomorphism.  In order to prove the preservation of
multiplication, we may use linearity to reduce our task to proving
only that
  $$
  (\irep\times\pi)(\d_sf\star\d_tg) =
  (\irep\times\pi)(\d_sf)\ (\irep\times\pi)(\d_tg),
  $$
  for every $f\in C_c(D_{\q s})$ and $g\in C_c(D_{\q t})$.  By
\lcite{\PartialProducts.i} the left-hand side equals
  $$
  (\irep\times\pi)\Big(\d_{st}\a_{t^*}\big(f\a_t(g)\big)\Big) =
  \irep_{st} \pi\Big(\a_{t^*}\big(f\a_t(g)\big)\Big) =
  \irep_{st} \irep_{t^*}\pi\big(f\a_t(g)\big)\irep_t \$=
  \irep_{st} \irep_{t^*}\pi(f)\pi\big(\a_t(g)\big)\irep_t =
  \irep_s\irep_t \irep_{t^*}\pi(f)\irep_t\pi(g)\irep_{t^*}\irep_t 
  \={(\ComutingPiAndRep)}
  \irep_s\pi(f)\irep_t \irep_{t^*}\irep_t\irep_{t^*}\irep_t\pi(g) \$=
  \irep_s\pi(f)\irep_t \pi(g) =
  (\irep\times\pi)(\d_sf)\ (\irep\times\pi)(\d_tg).
  $$
  Our claim will then be proved once we show that
  $$ 
  (\irep\times\pi)\big((\d_sf)^*\big) =
  \big((\irep\times\pi)(\d_sf)\big)^*.
  $$
  By \lcite{\PartialProducts.ii} the left-hand side equals
  $$
  (\irep\times\pi)\big(\d_{s^*}\,\a_s(\bar f)\big) =
  \irep_{s^*}\pi\big(\a_s(\bar f)\big) =
  \irep_{s^*}\irep_s\pi(\bar f)\irep_{s^*} =
  \pi(\bar f)\irep_{s^*}\irep_s\irep_{s^*} \$=
  \pi(f)^*\irep_{s^*} = 
  \big(\irep_s\pi(f)\big)^* =
  \big((\irep\times\pi)(\d_sf)\big)^*.
  $$
  This proves our claim. 
  By \lcite{\DefineTripleNorm} we then conclude that 
  $$
  \|(\irep\times\pi)(f)\| \leq \tn f
  \for f\in C_c(\G),
  $$
  which implies that $\irep\times\pi$ factors through $\iota$,
producing a *-representation of $C^*(\G)$, by abuse of language also
denoted by $\irep\times\pi$, clearly satisfying the first identity in
the statement.

In order to prove the second identity let $s\in\S$ and 
  $g\in C_c(D_{\p s})$.  Set $f=\a_{s^*}(g)$, so that 
  $g=\a_{s}(f)$, and 
  $f\in C_c(D_{\q s})$.  Therefore
  $$
  g\d_s=
  \a_s(f)\d_s \={(\ChangeConv)}
  \d_sf,
  $$
  so
  $$
  (\irep\times\pi)\big(\iota(g\d_s)\big) = 
  (\irep\times\pi)\big(\iota(\d_sf)\big) = 
  \irep_s\pi(f) \$=
  \irep_s\irep_s^*\irep_s\pi(f) \={(\ComutingPiAndRep)}
  \irep_s\pi(f)\irep_s^*\irep_s =
  \pi\big(\a_s(f)\big)\irep_s =
  \pi(g)\irep_s.
  \proofend
  $$


\section{Inverse semigroup crossed products}
  The main goal of this section is to show that, in the context of the
previous section, $C^*(\G)$ is naturally isomorphic to the inverse
semigroup crossed product $C_0(X)\rtimes_\a S$.

In the first part of this section we shall therefore briefly review the
theory of inverse semigroup crossed products based on \cite{\Sieben} and
\cite{\PatBook}, not only for the convenience of the reader, but also
because we will present a few improvements.

\long\def\RemarkOnUnits{%
  Sieben assumes that $\S$ is unital \scite{\Sieben}{4.1} and that
$\a_e$ is the identity map on $A$.  Attempting to avoid units Paterson
instead assumes that the family of the domains of the $\a_s$ forms an
upward directed chain \scite{\PatBook}{Definition 3.3.1.ii}.  These
assumptions are designed to be used in proving the equivalence between
covariant representations of the system and *-representations of the
covariance algebra.  See \scite{\Sieben}{5.6} and
\scite{\PatBook}{Proposition 3.3.3}.  As we will see below there is a
way to get around this problem without assuming either of this extra
conditions.}

\definition
  \label DefineActionOnAlg
  An \stress{action} of an inverse semigroup\fn{\RemarkOnUnits} $\S$
on a C*-algebra $A$ is a semigroup homomorphism
  $$
  \a:\S \to \I(A),
  $$ 
  (see \lcite{\DefineIX} for a definition of $\I(A)$) such that 
  \izitem 
  \zitem for every $s\in\S$, the domain (and hence also the range) of
$\a_s$ is a closed two sided ideal of $A$, and $\a_s$ is a
*-homomorphism,
  \zitem the linear span of the union of the domains of all the $\a_s$
is dense in $A$.

As in the case of actions on locally compact spaces, defined in
\lcite{\DefineAction}, for every $e\in E(\S)$, we denote by $J_e$ the
domain of $\a_e$.  For each $s\in\S$ one therefore has that $\a_s$ is a
*-isomorphism from $J_{\p s}$ to $J_{\p s}$.  See also footnote
\lcite{\FootnoteOnDe}.

Given an action of $\S$ on a locally compact space $X$ in the sense of
\lcite{\DefineAction}, it is easy to produce an action of $\S$ on
  $A=C_0(X)$, this time in the sense of \lcite{\DefineActionOnAlg}:
observing that $J_e:=C_0(D_e)$ is an ideal in $C_0(X)$, for each
$s\in\S$, one takes $\a_s:J_{\q s}\to J_{\p s}$ to be given by
\lcite{\BunchOfNotations.i}.  To check that
\lcite{\DefineActionOnAlg.ii} holds one uses \lcite{\DefineAction.ii}
and the Stone--Weierstrass Theorem.

 \cryout{From now on we fix an action of $\S$ on a C*-algebra $A$.}

One then considers the linear space
  $$
  L = \bigoplus_{s\in\S}J_{\p s}.
  \eqmark DefineL
  $$
  If $e$ is an idempotent notice that $J_e$ appears in the above
direct sum
as many times as there are elements  $s\in\S$ with $\p s=e$.

  Any element $x$ in $L$ is of the form $x = (a_s)_{s\in\S}$, where
  $a_s\in J_{\p s}$, and $a_s=0$ for all but finitely many $s$.  Given
$s\in\S$ and $a\in J_{\p s}$, we shall denote by $a\d_s$ the element of $L$
which is identically zero except for its $s^{th}$ component which is
equal to $a$.  Any element of $L$, say $x = (a_s)_{s\in\S}$, is
therefore given by
  $$
  x = \sum_{s\in\S} a_s\d_s,
  \eqmark GenEltInL
  $$  
  where the sum has finitely many nonzero terms.  Based on
\cite{\newpim} and \cite{\McClanahan}, Sieben defines a *-algebra
structure on $L$ according to which
  $$
  (a\d_s) (b\d_t) = \a_s\big(\a_{s^*}(a)b\big)\d_{ts}
  \and 
  (a\d_s)^* = \a_{s^*}(a^*)\d_{s^*},
  $$
  for every $s,t\in\S$, $a\in J_{\p s}$, and $b\in J_{\p t}$.

The crossed product $A\rtimes_\a L$ is then defined (see below) as a
certain completion of $L$.  However, contrary to what happens with
similar constructions, one does not expect $L$ to survive the completion
process intact: if $e$ and $f$ are idempotents and $a\in J_e\cap J_f$,
so that $a\d_e$ and $a\d_f$ are elements of $L$, the construction is
such that $a\d_e-a\d_f=0$, when passing to the crossed product.

Let us now review the construction of the crossed product. Sieben
first defines \scite{\Sieben}{4.5} a covariant representation of the
system $(\a,\S,A)$ on a Hilbert space $H$ (up to the fact that our
$\S$ needs not have a unit) precisely as in \lcite{\DefineCovarRep},
except that $C_0(X)$ is replaced by $A$, and $C_0(D_e)$ is replaced by
$J_e$.  Risking being a bit monotonous the definition is:

\definition
  \label DefineAlgCovarRep
  A \stress{covariant representation} of the system
$(\a,\S,A)$ on a Hilbert space $H$ is a pair $(\pi,\irep)$, where
$\pi$ is a nondegenerate *-representation of $A$ on $H$, and
  \ $
  \irep:\S\to B(H)
  $ \
  satisfies
  \izitem
  \zitem $\irep_{st}=\irep_s\irep_t$,
  \zitem $\irep_{s^*}=\irep_s^*$,
  \zitem $\pi(\a_s(a)) = \irep_s\pi(a)\irep_{s^*}$,
  \zitem $\overline{\pi\big(J_e\big)H} = \irep_e(H)$,
  \medskip\noindent
  for every
  $s,t\in\S$, 
  $a\in J_{\p s}$,
  and 
  $e\in E(\S)$.

It is then easy to see \scite{\Sieben}{5.3} that for every covariant
representation $(\pi,\irep)$ the formula
  $$
  (\pi\times\irep)\Big(
  \sum_{s\in\S} a_s\d_s
  \Big) = \sum_{s\in\S}\pi(a_s)\irep_s
  $$
  defines a nondegenerate *-representation of $L$ on $H$.

If $e,f\in E(S)$ are such that $e\leq f$ (meaning that $ef=e$), then
$\a_e\a_f=\a_e$, which gives $J_e\c J_f$.  If moreover and $a\in J_e$ we
may speak of two different elements of $L$, namely $a\d_e$ and $a\d_f$.
Moreover notice that
  $$
  (\pi\times\irep)(a\d_e) =
  \pi(a)\irep(e) =
  \pi(a)\irep(ef) =
  \pi(a)\irep(e)\irep(f) =
  \pi(a)\irep(f) =
  (\pi\times\irep)(a\d_f),
  $$
  where our use of the identity 
  $\pi(a)\irep(e) = \pi(a)$ is justified by
\lcite{\DefineAlgCovarRep.iv}.

Restricting one's attention to representations of $L$ which behave as
  $\pi\times\irep$ in the above respect is an important insight due to
Paterson.

\definition \scite{\PatBook}{3.87}
  A *-homomorphism $\phi$ from $L$ into another *-algebra will be called
\stress{admissible} if for every $e,f\in E(\S)$, with $e\leq f$, and
every $a\in J_e$, one has that $\phi(a\d_e)=\phi(a\d_f)$.

Following Sieben \scite{\Sieben}{5.6}, Patterson
\scite{\PatBook}{Proposition 3.3.3} proves that every admissible
nondegenerate *-representation $\Pi$ of $L$ on a Hilbert space $H$ is
given as above for a covariant representation $(\pi,\irep)$ of
$(\a,\S,A)$.  Under the assumption that $\S$ is unital the
construction of the first component of the covariant representation,
namely $\pi$, is a breeze: for every $a$ in $A$ one simply defines
$\pi(a) = \Pi(a\d_1)$.  Paterson avoids units by requiring that the
$J_e$ be upward directed.  However it is possible to get around this
problem with bare hands:

\state Lemma
  Given an admissible nondegenerate *-representation $\Pi$ of\/ $L$ on a
Hilbert space $H$, there exists a *-representation $\pi$ of $A$ on $H$
such that
  $$
  \pi(a) = \Pi(a\d_e)
  \for e\in E(\S) \for a\in J_e.
  $$

\proof
  We first claim that for every $s\in\S$, $e\in E(\S)$, and $a\in
J_e\cap J_{\p s}$ one has that 
  $$
  \Pi(a\d_{es})=  \Pi(a\d_s).
  \subeqmark PiasIsPiaes
  $$
  Since 
  $
  J_{es(es)^*}=J_{e\p s} = J_e\cap J_{\p s},
  $
  both elements appearing as arguments to $\Pi$ in
\lcite{\PiasIsPiaes} are indeed in $L$.  We have 
  $$
  (a\d_{es}- a\d_s)(a\d_{es} - a\d_s)^* =
  (a\d_{es}- a\d_s)\big(\a_{s^*e}(a^*)\d_{s^*e} -
\a_{s^*}(a^*)\d_{s^*}\big) =
  - aa^*\d_{ss^*e}
  + aa^*\d_{ss^*},
  $$
  so admissibility implies that
  $\Pi(a\d_{es}- a\d_s)\Pi(a\d_{es} - a\d_s)^*=0$, from which
\lcite{\PiasIsPiaes} follows.
  Let 
  $$
  A_0 = \sum_{e\in E(\S)}J_e,
  $$
  so that $A_0$ is a dense *-subalgebra of $A$.  Given $a$ in $A_0$,
write it as a finite sum $a = \sum\limits_{e\in E(\S)}a_e$, with $a_e\in
J_e$, and define
  $$
  \pi(a) = \sum_{e\in E(\S)}\Pi(a_e\d_e).
  $$
  We claim that $\pi(a)$ does not depend on the choice of the $a_e$'s.
Proving this claim is tantamount to proving that when $a$ vanishes, so
does the right-hand side above.  Since $\Pi$ is nondegenerate it is in
fact enough to prove that
  $$
  \sum_{e\in E(\S)}\Pi(a_e\d_e)\Pi(b\d_s) = 0,
  $$
  for every $s\in\S$ and $b\in J_{\p s}$.  The left-hand side above equals
  $$
  \sum_{e\in E(\S)}\Pi\big((a_e\d_e)(b\d_s)\big) =
  \sum_{e\in E(\S)}\Pi(a_eb\d_{es}) \={(\PiasIsPiaes)}
  \sum_{e\in E(\S)}\Pi(a_eb\d_s) \$=
  \Pi\Big(\!\sum_{e\in E(\S)}a_eb\d_s\Big) =
  \Pi(ab\d_s)=0.
  $$
  This proves that $\pi$ is a well defined map on $A_0$.  To prove that
$\pi$ is a *-representation let $e,f\in E(\S)$, $a\in J_e$, and $b\in
J_f$.  Then 
  $$
  \pi(a)\pi(b) =
  \Pi(a\d_e)  \Pi(b\d_f) =
  \Pi(ab\d_{ef}) =\pi(ab).
  $$
  We leave it for the reader the easy proof that $\pi$ preserves the
star operation.  Summarizing, $\pi$ is a *-representation of the dense
subalgebra $A_0\c A$ on $H$.
  Any finite sum of ideals among the $J_e$ gives a closed *-subalgebra
of $A$.
  This implies that $\pi$ is norm-decreasing on $A_0$ and hence
extends to a *-representation
  of $A$, which clearly satisfies the required conditions.
  \proofend

Inserting the result above into Sieben's proof of
\scite{\Sieben}{5.6}, or Paterson's proof of
\scite{\PatBook}{Proposition 3.3.3}, we arrive 
  at the following:

\state Proposition
  \label CoressponReps
 Let $\S$ be a (not necessarily unital) inverse semigroup and let $\a$
be an action of $\S$ on a C*-algebra $A$.  Then the association
  $$
  (\pi,\irep)\ \longmapsto\ \Pi=\pi\times\irep
  $$
  is a one-to-one correspondence between covariant representations
$(\pi,\irep)$ of $(\a,\S,A)$ and admissible nondegenerate
*-representations $\Pi$ of $L$.

Recall from \scite{\Sieben}{5.4} that the \stress{crossed product of
$A$ by $\S$ relative to the action $\a$}, denoted $A\rtimes_\a\S$, is
defined to be the Hausdorff completion of $L$ in the norm
  $$
  \tn x = \sup_\Pi \|\Pi(x)\|,
  $$
  where the supremum is taken over all representations of $L$ of the
form $\Pi=\pi\times\irep$ (equivalently over all admissible
nondegenerate representations).
  As such, it is evident that to the classes of objects put in
correspondence by \lcite{\CoressponReps}, one can add the
nondegenerate *-representations of $A\rtimes_\a\S$.

This concludes our review of inverse semigroup crossed products, so we
will now return to considering actions of inverse semigroups on
topological spaces.

\state Theorem 
  \label ISGCrossProdAndGpd
  Let $\S$ be a countable inverse semigroup, let $X$ be a second
countable locally compact Hausdorff space, and let $\act$ be an action
of\/ $\S$ on $X$ in the sense of \lcite{\DefineAction}.  Denoting by
$\G$ the groupoid of germs of $(\act,\S,X)$ one has that $C^*(\G)$ is
  isomorphic to $C_0(X)\rtimes_\a S$, where $\a$ is the action of $\S$
on $C_0(X)$ given by \lcite{\BunchOfNotations.i}.

\proof
  Choose a faithful nondegenerate *-representation 
  $$
  \Psi: C_0(X)\rtimes_\a S \to B(H),
  $$
  where $H$ is a Hilbert space.
  That representation, once composed with the natural map
  $$
  j: L \to C_0(X)\rtimes_\a S,
  $$
  yields a *-representation $\Pi=\Psi\circ j$, of $L$ on $H$ which is
clearly admissible and nondegenerate.  By \lcite{\CoressponReps} there
exists a covariant representation $(\pi,\irep)$ of
$\big(\a,\S,C_0(X)\big)$ on $H$ such that
  $$
  \Pi(f\d_s) = \pi(f)\irep_s,
  $$
  for every $s\in\S$, and $f\in J_{\p s}=C_0(D_{\p s})$.  Invoking
\lcite{\UnivPropGpd} we deduce that there exists a *-representation
$\Phi = \pi\times\irep$ of $C^*(\G)$ on $H$ such that
  $$
  \Phi\big(\iota(f\d_s)\big) = \pi(f)\irep_s =
  \Pi(f\d_s) = \Psi\big(j(f\d_s)\big),
  $$
  for every $s\in\S$, and every $f\in C_c(D_{\p s})$.

Observe that the notation ``$f\d_s$" means different things here: an
element of $L$ as in \lcite{\GenEltInL}, or an element of $C_c(\G)$ as
in \lcite{\BunchOfNotations.v}. However the context should suffice to
distinguish between these uses.

It follows that $\Phi$ maps $C^*(\G)$ into the image of
$C_0(X)\rtimes_\a S$ through $\Psi$ in $B(H)$, and since $\Psi$ is
faithful we can produce a *-homomorphism
  $$
  \phi: C^*(\G)\to C_0(X)\rtimes_\a S,
  $$
  such that 
  $$
  \phi\big(\iota(f\d_s)\big) = j(f\d_s)
  \for s\in\S \for f\in C_0(D_{\q s}).
  $$
 
Leaving this aside for a moment consider the map
  $$
  \gamma: \sum_{s\in\S}f_s\d_s \in L \longmapsto 
  \sum_{s\in\S}f_s\d_s \in C_c(\G),
  $$
  where again the double meaning of $f_s\d_s$ should bring no
confusion.  Using \lcite{\PartialProducts.iii--iv} it is immediate
that $\gamma$ is a *-homomorphism, and by \lcite{\OrderInBundle.iii} one
sees that it is admissible.  Therefore the composition
$\iota\circ\gamma$ extends to give a *-homomorphism
  $$
  \psi:  C_0(X)\rtimes_\a S \to C^*(\G),
  $$
  satisfying 
  $$
  \psi(j(f\d_s)) = \iota(f\d_s) 
  \for s\in\S \for f\in C_0(D_{\q s}).
  $$
  This proves that $\psi$ and $\phi$ are each other's inverse, and
hence isomorphisms.
  \proofend

We may use our methods to obtain the following generalization of
  \scite{\PatBook}{Theorem 3.3.1} and
  \scite{\QuigSieb}{8.1}.

\state Proposition
  \label GenQSP
  Let $\G$ be a \'etale groupoid with second countable unit space and
let $\S$ be a countable%
  \fn{If such an $\S$ exists then $\G$ itself is second countable.}
  *-subsemigroup of $\S(\G)$ satisfying 
  \lcite{\TwoGroupoids.i--ii}.  Let moreover $\act$ be the restriction
to $\S$ of the action of $\SG$ on $\Gexp0$ given by
\lcite{\TechDefActU}, and denote by $\a$ the induced action of $\S$ on
$C_0(\Gexp0)$, as in \lcite{\BunchOfNotations.i}.  Then
  $$
  C^*(\G) \simeq C_0(\Gexp0)\rtimes_\a\S.
  $$

\proof
  Let $\H$ be the groupoid of germs for the given action of $\S$ on
$\Gexp0$.  Applying \lcite{\ISGCrossProdAndGpd} we conclude that 
  $$
  C^*(\H) \simeq C_0(\Gexp0)\rtimes_\a\S,
  $$
  but we also have that 
  $$
  \H\simeq\G,
  $$
  by \lcite{\TwoGroupoids}, so the statement follows.
  \proofend

\section{Action on the spectrum}
  \label ActOnSpecSect
  As before we will let $\S$ be an inverse semigroup, but we will no
longer postulate the existence of actions of $\S$ on exogenous
topological spaces.  Instead we will construct actions on spaces which
are intrinsic to $\S$.  These spaces will actually be constructed from
the idempotent semilattice of $\S$, which we will denoted simply by
$E$.

\definition
  \label BasicSpectrum
  Let $E$ be any semilattice.
  A \stress{semicharacter} of $E$ is a nonzero map
  $$
  \phi:E\to\{0,1\},
  $$
  such that $\phi(ef) = \phi(e)\phi(f)$, for all $e,f\in E$.
The set of all semicharacters equipped with the topology of
pointwise convergence (equivalently the relative topology from the
product space $\{0,1\}^E$) is called the \stress{spectrum} of $E$ and
is denoted $\spec$. 

 It is easy to see that $\spec$ is a locally compact Hausdorff
topological space.

\definition
  For every $e\in E$ we will denote by $D_e$ the subset of
$\spec$ formed by all semicharacters $\phi$ such that $\phi(e)=1$.

Given that the correspondence $\phi\mapsto\phi(e)$ is continuous in
the topology of pointwise convergence, we see that $D_e$ is a clopen
subset of $\spec$. 

Notice that $\spec$ may fail to be compact since there may exist a net
of semicharacters converging pointwise to the identically zero map
(which is not a character by definition).  No such net may exist
inside $D_e$ because its semicharacters take the value 1 at $e$.  So
$D_e$ is actually closed in $\{0,1\}^E$, hence compact.

\state Proposition
  \label MovingSemicharacters
  Let $s\in\S$ and $\phi\in D_{\q s}$. 
  \izitem
  \zitem  The map
  $\act_s(\phi):e\in E\mapsto \phi(s^*es)\in \{0,1\}$
  is a semicharacter in $D_{\p s}$.
  \zitem The map $\act_s:\phi\in D_{\q s}\mapsto \act_s(\phi)\in D_{\p
s}$ is a homeomorphism.
  \zitem The map $\act:s\in\S\mapsto \act_s\in\I(\spec)$ is a
semigroup homomorphism.
  \zitem $\act$ is an action of $\S$ on $\spec$, as defined in
\lcite{\DefineAction}.

\proof For $e,f\in E$ we have
  $$
  \act_s(\phi)(ef)=\phi(s^*efs)=\phi(s^*ess^*fs)=\phi(s^*es)\phi(s^*fs)=
  \act_s(\phi)(e)\  \act_s(\phi)(f),
  $$
  so $\act_s(\phi)$ is multiplicative. In addition
  $$
  \act_s(\phi)(\p s) = \phi(s^*\p ss) = \phi(\q s)=1,
  $$
  so $\act_s(\phi)\in D_{\p s}$.  For every net $\{\phi_i\}_i$
converging to $\phi$ in $D_{\q s}$, and for every $e\in E$, one has
that
  $$
  \lim_i\act_s(\phi_i)(e) =
  \lim_i\phi_i(s^*es) =
  \phi(s^*es) =
  \act_s(\phi)(e),
  $$
  so we see that $\{\act_s(\phi_i)\}_i$ converges to $\act_s(\phi)$,
proving that $\act_s$ is continuous.  We next claim that $\act_s$ is
bijective and $\act_s\inv=\act_{s^*}$.  In fact, for all $\phi\in
D_{\q s}$ and all $e\in E$, we have
  $$
  \act_{s^*}\big(\act_s(\phi)\big)(e) =
  \act_s(\phi)(ses^*) =
  \phi(s^*ses^*s) =
  \phi(s^*s) \phi(e) \phi(s^*s) =
  \phi(e),
  $$
  so $\act_{s^*}\circ\act_s$ is the identity on $D_{\q s}$.  By
exchanging $s$ and $s^*$ we have that $\act_s\circ\act_{s^*}$ is
also the identity on $D_{\p s}$, verifying our claim, and also giving 
  $$
  \act_{s^*} = \act_s\inv.
  $$
  This proves also that $\act_s\inv$ is continuous, so $\act_s$ is a
homeomorphism as required by (ii). 

Before we tackle (iii) observe that for every $e,f\in E$ one has that
$D_e\cap D_f = D_{ef}$.  In addition we claim that
  $$
  \act_s(D_{\q s}\cap D_e)= D_{ses^*}.
  \subeqmark MovingDsInSpecAction
  $$
  In fact, a semicharacter $\phi$ lies in the set displayed on the
left-hand side above if and only if
  $$
  \act_s\inv(\phi) \in D_{\q s}\cap D_e = D_{\q se} \iff
  \act_{s^*}(\phi)(\q s e)=1 \iff
  \phi(s\q s es^*)=1 \$\iff
  \phi(ses^*)=1 \iff
  \phi\in D_{ses^*}.
  $$ 
  In particular, given $s,t\in\S$, the domain of $\act_t\circ\act_s$
is given by
  $$
  \act_s\inv(D_{\p s}\cap D_{\q t}) =
  \act_{s^*}(D_{\p s}\cap D_{\q t}) =
  D_{s^*\q ts} =
  D_{(ts)^*ts},
  $$
  which is precisely the domain of $\act_{ts}$.  Moreover for every
$\phi\in D_{(ts)^*ts}$, and every $e\in E$, we have
  $$
  \act_s\big(\act_t(\phi)\big)(e) = 
  \act_t(\phi)(s^*es) =
  \phi(t^*s^*est) =
  \act_{st}(\phi)(e),
  $$
  proving that 
  $\act_s\circ\act_t= \act_{st}$, and (iii) follows.  To prove (iv) it
is now enough to check \lcite{\DefineAction.ii}.  For this it suffices
to observe that if $\phi\in\spec$ then $\phi$ is nonzero by
definition, and hence there exists $e\in E$ such that $\phi(e)=1$.
Thus $\phi$ lies in the domain of any $\act_s$ for which $\q s=e$, for
example $s=e$.
  \proofend

\definition
  \label DefineRepSoDoCara
  Let $H$ be a Hilbert space.  A map $\irep:\S\to B(H)$ will
be called a \stress{representation of $\S$ on $H$} if for every
$s,t\in\S$ one has
  \izitem
  \zitem $\irep_{st}=\irep_s\irep_t$,
  \zitem $\irep_{s^*}=\irep_s^*$.

We have already encountered such objects when we studied covariant
representations, as defined in \lcite{\DefineCovarRep}.  The
difference is that here there is no representation $\pi$ of $C_0(X)$
to go along with $\irep$.

  For every $e\in E$, denote by $1_e$ the characteristic function of
$D_e\c\spec$.  A concrete description of $1_e$ may be given e.g.~by
  $$
  1_e(\phi) = \phi(e)
  \for \phi\in\spec.
  \eqmark OneeOfPhi
  $$
  Since $D_e$ is clopen we have that $1_e$ is continuous. 
Moreover $D_e$ is compact so $1_e\in C_c(\spec)\c
C_0(\spec)$.

\state Proposition
  \label IntroducePiu
  Let $\irep$ be a representation of $\S$ on a Hilbert space $H$.
Then there exists a unique *-representation $\pi_\irep$ of
$C_0(\spec)$ on $H$ such that $\pi_\irep(1_e)=\irep_e$, for every
$e\in E$.  In addition the pair $(\pi_\irep,\irep)$ is a covariant
representation of the system $(\act,\S,\spec)$

  \proof The Stone--Weierstrass theorem readily implies that the set
of all $1_e$'s span a dense subalgebra of $C_0(\spec)$, from where
uniqueness follows.  To prove existence, let $A$ be the closed
*-subalgebra of $B(H)$ generated by $\{\irep_e: e\in E\}$.  It is
immediate that $A$ is commutative, so let us denote the spectrum of
$A$ by $\hat A$\pilar{11pt}.
  Given $\psi\in\hat A$, observe that the map
  $$
  \phi: e\in E \mapsto \psi(\irep_e)\in\{0,1\}
  $$
  is a semicharacter of $E$ (it is nonzero because $\psi$ is
nonzero).  This allows us to define a map
  $$
  j: \psi\in \hat A \ \longmapsto \ \phi=\psi\circ\irep\in\spec,
  $$
  which is obviously continuous and injective.  If we temporarily (and
heretically) alter the definition of both $\hat A$ and $\spec$ by
dropping the requirement that characters (in the case of $\hat A$) and
semicharacters (in the case of $\spec$) be nonzero, then the map $j$
above will satisfy $j(0)=0$.  This means that, returning to the usual
(and sacrosanct) notion of spectrum, $j$ is a proper map.  It follows
that
  $$
  \pi_\irep:f\in C_0(\spec) \ \longmapsto \ f\circ j\in C_0(\hat A) =
A
  $$
  is a well defined surjective *-homomorphism.  Since $A\c B(H)$, we
may view $\pi_\irep$ as a representation of $A$ on $H$.
  Let us next prove that 
  $$
  \pi_\irep(1_e) = \irep_e.
  \subeqmark PiOneE
  $$
  To prove it observe that for every $\psi\in\hat A$ we have
  $$
  \psi(\pi_\irep(1_e)) = 
  \widehat{\pi_\irep(1_e)}(\psi)=
  1_e(j(\psi)) =
  1_e(\psi\circ \irep) =
  (\psi\circ \irep)(e) =
  \psi(\irep_e),
  $$
  proving \lcite{\PiOneE}.
  In order to prove \lcite{\DefineCovarRep.iii} let $s\in\S$ and $f\in
C_0(\spec)$.  Since the algebra generated by the $1_e$ is dense in
$C_0(\spec)$, we may assume that $f=1_e$, for some $e\in E$.  Denoting by
$\a_s(f)=f\circ\act_{s^*}$, for $f\in C_0(D_{\q s})$, as in
\lcite{\BunchOfNotations.i}, notice that
  \ $
  \a_s(1_e) = 1_e\circ \act_{s^*}
  $ \
  is the characteristic function of
  $$
  \{\phi\in\spec: \act_{s^*}(\phi)\in D_e\} =
  \{\phi\in\spec: \phi(s^*es)=1\} =
  D_{s^*es},
  $$
  that is,
  $\a_s(1_e) = 1_{s^*es}$.  Therefore 
  $$
  \irep_s\pi_\irep(1_e)\irep_{s^*} \={(\PiOneE)}
  \irep_s\irep_e\irep_{s^*} =
  \irep_{ses^*} =
  \pi_\irep(1_{s^*es}) =
  \pi_\irep(\a_s(1_e)).
  $$
  Addressing \lcite{\DefineCovarRep.iv} notice that the compacity of
$D_e$ implies that $C_0(D_e) = C(D_e)$ is a unital algebra with unit
$1_e$.  It follows that
  $$
  \overline{\pi_\irep\big(C_0(D_e)\big)(H)} =
  \pi_\irep(1_e)(H) \={(\PiOneE)}
  \irep_e(H),
  $$
  and the proof is complete.
  \proofend

We do not want to be restricted to studying only the action of $\S$ on
$\spec$.  In fact the most interesting intrinsic actions take place on
subsets of $\spec$.  But of course only invariant subsets matter.

\definition
  \label DefineRestriction
  We say that a subset $X\c\spec\,$ is \stress{invariant} if for every
$s\in\S$ one has that
  $$
  \act_s(D_{\q s}\cap X)\c X.
  $$
  In that case, for every $e\in E$ we denote
  $$
  D^X_e = D_e\cap X, 
  $$
  and for every $s\in\S$ we let 
  $$
  \act^X_s: D_{\q s}^X \to D_{\p s}^X
  $$
  be given by restricting $\act_s$.

It is then elementary to prove that the correspondence 
  $$
  \act^X : s\in\S \mapsto \I(X)
  \eqmark DefineSigmaX
  $$
  is an action of $\S$ on $X$.


These actions, for suitably chosen subsets $X\c\spec$, will dominate
our attention throughout this work.  It is therefore interesting that
we can sometimes guarantee that its groupoid of germs is Hausdorff:

\state Corollary
  Let $\S$ be an inverse semigroup which is a semilattice with respect
to its natural order (such as an $E^*$-unitary inverse semigroup).  If
$\E$ denotes the idempotent semilattice of $\S$, and if $X\c\spec$ is
a closed invariant subspace, then the groupoid of germs for the action
$\act^X$ of $\S$ on $X$, as defined in \lcite{\DefineRestriction}, is
a Hausdorff groupoid.

\proof 
  As pointed out shortly after \lcite{\OneeOfPhi}, we have that $D_e$
is a compact subset of $\spec$, for every $e\in E$.  Hence $D_e^X =
D_e\cap X$ is closed.  The statement then follows from
\lcite{\HausdorffGPG}.
  \proofend

The following result shows that invariant subsets may be found
underlying Hilbert space representations of $\S$.

\state Proposition
  \label KernelInvariant
  Given a representation $\irep$ of $\S$ on a Hilbert space $H$, write
the kernel of $\pi_\irep$ as $C_0(U)$, where $U$ is an open subset of
$\spec$.  Then $X:=\spec\setminus U$ is a closed invariant subset.

\proof Given $s\in\S$ and $\phi\in D_{\q s}\cap X$, suppose by
contradiction that $\act_s(\phi)\notin X$.  Then $\act_s(\phi)\in
U\cap D_{\p s}$, and by Urysohn's Theorem, there exists $f\in
C_0(U\cap D_{\p s})$ such that $f(\act_s(\phi))=1$.  Then
  $$
  0 = \irep_{s^*}\pi_\irep(f)\irep_s = \pi_\irep(\a_{s^*}(f)),
  $$
  which implies that $\a_{s^*}(f)\in C_0(U)$.  Since
  $$
  \a_{s^*}(f)(\phi) = f(\act_s(\phi)) = 1,
  $$
  we conclude that $\phi\in U$, but we have taken $\phi$ in $X$.  This
is a contradiction and hence $X$ is indeed invariant.
  \proofend

\definition Let $\irep$ be a representation of $\S$ on a Hilbert space
$H$.  We will say that $\irep$ is \stress{supported} on a given subset
$X\c\spec$ if the representation $\pi_\irep$ of \lcite{\IntroducePiu}
vanishes on $C_0(\spec\setminus X)$.

Fix for the time being a closed invariant set $X\c\spec$ and let 
  $$
  \G^X = \G(\act^X,\S,X),
  \eqmark DefineGx
  $$
  be the groupoid of germs associated to the system $(\act^X,\S,X)$.
Observe that for every $e\in E$ one has that
  $D^X_e$
  (defined in \lcite{\DefineRestriction})
  is a compact open subset of $X$.  Denoting by $1^X_e$ the
characteristic function of $D^X_e\c X$, we then have that $1^X_e\in
C_c(X)$.  Employing the notation introduced in
\lcite{\BunchOfNotations.v} we see that for every $s\in\S$,
  $$
  1^X_{\p s}\d_s \in C_c(\O_s)\c C_c(\G^X),
  $$
  where $\O_s$ was defined in \lcite{\BunchOfNotations.ii}.

\state Proposition 
  \label RepXIsRep
  Let $X\c\spec$ be a closed invariant set.  Then   the correspondence
  $$
  \irep^X: s\in\S\mapsto \iota\big(1^X_{\p s}\d_s\big) \in C^*(\G^X),
  $$
  (recall that $\iota$ was defined in \lcite{\NaturalMap})
  is a representation of $\S$ (where we imagine $C^*(\G)$ as an
operator algebra via any faithful *-representation) which is supported
on $X$.  In fact, the set $U$ referred to in \lcite{\KernelInvariant}
is precisely equal to $\spec\setminus X$.

\proof For simplicity in this proof we will occasionally drop the
superscripts ``$X$", as it will cause no confusion.  We will moreover
identify $C_c(\G^X)$ with its copy within $C^*(\G^X)$, hence dropping
``$\iota$" as well.  For $s,t\in\S$ we have by
\lcite{\PartialProducts.iii}
  $$
  \irep_s\irep_t=
  (1_{\p s}\d_s)\star (1_{\p t}\d_t) =
  \a_s\big(\a_{s^*}(1_{\p s})1_{\p t}\big) \d_{st} =
  \a_s(1_{\q s}1_{\p t}) \d_{st} =
  1_{st(st)^*} \d_{st},
  $$
  where the last step follows easily from
\lcite{\MovingDsInSpecAction}.  This proves
\lcite{\DefineRepSoDoCara.i}, and \lcite{\DefineRepSoDoCara.ii} may
easily be proved with the aid of \lcite{\PartialProducts.iv}.
  Thus $\irep^X$ is indeed a representation of $\S$ in $C^*(\G^X)$, but
we must still identify the set $U$ of \lcite{\KernelInvariant}.  The
relevant representation of $C_0(\spec)$ should really be denoted
  $\pi_{\irep^X}$,
  but we will simply denote it by $\pi$.
  By \lcite{\PiOneE} we have 
  $$
  \pi(1_e)=\irep_e = 1^X_{e}\d_e.
  \eqno{(\dagger)}
  $$
  Identifying the unit space of $\G^X$ with $X$ as in
\lcite{\GZeroIsX}, and hence identifying $C_0(X)$ as a subalgebra of
$C^*(\G^X)$, we may write ($\dagger$) as
  $$
  \pi(1_e)=1_e^X=1_e|_X,
  $$
  so we conclude that $\pi(f) = f|_X$, for all $f\in C_0(\spec)$.  The
kernel of $\pi$ is therefore seen to be $C_0(\spec\setminus X)$.
  \proofend

The above representation is {universal} in the following sense: 

\state Theorem
  \label FromRepISGtoGroupoid
  Let $\S$ be a countable inverse semigroup and let $\irep:\S\to
B(H)$ be a representation which is supported on a closed invariant
subset $X\c\spec$.  Then there exists a *-representation $\rho$ of
$C^*(\G^X)$ on $H$ such that
  $
  \rho\big(\iota(1^X_{\p s}\d_s)\big) =
  \irep_s,
  $
  for all $s\in\S$, and hence the diagram

  \medskip\beginmypicture
  \setcoordinatesystem units <0.0010truecm, -0.0013truecm> point at 0 0
  \pouquinho = 600
  \morph {-1500}{0000}{1500}{0000} \arwlabel{\irep}{0000}{-200}
  \put {$\S$} at -1500  0000
  \put {$B(H)$} at 1800  0000
  \put {$C^*(\G^X)$} at 0300  1500
  \morph {-1400}{0100}{0000}{1400} \arwlabel{\irep^X}{-1200}{0900} 
  \morph {0600}{1400}{1800}{0100} \arwlabel{\rho}{1400}{0900} 
  \endmypicture
  \bigskip

\noindent  commutes.  
  
\proof
  Let $\pi_\irep$ be as in \lcite{\IntroducePiu}.  Since $\pi_\irep$ vanishes on
$C_0(\spec\setminus X)$ we may factor $\pi_\irep$ through $C_0(X)$ obtaining a
representation $\pi$ of $C_0(X)$ on $H$ such that the diagram

  \bigskip\beginmypicture
  \setcoordinatesystem units <0.0010truecm, -0.0013truecm> point at 0 0
  \pouquinho = 600
  \morph {-1500}{0000}{1500}{0000} \arwlabel{\pi_\irep}{0000}{-200}
  \put {$C_0(\spec)$} at -1800  0000
  \put {$B(H)$} at 1800  0000
  \put {$C_0(X)$} at 0300  1500
  \morph {-1400}{0100}{0000}{1400} 
  \morph {0600}{1400}{1800}{0100} \arwlabel{\pi}{1400}{0900} 
  \endmypicture
  \bigskip

\noindent commutes, where the southeast arrow is given by restriction.
We then claim that $(\pi,\irep)$ is a covariant
representation of the system $(\act^X,\S,X)$.  In order to prove it
let $f\in C_0(X)$ and choose $g\in C_0(\spec)$ whose restriction to $X$
gives $f$.  Then for every $s\in\S$ we have
  $$
  \irep_s\pi(f)\irep_{s^*} =
  \irep_s\pi_\irep(g)\irep_{s^*} =
  \pi_\irep\big(\a_s(g)\big) =
  \pi\big(\a_s(g)|_X\big) =
  \pi\big(\a^X_s(g|_X)\big) =
  \pi\big(\a^X_s(f)\big),
  $$
  where we are denoting by $\a^X$ the action of $\S$ on $C_0(X)$
associated to $\act^X$, as in \lcite{\BunchOfNotations.i}.
This proves \lcite{\DefineCovarRep.iii}.  To check
\lcite{\DefineCovarRep.iv} observe that for every $e\in E$ we have
  $$
  \pi\big(C_0(D^X_e)\big)H =
  \pi(1_e^X)(H) = 
  \pi_\irep(1_e)(H) =
  \irep_e(H),
  $$
  concluding the proof that $(\pi,u)$ is covariant.  

We next wish to
apply Theorem \lcite{\UnivPropGpd} to this covariant representation, so we
must address the countability restrictions:  since $E$ is countable,
the product space $\{0,1\}^E$ is metrizable and hence second countable.
We are then given the green light to apply the said Theorem and
hence there exists a *-representation $\rho$ of $C^*(\G^X)$  on $H$
such that 
  $$
  \rho\big(\iota(f\d_s)\big)= \pi(f)\irep_s,
  $$
  for every $s\in\S$, and every $f\in C_c(D_{\p s})$. 
We then conclude that 
  $$
  \rho(\irep^X_s) =
  \rho\big(\iota(1^X_{\p s}\d_s)\big) =
  \pi(1^X_{\p s})\irep_s =
  \pi_\irep(1_{\p s})\irep_s \={(\PiOneE)}
  \irep_{\p s}\irep_s =
  \irep_{\p ss} = \irep_s.
  \proofend
  $$

The following is a main result:

\state Corollary
  \label TheBigCorollary
  Let $\S$ be a countable inverse semigroup and let $X$ be a closed
invariant subset of $\spec$.  Then there is a one-to-one
correspondence between representations $\irep$ of $\S$ supported on
$X$ and representations $\rho$ of the C*-algebra of the groupoid of
germs
  for the action $\act^X\!$ of\/ $\S$ on $X$.  If $\irep$ and $\rho$
correspond to each other then
  $
  \rho\big(\iota(1^X_{\p s}\d_s)\big) =
  \irep_s,
  $
  for all $s\in\S$.

\proof 
  Follows immediately from \lcite{\RepXIsRep} and
\lcite{\FromRepISGtoGroupoid}.
  \proofend

One should notice that any representation of $\S$ is supported on
$\spec$, so we obtain the following version of
\scite{\PatBook}{Theorem 4.4.1}:

\state Corollary
  If $\S$ is a countable inverse semigroup then there is a one-to-one
correspondence between representations of $\S$ and representations of
the C*-algebra of the groupoid of germs for the natural action $\act$
of $\S$ on $\spec$.


\section{Representations of semilattices}
  \label RepSemilatSec
  We have intentionally postponed until now a very delicate and subtle
conceptual problem.
  If $\S$ contains a zero element $0$, and $\irep$ is a Hilbert space
representation of $\S$,  is it not natural to expect that $\irep_0=0$?
  However, including our development so far, most treatments of
inverse semigroups completely ignore this issue.  In fact some of the
better known examples of inverse semigroup representations, such as
Wordingham's Theorem \scite{\PatBook}{Theorem 2.2.2}, do send zero to
a nonzero element!

The problem with zero is but the tip of an iceberg which we will now
explore.  The issue apparently only concerns the idempotent
semilattice of $\S$.

\cryout{So we will now fix an abstract semilattice%
  \fn{A \stress{semilattice} is by definition a partially ordered set
$E$ such that for every $x,y\in E$, there exists a maximum among the
elements which are smaller than $x$ and $y$.  Such an element is said
to be the \stress{infimum} of $x$ and $y$, and is denoted $x\wedge
y$.}
  $E$, which will always be assumed to contain a smallest element
$0$.}

\definition 
  \label ItsAndPerpInSL
  Given a partially ordered set $X$ with smallest element $0$, we
shall say that two elements $x$ and $y$ in $X$ \stress{intersect}, in
symbols
  \ $
  x\its y,
  $ \
  if there is a nonzero $z\in X$ such that $z\leq x,y$.  Otherwise we
will say that $x$ and $y$ are \stress{disjoint}, in symbols
  \ $
  x\perp y.
  $

If $\E$ is a semilattice it is easy to see that $x$ and $y$ intersect
if and only if $x\inf y \neq 0$.

\definition 
  \label DefineRepSemilat
  Let $\E$ be a semilattice and let 
  $\B=(\B,0,1,\wedge,\vee,\nega)$
  be a Boolean algebra.  By a \stress{representation} of $\E$ in $\B$ we
shall mean a map
  $
  \brep:\E\to\B,
  $
  such that \izitem
  \zitem $\brep(0) = 0$, and
  \zitem $\brep(x\inf y) = \brep(x)\inf\brep(y),$ for every
$x,y\in\E$.

\medskip Recall that a Boolean algebra $\B$ is also a semilattice
  under the standard order relation given by
  $$
  x\leq y \iff x= x\inf y
  \for x,y\in\B.
  $$

Fix for the time being a representation $\brep$ of a semilattice $\E$
in a Boolean algebra $\B$.  For every $x,y\in\E$, such that $x\leq y$,
one has that $x = x\inf y$, and hence
  $$
  \brep(x) = \brep(x\inf y) = \brep(x)\inf\brep(y),
  $$
  which means that $\brep(x)\leq \brep(y)$.  In other words,
$\brep$ preserves the respective order relations. On the other hand
if $x,y\in\E$ are such that $x\perp y$, one has that
$\brep(x)\perp\brep(y)$, which may also be expressed in $\B$ by
saying that
  $$
  \brep(x)\leq \nega \brep(y).
  $$
  More generally, if $X$ and $Y$ are finite subsets of $\E$, and one is
given an element $z\in\E$ such that $z\leq x$ for every $x\in X$, and
$z\perp y$ for every $y\in Y$, it follows that
  $$
  \brep(z) \leq 
  \bigwedge_{x\in X} \brep(x) \wedge 
  \bigwedge_{y\in Y} \nega{\brep(y)}.
  \eqmark CrucialInequality
  $$
  The set of all such $z$'s will acquire an increasing importance, so we
make the following:

\definition
  \label DefineSXY
  Given finite subsets $X,Y\c\E$, we shall denote by $\E^{X,Y}$
the subset of\/ $\E$ given by 
  $$
  \E^{X,Y} = \{z\in\E: z\leq x,\ \forall x\in X,\hbox{ and } z\perp
y,\ \forall y\in Y\}.
  $$

Notice that if $X$ is nonempty and
  $
  x_{min}=  \bigwedge_{x\in X}x,
  $
  one may replace $X$ in \lcite{\DefineSXY} by the singleton
$\{x_{min}\}$, without altering $\E^{X,Y}$.  However there does not
seem to be a similar way to replace $Y$ by a smaller set.

\definition 
  \label DefineCoverInSLat
  Given any subset $F$ of the semilattice $\E$, we shall say that a
subset $Z\c F$ is a \stress{cover for} $F$, if for every
nonzero $x\in F$, there exists $z\in Z$ such that $z\its x$.  If
$y\in\E$ and $Z$ is a cover for $F=\{x\in \E: x\leq y\}$, we will say
that $Z$ is a \stress{cover for} $y$.

The notion of covers is relevant to the introduction of the following
central concept (compare \scite{\infinoa}{1.3}):

\definition 
  \label DefLatTightRep  
  Let $\brep$ be a representation of the semilattice $\E$ in the
Boolean algebra $\B$.  We shall say that $\brep$ is \stress{tight} if
for every finite subsets $X,Y\c\E$, and for every finite cover
$Z$ for
  $\E^{X,Y}$,  one has that 
  $$
  \bigvee_{z\in Z}\brep(z) \geq  \bigwedge_{x\in X} \brep(x) \wedge
\bigwedge_{y\in Y} \nega{\brep(y)}.
  $$

Notice that the reverse inequality ``$\leq$" always holds by
\lcite{\CrucialInequality}.  Thus, when $\brep$ is tight, we actually
get an equality above.  We should also remark that in the absence of
any finite cover $Z$ for any $\E^{X,Y}$, every representation is
considered to be tight by default.

It should be stressed that the definition above is meant to include
situations in which $X$, $Y$, or $Z$ are empty, and in fact this will
often be employed in the sequel.  It might therefore be convenient to
reinforce the convention according to which the supremum of the empty
subset of a Boolean algebra is zero, and that its infimum is 1.

For example, if $X=Y=\emptyset$, then $\E^{X,Y}=\E$, and hence a cover
$Z$ for $\E^{X,Y}$ must contain quite a lot of elements.  If a
representation $\brep$ is tight then the supremum of $\brep(z)$ over
such a cover is required to coincide with 1.  This may be considered
as a \stress{nondegeneracy} condition for tight representations
(applicable only when $\E$ admits a finite cover).

In certain cases the verification of tightness may be simplified by
assuming that $X\neq\emptyset$:

\state  Lemma
  \label ChaeckWithNonEmptX
  Let $\brep:\E\to\B$ be a representation of the semilattice $\E$ in
the Boolean algebra $\B$ and suppose that $\brep$ is known to satisfy
the tightness condition \lcite{\DefLatTightRep} only when $X$ is
nonempty.  If moreover 
  \izitem 
  \zitem $\E$ contains a finite set $X$ such that $\bigvee_{x\in X}
\brep(x) = 1$, or
  \zitem $\E$ does not admit any finite cover,
  \medskip \noindent then $\brep$ satisfies \lcite{\DefLatTightRep}
in full, i.e., $\brep$ is tight.

\proof
  Our task is therefore to prove the tightness condition even when
$X=\emptyset$.  So, let $Y\c\E$ be a finite set and let $Z$ be a
finite cover for $\E^{\emptyset,Y}$.  Notice that for every $u\in\E$,
either $u\its y$, for some $y\in Y$, or $u\in \E^{\emptyset,Y}$, in
which case $u\its z$, for some $z\in Z$.  Therefore $Y\cup Z$ is a
finite cover for $\E$.  Under hypothesis (ii) this is impossible,
meaning that there are no finite covers for $\E^{\emptyset,Y}$, 
there is nothing to be done.  We therefore
assume (i), and we must show that
  $$
  \bigvee_{z\in Z} \brep(z) \geq \bigwedge_{y\in Y} \nega\brep(y).
  \subeqmark GoalOfTighWithoutX
  $$
  Let $X$ be as in the statement.  We claim that for each $x\in X$,
the set
  $x\inf Z := \{x\inf z: z\in Z\}$ is a cover for $\E^{\{x\},Y}$.  In
fact, given a nonzero
  $$
  w\in\E^{\{x\},Y}\c\E^{\emptyset,Y},
  $$
  there exists some $z\in Z$ such that $z\its w$.  Since $w\leq x$, we
have
  $$
  w\inf  x\inf z  =   w\inf z  \neq 0,
  $$
  so $w\its (x\inf z)$, concluding the proof of our claim.  By
hypothesis $\brep$ satisfies the tightness condition with respect to
the cover $x\inf Z$ \ for $\E^{\{x\},Y}$, and hence
  $$
  \bigvee_{z\in Z}\brep(x\inf z) \geq
  \brep(x) \wedge \bigwedge_{y\in Y} \nega\brep(y).
  \subeqmark TightUnderX
  $$
  We therefore have
  $$
  \bigvee_{z\in Z}\brep(z) \={(i)}
  \bigvee_{z\in Z}\Big(\bigvee_{x\in X}\brep(x)\Big)\wedge\brep(z) =
  \bigvee_{x\in X} \bigvee_{z\in Z}\brep(x\inf z) 
  \buildrel (\TightUnderX)\over \geq
  \bigvee_{x\in X} \Big(\brep(x) \wedge \bigwedge_{y\in Y}
\nega\brep(y)\Big) \$=
  \Big(\bigvee_{x\in X} \brep(x)\Big) \wedge \bigwedge_{y\in Y}
\nega\brep(y) =
  \bigwedge_{y\in Y} \nega\brep(y),
  $$
  proving \lcite{\GoalOfTighWithoutX}.
  \proofend

The following alternative characterization of tightness is apparently
even weaker than the above:

\state Proposition
  \label AltLatTightRep
  Let $\brep$ be a representation of the semilattice $\E$ in the
Boolean algebra $\B$, satisfying either (i) or (ii) of
\lcite{\ChaeckWithNonEmptX}.  Then $\brep$ is tight if and only if
for every $x\in\E$ and for every finite 
  cover $Z$ for $x$, one has that
  $$
  \bigvee_{z\in Z}\brep(z) \geq \brep(x).
  $$

\proof
  The only if part is immediate since $\{u\in\E: u\leq x\}=
\E^{\{x\},\emptyset}$.  To prove the converse implication let
$X,Y\c\E$ be finite subsets and let $Z$ be a cover for $\E^{X,Y}$.
Using \lcite{\ChaeckWithNonEmptX} we may assume that $X$ is nonempty,
so let $x_{min}= \bigwedge_{x\in X}x$. We claim that $Y\cup Z$ is a
cover for $\E^{\{x_{min}\},\emptyset}$.  In order to prove it pick
$u\leq x_{min}$.  Then clearly $u\leq x$, for every $x\in X$.

Suppose first that $u\notin \E^{X,Y}$.  Then $u$ necessarily fails to
be disjoint from some $y\in Y$, meaning that $x\its y$, and thus
proving that $u$ intersects some element of $Y\cup Z$.  On the other
hand, if $u\in \E^{X,Y}$, then our assumption guarantees that there
exists some element $z$ in $Z$, and hence also in $Y\cup Z$, which
intersects $x$.  This proves our claim, and so the hypothesis gives
  $$
  \brep(x_{min}) \leq \bigvee_{u\in Y\cup Z}\brep(u),
  $$
  and hence also 
  $$
  \brep(x_{min}) \wedge 
  \Big(\bigwedge_{y\in Y} \nega{\brep(y)}\Big) \leq
  \Big(\bigvee_{u\in Y\cup Z}\brep(u)\Big) \wedge
  \Big(\bigwedge_{y\in Y} \nega{\brep(y)}\Big) =
  \bigvee_{u\in Y\cup Z}\Big(\brep(u)\wedge
  \bigwedge_{y\in Y} \nega{\brep(y)}\Big).
  $$
  Referring to the term $\brep(u)\wedge \bigwedge_{y\in Y}
\nega{\brep(y)}$, appearing above, notice that it is zero for every
$u\in Y$.  In case $u\in Z$, then because $Z\c \E^{X,Y}\pilar{11pt}$,
we see that $\brep(u)\leq\nega{\brep(y)}$, for all $y\in Y$, and
hence the alluded term coincides with $\brep(u)$.  The right-hand
side of the expression displayed above thus becomes simply
  $
  \bigvee_{u\in Z}\brep(u),
  $
  and since 
  $\brep(x_{min}) = \bigwedge_{x\in X}\brep(x)$, the left-hand side
is
  $$
  \Big(\bigwedge_{x\in X} \brep(x)\Big) \wedge \Big(\bigwedge_{y\in
Y} \nega{\brep(y)}\Big).
  \proofend
  $$

When $\E$ happens to be a Boolean algebra there is a very
elementary characterization of tight representations:

\state  Proposition
  \label BooleToBoole
  Suppose that $\E$ is a semilattice admitting the structure of a
Boolean algebra which induces the same order relation as that of $\E$,
and let 
  \ $
  \brep:\E\to\B 
  $ \
  be a representation of $\E$ in some Boolean algebra $\B$.  Then
$\brep$ is tight if and only if it is a Boolean algebra homomorphism.

\proof
  Supposing that $\brep$ is tight, notice that $\{1\}$ is a cover for
$\E^{\emptyset,\{0\}}$, so
  $$
  \brep(1) = \nega \brep(0) = \nega 0 = 1.
  $$
  Given $x\in\E$ notice that $\{\nega x\}$ is a cover for
$\E^{\emptyset,\{x\}}$, therefore 
  $$
  \brep(\nega x) = \nega\brep(x).
  $$
Since 
  $
  x\vee y =   \nega(\nega x\inf \nega y),
  $
  for all $x,y\in\E$,  we may easily prove that $\brep(x\vee
y)=\brep(x)\vee \brep(y)$.  Thus $\brep$ is a Boolean algebra
homomorphism, as required.

In order to prove the converse implication let $X,Y\c\E$ be
finite sets and let $Z$ be a finite cover for $\E^{X,Y}$.  Let
  $$
  z_0 = \bigvee_{z\in Z} z, \quad
  x_0 = \bigwedge_{x\in X} x \and
  \bar y_0 = \bigwedge_{y\in Y} \nega y.
  $$
  It is obvious that $z_0\leq x_0\inf \bar y_0$, and we claim that in
fact $z_0= x_0\inf \bar y_0$.  We will prove it by checking that
  $$
  \nega z_0\inf x_0\inf \bar y_0 = 0.  
  $$
  Let $u = \nega z_0\inf x_0\inf \bar y_0$, and notice that the fact that
$u\leq x_0\inf \bar y_0$ implies that $u\in\E^{X,Y}$. Arguing by
contradiction, and hence supposing that $u$ is nonzero, we deduce that
$u\its z$, for some $z\in Z$, but this contradicts the fact that $u\leq
\nega z_0$.  This proves our claim so, assuming that $\brep$ is a
Boolean algebra homomorphism, we have
  $$
  \bigvee_{z\in Z} \brep(z) =
  \brep\Big( \bigvee_{z\in Z} z\Big) =
  \brep(z_0) =
  \brep(x_0\inf \bar y_0) =
  \bigwedge_{x\in X} \brep(x) \inf 
  \bigwedge_{y\in Y} \nega \brep(y),
  $$
  showing that $\brep$ is tight.
  \proofend

Not all semilattices admit tight injective representations.  In order to
study this issue in detail it is convenient to introduce the following:

\definition 
  \label DefineDense
  Let $\E$ be a semilattice and let $x,y\in\E$ be such that $y\leq x$.
We shall say that $y$ is \stress{dense} in $x$ if there is no nonzero
$z\in\E$ such that $z\perp y$ and $z\leq x$.  Equivalently, if
$\E^{\{x\},\{y\}}=\{0\}$.

Obviously each $x\in\E$ is dense in itself but it is conceivable that
some $y\neq x$ is dense in $x$.  For a concrete example notice that in
the semilattice $\E=\{0,{1\over 2},1\}$, where $0\leq {1\over
2}\leq1$, one has that ${1\over 2}$ is dense in $1$.

In the general case, whenever $y$ is dense in $x$ we have that
  $
  \E^{\{x\},\{y\}}=\{0\},
  $
  and hence the empty set is a cover for 
  $\E^{\{x\},\{y\}}$.  Therefore for every tight representation $\brep$
of $\E$ one has that
  $$
  0=\brep(x) \wedge \nega{\brep(y)},
  $$
  which means that $\brep(x) \leq \brep(y)$.  Since the opposite
inequality also holds, we have that $\brep(x)=\brep(y)$.  Thus no
tight representation of $\E$ can possibly separate $x$ and $y$.
The reader is referred to \cite{\ExelAlgebra} for a thorough study of this
and related problems.
For future reference we record this conclusion in the next:

\state Proposition
\label ThouShallNotSeparate
If $y\leq x$ are elements in the semilattice $\E$, such that $y$
is dense in $x$, then $\brep(y)=\brep(x)$ for every tight
representation $\brep$ of\/ $\E$.

We will have a lot more to say about tight representations in the
following sections.


\section{Filters and characters}
   As in the previous section we fix a semilattice $E$ with smallest
element $0$.
  A fundamental tool for the study of tight representations of
$E$ is the notion of filters, which we shall discuss in this
section.

\definition
  \label DefineFilter
  Let $X$ be any partially ordered set with minimum element $0$.
  A \stress{filter} in $X$ is a nonempty
subset $\xi\c X$, such that
  \izitem
  \zitem $0\notin\xi$,
  \zitem if $x\in\xi$ and $y\geq x$, then $y\in\xi$,
  \zitem if $x,y\in\xi$, there exists $z\in\xi$, such that $x,y\geq z$.
  \medskip \noindent An \stress{ultra-filter} is a filter which is not
properly contained in any filter. 

Given a partially ordered set $X$ and any nonzero element $x\in X$, it is
elementary to prove that
  $$
  \xi=\{y\in X: y\geq x\}
  $$
  is a filter containing $x$.  By Zorn's Lemma there exists an
ultra-filter containing $\xi$,   thus every nonzero element in $X$ belongs
to some ultra-filter.

When $\E$ is a semilattice, given the existence of $x\inf y$ for every
$x,y\in\E$, condition \lcite{\DefineFilter.iii} may be replaced by 
  $$
  x,y\in\xi \  \imply  \ x\inf y\in\xi.
  \eqmark NewCondForFilInSL
  $$

The following is an important fact about filters in semilattices which
also benefits from the existence of $x\inf y$.

\state Lemma
  \label UltrafilterCriterium
  Let $\E$ be a semilattice and let $\xi$ be a filter in $\E$.  Then
$\xi$ is an ultra-filter if and only if $\xi$ contains every element
$y\in\E$ such that $y\its x$, for every $x\in\xi$.

\proof
  In order to prove the ``if" part let $\eta$ be a filter such that
$\xi\c\eta$.  Given $y\in\eta$ one has that for every
$x\in\xi$, both $y$ and $x$ lie in $\eta$, and hence 
\lcite{\NewCondForFilInSL} implies that
  $y\inf x\in\eta$, so $y\inf x\neq 0$, and hence $y\its x$.  By
hypothesis $y\in\xi$, proving that $\eta=\xi$, and hence that $\xi$ is
an ultra-filter.

Conversely let $\xi$ be an ultra-filter and suppose that $y\in\E$ is
such that $y\its x$, for every $x\in\xi$.  Defining
  $$
  \eta = \{ u\in\E: u \geq y\inf x, \hbox{ for some } x\in \xi\},
  $$
  we claim that $\eta$ is a filter.
  By hypothesis $0\notin\eta$. 
  Also if $u_1,u_2\in\eta$, choose for every $i=1,2$ some $x_i\in\xi$
such that $u_i\geq y\inf x_i$.  Then
  $$
  u_1\inf u_2 \geq (y\inf x_1)\inf(y\inf x_2) = y\inf (x_1\inf x_2),
  $$
  so $u\in\eta$.  Given that \lcite{\DefineFilter.ii} is obvious we
see that $\eta$ is indeed a filter, as claimed.  Noticing that
$\xi\c\eta$ we have that $\eta=\xi$, because $\xi$ is an
ultra-filter.  Since $y\in\eta$, we deduce that $y\in\xi$.
  \proofend

The study of representations of our semilattice $\E$
in the most elementary Boolean algebra of all, namely $\{0,1\}$, leads
us to the following specialization of the notion of semicharacters:

\definition
  \label DefineCharacter
  By a \stress{character} of $\E$ we shall mean any nonzero
representation of $\E$ in the Boolean algebra
  $
  \{0,1\}.
  $
  The set of all characters will be denoted by $\specz$.

Thus, a character is nothing but a semicharacter which vanishes at
$0$.  Perhaps the widespread use of the term \stress{semicharacter} is
motivated by the fact that it shares prefix with the term
\stress{semilattice}.  If this is really the case then our choice of
the term \stress{character} may not be such a good idea but alas, we
cannot think of a better term.

It is easy to see that $\specz$ is a closed subset of $\spec$, and
hence that $\specz$ is locally compact.

Given a character $\phi$, observe that
  $$
  \xi_\phi  = \{x\in\E: \phi(x) = 1\},
  \eqmark XiPhi
  $$
  is a filter in $\E$ (it is nonempty because $\phi$ is assumed not to
be identically zero).
   Conversely, given a filter $\xi$, define for every $x\in\E$,
  $$
  \phi_\xi(x) = \left\{\matrix{
  1, & \hbox{ if } x\in\xi,\hfill\cr
  \pilar{12pt}
  0, & \hbox{ otherwise.}}\right.
  \eqmark PhiXi
  $$
  It is then easy to see that $\phi_\xi$ is a character.  Therefore we
see that \lcite{\XiPhi} and \lcite{\PhiXi} give one-to-one
correspondences between $\specz$ and the set of all filters.

\state  Proposition
  \label UltraThenTight
  If $\xi$ is an ultra-filter then $\phi_\xi$ is a tight representation
of\/ $\E$ in $\{0,1\}$.

  \proof
  Let $X,Y\subset \E$ be finite subsets and let $Z$ be a cover for
$\E^{X,Y}$.  In order to prove that
  $$
  \bigvee_{z\in Z}\phi(z) \geq 
  \prod_{x\in X} \phi(x) \prod_{y\in Y} (1-\phi(y)),
  $$
  it is enough to show that if the right-hand side equals 1, then so
do the left-hand side.  This is to say that if $x\in\xi$ for every
$x\in X$, and $y\notin\xi$ for every $y\in Y$, then there is some
$z\in Z$, such that $z\in\xi$. 

By \lcite{\UltrafilterCriterium}, for each $y\in Y$ there exists some
$x_y\in\xi$ such that $y\perp x_y$.
Supposing by contradiction that
$Z\cap\xi = \emptyset$, then for every $z\in Z$ there exists, again by
\lcite{\UltrafilterCriterium}, some $x_z\in\xi$, such that $z\perp
x_z$.  Set 
  $$\hbox{$
  w =
  \bigwedge\limits_{x\in X} x \wedge
  \bigwedge\limits_{y\in Y} x_y \wedge
  \bigwedge\limits_{z\in Z} x_z.
  $}$$
  Since $w\in\xi$ we have that $w\neq 0$.  Obviously $w\leq x$ for
every $x\in X$, and $w\perp y$ for every $y\in Y$, and hence
$w\in\E^{X,Y}$.  Since $Z$ is a cover there exists some $z_1\in Z$
such that $w\its z_1$.  However, since $w\leq x_{z_1}\perp z_1$, we
have that $w\perp z_1$, a contradiction.
  \proofend

\definition
  \label AltSpectrums
  We shall denote by $\speci$ the set of all characters
$\phi\in\specz$ such that $\xi_\phi$ is an ultra-filter.  Also we will
denote by $\spect$ the set of all tight characters.

Employing the terminology just introduced we may rephrase
\lcite{\UltraThenTight} by saying that $\speci\c\spect$.
The following main result further describes the relationship
between $\speci$ and $\spect$.

\state Theorem
  \label ClosureOfUltrafilters
  Let $\E$ be a semilattice with smallest element $0$, and let
$\speci$ and $\spect$ be as defined in \lcite{\AltSpectrums}.  Then
the closure of $\speci$ in $\specz$ coincides with $\spect$.

\proof Since the condition for any given  $\phi$ in $\specz$ to belong to
$\spect$ is given by equations 
it is easy to prove that $\spect$ is closed within $\specz$,
and since $\speci\c\spect$ by \lcite{\UltraThenTight}, we
deduce that
  $$
  \overline{\speci}\c\spect.
  $$
  To prove the reverse inclusion let us be given $\phi\in\spect$.  We
must therefore show that $\phi$ can be arbitrarily approximated by
elements from $\speci$.  Let $U$ be a neighborhood of $\phi$
within $\specz$.  By definition of the product topology, $U$ contains
a neighborhood of $\phi$ of the form
  $$
  V = V_{X,Y} = 
  \{\psi\in\specz: 
  \psi(x) = 1, \hbox{ for all } x\in X, \hbox{ and }
  \psi(y) = 0, \hbox{ for all } y\in Y\},
  $$
  where $X$ and $Y$ are finite subsets of $\E$.  
We next claim that $\E^{X,Y}\neq \{0\}$.  In order to prove this
suppose the contrary, and hence $Z=\emptyset$ is a cover for
$\E^{X,Y}$.  Since $\phi$ is tight we conclude that 
  $$
  0 = \bigvee_{z\in Z}\phi(z) =
  \prod_{x\in X} \phi(x) \prod_{y\in Y} (1-\phi(y)).
  $$
  However, since $\phi$ is supposed to be in $V$, we have that 
$\phi(x)=1$ for all $x\in X$, and $\phi(y)=0$ for all $y\in Y$, which
means that the right-hand side of the expression displayed above
equals 1.  This is a contradiction and hence our claim is proved.

We are therefore allowed to choose a nonzero $z\in \E^{X,Y}$, and
further to pick
an ultra-filter $\xi$ such that $z\in\xi$.  Observe that
$\phi_\xi\in\speci$, and the proof will be concluded once we
show that $\phi_\xi\in U$.

For every $x\in X$ and $y\in Y$, we have that $z\leq x$ and $z\perp
y$, hence $x\in\xi$ and $y\notin \xi$.  This entails $\phi_\xi(x)=1$
and $\phi_\xi(y)=0$, so $\phi_\xi\in V\c U$, as required.
  \proofend

Before we close this section let us discuss the issue of tight filters
in the idempotent semilattice of an inverse semigroup.  We
specifically want to prove that the correspondence described by
\lcite{\MovingSemicharacters.ii} preserves tight characters.  For this
we need an auxiliary result:

\state Lemma
  \label MovingCovers
  Let $\S$ be an inverse semigroup with zero and let $E$ be the
idempotent semilattice of $\S$.
  Given finite subsets $X$ and $Y$ of $E$, with $X$ nonempty, let $Z$
be a finite cover for $E^{X,Y}$.  Then for every $s\in\S$ one has that
$sZs^*$ is a cover for $\pilar{14pt} E^{sXs^*,sYs^*}$.

\proof
  Let $w$ be a nonzero element of $E$ such that $w\leq sxs^*$ for
every $x\in X$, and $w\perp sys^*$ for every $y\in Y$.  Then
  $$
  (s^*ws)y = 
  s^*wss^*sy = 
  s^*wsys^*s = 0,
  $$
  so $s^*ws\perp y$, for every $y\in Y$.  For every  $x\in X$ we have
that 
  $$
  (s^*ws) x = 
  s^*wss^*s x = 
  s^*ws x s^*s =
  s^*ws,
  $$
  so 
  $s^*ws\leq x$. This shows that $s^*ws\in E^{X,Y}$, and we claim
that $s^*ws\neq 0$.  For this choose $x\in X$ (allowed because $X$ is
nonempty) and observe that
  $w\leq sxs^*\leq ss^*$.  So
  $$
  0\neq w = ss^* w ss^*,
  $$
  which implies our claim.
  By hypothesis there exists some $z\in Z$ such that
$s^*ws\its z$.  Noticing that 
  $$
  0\neq 
  (s^*ws)z =
  s^*wss^*sz =
  s^*wszs^*s,
  $$
  we deduce that $wszs^*\neq 0$, so $w\its szs^*$. \proofend

The promised preservation of tightness is in order:
 
\state Proposition
  \label TightInvariance
  Let $\S$ be an inverse semigroup with zero and let $E$ be the
idempotent semilattice of $\S$.
  Given $s\in\S$ and a tight character $\phi$ on $E$ such that
$\phi(\q s)=1$, one has that the character $\act_s(\phi)$ defined in
\lcite{\MovingSemicharacters.ii} is also tight.

  \proof
  In view of the requirement that $X$ be nonempty in
\lcite{\MovingCovers} we will use \lcite{\ChaeckWithNonEmptX} for the
characterization of tight characters.  We may do so for $\act_s(\phi)$
because $\act_s(\phi)(\p s)=1$.

So let $X$ and $Y$ be finite subsets of $E$, with $X$ nonempty, and
let $Z$ be a cover for $E^{X,Y}$. Then
  $$
  \bigvee_{z\in Z}\act_s(\phi)(z) =  
  \bigvee_{z\in Z}\phi(s^*zs) =  
  \bigvee_{z'\in s^*Zs}\phi(z') =  
  \prod_{x'\in s^*Xs} \phi(x') 
  \prod_{y'\in s^*Ys} 1-\phi(y') \$=
  \prod_{x\in X} \phi(s^*xs) 
  \prod_{y\in Y} 1-\phi(s^*ys) =
  \prod_{x\in X} \act_s(\phi)(x) 
  \prod_{y\in Y} 1-\act_s(\phi)(y),
  $$
  where we have used \lcite{\MovingCovers} and the hypothesis that
$\phi$ is tight in walking through the third equal sign above.  This
concludes the proof.
  \proofend

If the content of this work is to be subsumed in a single idea, than
that idea is that the most natural intrinsic action of $\S$ on a
topological space is the restriction of the action $\act$ to $\spect$,
as defined by \lcite{\DefineSigmaX}.  In the following sections we
hope to convince the reader of its relevance.

\section{Tight representations of inverse semigroups}
  \label TightRepISGSect
  Throughout this section we will fix an inverse semigroup $\S$ with
$0$.  Suppose we are given a representation $\irep$ of $\S$ on a
Hilbert space $H$ and denote by $A$ the closed unital *-subalgebra of
$B(H)$ generated by the identity operator and $\{\irep_e: e\in
E(\S)\}$.  Since $A$ is abelian we see that the set
  $$
  \B_A = \{e\in A: e^2=e\}
  $$
  is a Boolean algebra relative to the operations
  $$
  e\wedge f = ef, \quad
  e\vee f = e+f-ef \and
  \nega e = 1-e,
  $$
  for all $e,f\in \B_A$.
Provided we assume that $\irep_0=0$, it is clear that the restriction
of $\irep$ to $E(\S)$ is a representation of $E(\S)$ in $\B_A$, in the
sense of Definition \lcite{\DefineRepSemilat}.

\definition 
  \label DefTightRepISG
  A representation $\irep$ of $\S$ on a Hilbert space $H$ is said to be
\stress{tight} if the restriction of $\irep$ to $E(\S)$ is a tight
representation of $E(\S)$ in the Boolean algebra $\B_A$, in the sense
of \lcite{\DefLatTightRep}.

Notice that, at the very least, tight representations are required to
satisfy $\irep_0=0$.

\state  Theorem 
  \label CharactTightReps
  A representation $\irep$ of $\S$ on a Hilbert
space $H$ is  tight if and only if it is supported in $\spect$.  

\proof
  Let $\pi_\irep$ be the *-representation of $C_0(\spec)$ on $H$ given by
\lcite{\IntroducePiu},  and write ${\rm Ker}(\pi_u)=C_0(U)$, for a
suitable open subset $U\c\spec$.  Fix finite subsets $X,Y\c\E$ and a finite
cover $Z$ for $E(\S)^{X,Y}$.  
  The condition for tightness of $\irep$  is that 
  $$
  \bigvee_{z\in Z} \irep_z =
  \prod_{x\in X} \irep_x 
  \prod_{y\in Y} (1- \irep_y),
  \eqno{(\dagger)}
  $$
  which, in view of \lcite{\PiOneE}, is equivalent to 
  $$
  \bigvee_{z\in Z} \pi_\irep(1_z) =
  \prod_{x\in X} \pi_\irep(1_x) 
  \prod_{y\in Y} \big(1- \pi_\irep(1_y)\big),
  $$
  or to
  $$
  \bigvee_{z\in Z} 1_z -
  \prod_{x\in X} 1_x
  \prod_{y\in Y} (1- 1_y) \in C_0(U).
  $$
  If $f = f_{X,Y,Z}$ is the function on the left-hand side of the
expression displayed above then to say that $f\in C_0(U)$ means that
$f(\phi)=0$, for every $\phi\notin U$.

Using \lcite{\OneeOfPhi} notice that for every $\phi\in\spec$, to say
that $f(\phi)=0$, is the same as saying that
  $$
  \bigvee_{z\in Z} \phi(z) =
  \prod_{x\in X} \phi(x)
  \prod_{y\in Y} (1- \phi(y)).
  \eqno{(\ddagger)}
  $$
  Summarizing, $\irep$ is tight if and only if for every $X$, $Y$ and
$Z$, as above, one has that ${(\ddagger)}$ holds for every
$\phi\in\spec\setminus U$.  But this is precisely expressing that 
  $$
  \spec\setminus U\c\spect,
  $$
  which is equivalent to 
  $
  \spec\setminus \spect\c U,
  $
  or to saying that $\pi_\irep$ vanishes on $C_0(\spec\setminus
\spect)$.  The last condition means, by definition, that $\irep$ is
supported in $\spect$.
  \proofend

The following result largely subsumes our main point so far:

\state Theorem
  \label EquivTightRepISGandGPD
  Let $\S$ be a countable inverse semigroup with zero and let $\Gt$ be
the groupoid of germs associated to the restriction of the action
$\act$ of \lcite{\MovingSemicharacters.iv} to the closed invariant
space $\spect\c\spec$.  Then there is a one-to-one correspondence
between tight Hilbert space representations of $\S$ and
*-representations of $C^*(\Gt)$.  An explicit form of this
correspondence is given by the formula at the end of
\lcite{\TheBigCorollary}.

\proof Follows immediately from \lcite{\TheBigCorollary} and
\lcite{\CharactTightReps}.
  \proofend


\section{The inverse semigroup associated to a semigroupoid}
  \label ISGFromGPDSec
  With this section we start to discuss an application of our methods
to semigroupoid C*-algebras, as defined in \cite{\ExSemiGpdAlg}.
Our task here will be to construct an inverse semigroup
$\SM$ from a given semigroupoid $\M$.  

We begin by recalling a few basic concepts from the theory of
semigroupoids.  See \cite{\ExSemiGpdAlg} for more details.
  A \stress{semigroupoid} is a triple $(\M,\Mt,\ \cdot\ )$ such that
$\M$ is a set, $\Mt$ is a subset of $\M\times\M$, and
  $$
  \cdot\ : \Mt \to \M
  $$
  is an operation which is associative in the following sense:  if
$f,g,h\in\M$ are such that either
  \bitem $(f,g)\in\Mt$ and $(g,h)\in\Mt$, or
  \bitem $(f,g)\in\Mt$ and $(f g,h)\in\Mt$, or
  \bitem $(g,h)\in\Mt$ and $(f,\ g h)\in\Mt$,

  \medskip\noindent 
  then all of $(f,g)$, $(g,h)$, $(f g,h)$ and $(f,g h)$ lie in
$\Mt$, and
  $$
  (f g) h = f (g h).
  $$
  Moreover, for every $f\in\M$, we will let
  $$
  \D f = \left\{g\in\M: (f,g)\in\Mt\right\}.
  $$

\cryout{From now on we fix a semigroupoid $\M$.}

\noindent If $f,g\in\M$ we will say that $f$ \stress{divides} $g$, or
that $g$ is a \stress{multiple} of $f$, in symbols $f\dil g$, if
either
  \bitem $f=g$, or
  \bitem  there exists $h\in\M$ such that
  $
  fh=g.
  $

\bigskip
We recall from \cite{\ExSemiGpdAlg} that division is reflexive, transitive and
invariant under multiplication on the left.

A useful artifice is to introduce a unit for $\M$, that is, pick some
element in the universe outside $\M$, call it $1$, set
$\Mu=\M\mathop{\dot\cup}\ \{1\}$, and for every $f\in\Mu$ put
  $$ 
  1f = f1 = f.
  $$
  Then, whenever $f\dil g$, regardless of whether $f=g$ or not, there
always exists $x\in\Mu$ such that $g=fx$.
 
We will find it useful to extend the definition of $\D f$, for
$f\in\Mu$, by putting
  $$
  \D 1 = \M.
  $$ 
  Nonetheless, even if $f1$ is a meaningful product for every
$f\in\M$, we \stress{will not} include $1$ in $\D f$.  In the few
occasions that we need to refer to the set of \stress{all} elements
$x$ in $\Mu$ for which $fx$ makes sense we shall use $\D f\cup\{1\}$.

It is interesting to notice that, as a consequence of the associative
axiom, for every $f\in\Mu$, and $g\in\D f$, one has 
  $$
  \D{fg} = \D g.
  \eqmark DfgDg
  $$
  Note that condition above does not allow for $g=1$, since $1$ is
never in $\D f$.  Besides, if $g=1$ then the above equality will most
likely fail.
  It is also easy to see that if $g\in\M$, and $h\in \D g\cup\{1\}$,
then
  $$
  g\in\D f\iff gh\in\D f,
  \eqmark DfgDgDois
  $$
  for every $f\in\Mu$.

Recall from \scite{\ExSemiGpdAlg}{Section 3} that a \stress{spring} is
an element $f\in\M$ such that
  $$
  \D f=\emptyset.
  $$
  If $f$ is a spring one is therefore not allowed to right-multiply it
by any element, that is, $fg$ is never a legal multiplication, unless
$g=1$.  In some key places below we will suppose that $\M$ has no
springs.

We should be aware that $\Mu$ is \stress{not} a semigroupoid.
Otherwise, since $f1$ and $1g$ are meaningful products, the
associativity axiom would imply that $(f1)g$ is also a meaningful
product, but this is clearly not always the case.  Nevertheless it is
interesting to understand precisely which one of the three clauses of
the associativity property is responsible for this problem.  As
already observed, the first clause does fail irremediably when $g=1$.
However it is easy to see that all other clauses do generalize to
$\Mu$.  This is quite useful, since when we are developing a
computation, having arrived at an expression of the form $(fg)h$, and
therefore having already checked that all products involved are
meaningful, we most often want to proceed by writing
  $$
  \ldots = (fg)h = f(gh),
  $$
  and this is fortunately meaningful and correct for all
$f,g,h\in\Mu$, because it does not rely on the delicate first clause
of associativity.

Given $f,g\in\M$, we say that $f$ and $g$ \stress{intersect} if they
admit a \stress{common multiple}, that is, an element $m\in\M$ such
that $f\dil m$ and $g\dil m$.  Otherwise we will say that $f$ and $g$
are \stress{disjoint}.  We will write
  $$
  f\its g 
  \eqmark DefineItsSGPD
  $$
  when $f$ and $g$ intersect, and
  $$
  f\disj g 
  \eqmark DefineDisjSGPD
  $$
  when $f$ and $g$ are disjoint.  Incidentally this notation employs
the same symbols ``$\its$" and ``$\disj$", defined in
\lcite{\ItsAndPerpInSL} in connection to semilattices, with different
(although deeply related)
meanings, and we will rely on the context to determine the correct
interpretation of our notation.  Employing the unitization $\Mu$
notice that $f\its g$ if and only if there are $x,y\in\Mu$ such that
$fx=gy$.

\definition 
  \label DefineMonic
  We shall say that an element $f\in\Mu$ is \stress{monic} if for every
$g,h\in\Mu$ we have
  $$
  fg=fh \imply g=h.
  $$

  Observe that the above includes the implication
  \ $
  fg=f \ \imply \ g=1.
  $
  Obviously 1 is monic.

  \definition 
  \label DefineLCM
  Let $f,g\in\M$ be such that $f\its g$.  We shall say that an element
$m\in\M$ is a \stress{least common multiple} of $f$ and $g$, if $m$ is
a common multiple of $f$ and $g$ and for every other common multiple
$h$, one has that $m\dil h$.

From now on we shall assume the following:

\sysstate{Standing Hypothesis}{\rm}{\label StandingHyp $\M$ is a
semigroupoid in which every element is monic, and moreover every
intersecting pair of elements admits a least common multiple.}

Observe that if $f,g,h\in\Mu$ then
  $$
  f\dil g \hbox{ \ and \ } g\dil f \ \imply \ f=g.
  \eqmark NoEquiv
  $$
  In fact, writing $f=gx$, and $g=fy$, for $x,y\in\Mu$, we deduce that
  $$
  g=fy = (gx)y  = g(xy),
  $$
  which implies that $xy=1$, but this can only happen if $x=y=1$, and
hence $f=g$.

If $f$ and $g$ are intersecting elements in $\M$ and if $m_1$ and
$m_2$ are both least common multiples of $f$ and $g$, then $m_1\dil
m_2$ and $m_2\dil m_1$, so $m_1=m_2$ by \lcite{\NoEquiv}.  Therefore
there is exactly one least common multiple for $f$ and $g$, which we
denote as
  $$
  \lcm(f,g).
  $$

We next relate the notion of least common multiples to the categorical
notion of \stress{pull-backs}.

\state Proposition
  \label PullBack
  Let $f$ and $g$ be intersecting elements in $\M$, and write
$\lcm(f,g)=fp=gq$.  Then $(p,q)$ is the  unique pair of elements in
$\Mu$ such that
  \izitem 
  \zitem $fp=gq$, and 
  \zitem for every other pair of elements $p',q'\in\Mu$ such that
$fp'=gq'$, there exists a unique $r\in\Mu$ such that $p'=pr$, and
$q'=qr$.

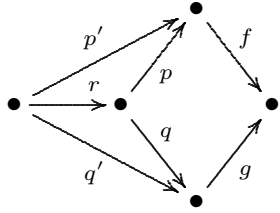
\beginmypicture
\setcoordinatesystem units <0.0010truecm, 0.0013truecm> point at 0 0
\pouquinho = 300
\put {$\bullet$} at 0000  0000
\put {$\bullet$} at 1000  1000
\put {$\bullet$} at 1000 -1000
\put {$\bullet$} at 2000  0000
\put {$\bullet$} at -1410  0000

\morph {0000}{0000}{1000}{1000} \arwlabel{p}{0600}{0300}
\morph {1000}{1000}{2000}{0000} \arwlabel{f}{1650}{0700}

\morph {0000}{0000}{1000}{-1000} \arwlabel{q}{0600}{-0300}
\morph {1000}{-1000}{2000}{0000} \arwlabel{g}{1650}{-0700}
\pouquinho = 200
\morph {-1400}{0000}{0000}{0000} \arwlabel{r}{-350}{200}

\pouquinho = 350
\morph {-1400}{0000}{1000}{1000} \arwlabel{p'}{-350}{700}
\morph {-1400}{0000}{1000}{-1000} \arwlabel{q'}{-350}{-700}

\endmypicture \bigskip

\proof
  We initially notice that the occurrence of black dots in our diagram
is not intended to give the idea of \stress{source} or \stress{range},
as we are not assuming that our semigroupoid is a category.

  Given $(p',q')$ as in (ii) notice that $m':=fp'$ is a common
multiple of $f$ and $g$.  Therefore $m\dil m'$, so there exists
$r\in\Mu$ such that $m'=mr$.  It follows that
  $$
  fp' = m' = mr = fpr.
  $$
  Since $f$ is monic we deduce that $p'=pr$, and a similar reasoning
gives $q'=qr$.  The uniqueness of $r$ follows from the fact that $p$
is monic.  Next let us address the uniqueness of $(p,q)$, by assuming
that $(p_1,q_1)$ and $(p_2,q_2)$ are two pairs satisfying (i) and
(ii).
  
Applying (ii) twice we conclude that there are $r$ and $s$ in $\Mu$
such that $p_2=p_1r$, and $q_2=q_1r$ on the one hand, and $p_1=p_2s$,
and $q_1=q_2s$, on the other.  Since
  $$
  p_1 = p_2s = p_1rs,
  $$
  we deduce that $rs=1$, whence $r=s=1$ and uniqueness follows.
  \proofend

We will now begin the actual construction of the inverse semigroup
$\SM$.  The first step is to consider a certain collection of subsets
of $\M$:

\definition 
  \label DefineMiddleSets
  We shall let $\Q$ denote the collection of all subsets of $\M$ of
the form 
  $$
  Q^F = \inters_{f\in F}\D f,
  $$
  where $F$ is a nonempty finite subset of $\M$.  By default the empty
set will also be included in $\Q$.

  Since we have prohibited 1 to be in any $\D f$, no member of $\Q$ is
allowed to contain 1.  In addition we have prohibited 1 to be in the
set $F$ above (recall that $F\c\M$), so $\D 1$ is never involved in
the intersection of sets making up $Q^F$, above.  Therefore $\M$ is
only a member of $\Q$ if there exists some $f\in\M$ for which $\D
f=\M$, which is not always the case, and rarely true in the examples
we wish to consider.

  It is noteworthy that $\Q$ is closed under intersections and hence
it is a semilattice with smallest element $\emptyset$.  As already
noticed it may or may not contain a largest element.

The underlying set of the inverse semigroup we wish to construct may
already be introduced:

\definition
  We will let $\SM$ denote the set 
  $$
  \SM = \{(f,A,g)\in \Mu\times\Q\times\Mu\ :\ \ A\c\D f\cap\D g\}.
  $$
  We will tacitly assume that all elements of $\SM$ of the form
$(f,A,g)$, with $A=\emptyset$, are identified with each other, forming
an \stress{equivalence class} which we will call \stress{zero} and
denote by 0.

Apart from the identification referred to above, no other
identifications will be implicitly or explicitly made.


We will now work towards defining the multiplication operation on
$\SM$.  The following rudimentary notation  will be extremely useful:

\definition Given $f\in\Mu$ and $A\in\Q$ we shall let
  $$
  f\inv(A) = \{g\in\D f: fg\in A\}.
  $$

The true meaning of $f\inv(A)$ is revealed next:

\state Proposition
  \label InverseImage
  Given $A\in\Q$ one has
  \izitem
  \zitem $1\inv(A)=A$,
  \zitem if $f\in A$, then $f\inv(A) =\D f$, 
  \zitem if $f\in\M\setminus A$, then $f\inv(A) =\emptyset$.

\proof Skipping the obvious first statement write
  $
  A = {\bigcap}_{h\in H}\D h,
  $
  where $H\c\M$ is a finite subset.  We begin by proving (iii) by
contradiction. So, supposing that $f\in\M$ and $f\inv(A)$ is nonempty
pick $g\in f\inv(A)$.  Then $fg\in A$, which means that $(h,fg)\in\Mt$
for all $h\in H$.  By the associativity property (and the fact that
$f\neq 1$) we deduce that $(h,f)\in\Mt$, and hence that $f\in A$,
proving (iii).

As for (ii) if $f\in A$, then again by the associativity property we
have that $(h,fg)\in\Mt$ for all $h\in H$ and $g\in\D f$, so $fg\in
A$, of $g\in f\inv(A)$.
  \proofend

A couple of elementary facts related to the above notation are:

  \state Proposition
  \label RulesInverImag
  Let $f,g\in\Mu$, and $A\in\Q$.
  \izitem 
  \zitem If $g\inv\big(f\inv(A)\big)$ is nonempty, then $g\in\D
f\cup\{1\}$,
  \zitem If $g\in\D f\cup\{1\}$, then $(fg)\inv(A) =
g\inv\big(f\inv(A)\big)$.

\proof
  (i) The result is obvious if either $f=1$ or $g=1$, so we suppose
$f,g\in\M$.  Clearly $f\inv(A)$ is nonempty, so $f\in A$ by
\lcite{\InverseImage.iii}, in which case $f\inv (A)=\D f$ by
\lcite{\InverseImage.ii}.  The hypothesis is then that $g\inv(\D f)$
is nonempty and, again by \lcite{\InverseImage.iii}, we conclude that
$g\in\D f$.

\medskip\noindent (ii) Left to the reader. 
  \proofend

We are now ready to describe the multiplication operation on $\SM$.

\definition
  \label DefineProduct
  Given $(f,A,g)$ and $(h,B,k)$ in $\SM$ we will let
  $$
  (f,A,g)(h,B,k) = \left\{\matrix{
  \Big(fu\ ,\ u\inv(A)\cap v\inv(B)\ ,\ kv\Big)\ , & \hbox{ if }
\lcm(g,h) = gu = hv,\cr
  \pilar{20pt}
  \hfill0\hfill\ , & \hbox{ if } g\disj h.\hfill 
  }\right.
  $$

There is a slight hitch in the above definition in the sense that
nothing guarantees that $fu$ and $kv$ are legal products.  However, if
$u$ is not in $\D f\cup \{1\}$, then $u$ is not in $A$ either, because
$A\c\D f$.  By \lcite{\InverseImage} we deduce that
$u\inv(A)=\emptyset$, and hence we define the product to be zero by
the rule that any triple with the empty set in the middle represents
zero in $\SM$, regardless of the fact that $fu$ is not defined.  The
same argument applies if $v$ is not in $\D g\cup
\{1\}$.

Rather than include a third clause in the definition above we shall
accept illegal products $fu$ and $kv$, only as long as the empty set
rule applies.

There is a diagrammatic interpretation for the product: in case
  $(f,A,g)(h,B,k)$ is nonzero,
  we have that $g\its h$, so we may write a \stress{pull-back diagram}
for $(g,h)$, as displayed in the diamond at the center of the diagram
below.

  \medskip
  \beginmypicture
  \setcoordinatesystem units <0.0010truecm, 0.0013truecm> point at 0 0
  \pouquinho = 300

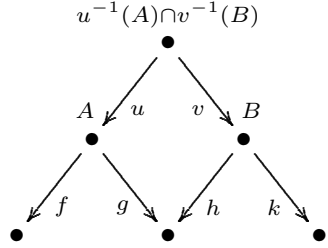
  \put {$\bullet$} at 0000 0000
  \put {$\bullet$} at 1000 1000
  \put {$\bullet$} at 2000 0000
  \put {$\bullet$} at 3000 1000
  \put {$\bullet$} at 4000 0000

  \put {$\bullet$} at 2000 2000

  \morph {1000}{1000}{0000}{0000} \arwlabel{f}{0600}{0300}
  \morph {1000}{1000}{2000}{0000} \arwlabel{g}{1400}{0300}
  \morph {3000}{1000}{2000}{0000} \arwlabel{h}{2600}{0300}
  \morph {3000}{1000}{4000}{0000} \arwlabel{k}{3400}{0300}

  \morph {2000}{2000}{1000}{1000} \arwlabel{u}{1600}{1300}
  \morph {2000}{2000}{3000}{1000} \arwlabel{v}{2400}{1300}

  \put{$\scriptstyle A$} at 900 1300
  \put{$\scriptstyle B$} at 3100 1300

  \put{$\scriptstyle u\inv(A)\cap v\inv(B)$} at 2000 2300
  \endmypicture
  \bigskip

Imagining that the element $(f,A,g)$ is represented by the triangle in
the lower left corner, including the decoration ``$A$" at its top
vertex, and similarly for $(h,B,k)$, the product is then represented
by the big triangle encompassing the whole diagram, with $u\inv(A)\cap
v\inv(B)$ as decoration.  This idea is used to prove the following:

\state Theorem
  \label TheInverseSGP
  Let $\S$ be an inverse semigroup satisfying \lcite{\StandingHyp}.
Then the multiplication on $\SM$ introduced above is well defined and
associative, and hence $\SM$ is a semigroup.  It is moreover an
inverse semigroup with zero, where the adjoint operation is given by
  $$
  (f,A,g)^* = (g,A,f).
  $$

\proof
  Since the middle coordinate of  $0$ is the empty set, it is clear
that
  $$
  0s =   s0 = 0
  \for s\in\SM.
  $$
  Next, 
  given $(f,A,g)$ and $(h,B,k)$ in $\SM$ we must show that their
product in fact lies in $\SM$.  Clearly this is so if the product
comes out to be zero, so we suppose otherwise, and hence
$\lcm(g,h)=gu=hv$, for suitable $u,v\in\M$.  We claim that
$u\inv(A)\c\D{fu}$.  This is obvious if $u=1$, and if $u\neq1$,
\lcite{\InverseImage} applies to give
  $$
  u\inv(A)\c\D u=\D {fu}.
  $$
  Similarly 
  $
  v\inv(B)\c\D {kv}.
  $
  Therefore 
  $u\inv(A)\cap v\inv(B)\c\D {fu}\cap \D {kv}$, proving that the
product does belong to $\SM$.  To prove associativity let $(l,C,m)$ be
a third element in $\SM$ and we shall prove that 
  $$
  \Big((f,A,g) (h,B,k)\Big) (l,C,m) =
  (f,A,g) \Big((h,B,k) (l,C,m)\Big).
  \subeqmark AssocEqtn
  $$
  We leave it up to the reader to show that if either
  $g\disj h$, or $k\disj l$, then both sides reduce to zero.  So we
assume instead that $g\its h$, and $k\its l$, and write
$\lcm(g,h)=gu=hv$, and $\lcm(k,l)=kx=ly$.

  We now claim that if $v\disj x$, then both sides of
\lcite{\AssocEqtn} vanish.  Assuming by contradiction that e.g.~the
left-hand side is nonzero then
  $$
  0\neq (f,A,g) (h,B,k)= (fu,u\inv(A)\cap v\inv(B),kv),
  $$
  and moreover $kv\its l$.  Write $\lcm(kv,l)=kvz=lw$.  By
\lcite{\PullBack} we deduce that $vz=xr$, and $w=yr$, for some
$r\in\Mu$.  This contradicts the assumption that $v\disj x$, and a
similar argument proves that the right-hand side vanishes as well, so
\lcite{\AssocEqtn} is proved under the hypothesis that $v\disj x$.  

We are then left to treat the case in which $v\its x$.  Write
  $$
  \lcm(v,x)=vp=xq,
  \subeqmark vpxq
  $$
  so we have built the diagram

  \medskip
  \beginmypicture
  \setcoordinatesystem units <0.0010truecm, 0.0013truecm> point at 0 0
  \pouquinho = 300

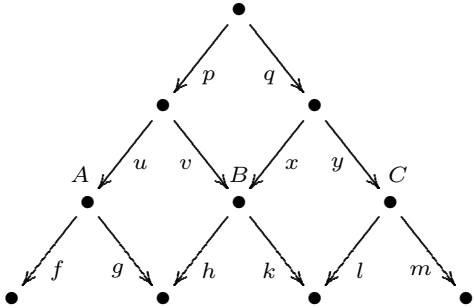
  \put {$\bullet$} at 0000 0000
  \put {$\bullet$} at 2000 0000
  \put {$\bullet$} at 4000 0000
  \put {$\bullet$} at 6000 0000
  \put {$\bullet$} at 1000 1000
  \put {$\bullet$} at 3000 1000
  \put {$\bullet$} at 5000 1000
  \put {$\bullet$} at 2000 2000
  \put {$\bullet$} at 4000 2000
  \put {$\bullet$} at 3000 3000

  \morph {1000}{1000}{0000}{0000} \arwlabel{f}{0600}{0300}
  \morph {1000}{1000}{2000}{0000} \arwlabel{g}{1400}{0300}
  \morph {3000}{1000}{2000}{0000} \arwlabel{h}{2600}{0300}
  \morph {3000}{1000}{4000}{0000} \arwlabel{k}{3400}{0300}
  \morph {5000}{1000}{4000}{0000} \arwlabel{l}{4600}{0300}
  \morph {5000}{1000}{6000}{0000} \arwlabel{m}{5400}{0300}

  \morph {2000}{2000}{1000}{1000} \arwlabel{u}{1700}{1400}
  \morph {2000}{2000}{3000}{1000} \arwlabel{v}{2300}{1400}
  \morph {4000}{2000}{3000}{1000} \arwlabel{x}{3700}{1400}
  \morph {4000}{2000}{5000}{1000} \arwlabel{y}{4300}{1400}

  \morph {3000}{3000}{2000}{2000} \arwlabel{p}{2600}{2300}
  \morph {3000}{3000}{4000}{2000} \arwlabel{q}{3400}{2300}

  \put{$\scriptstyle A$} at 0900 1300
  \put{$\scriptstyle B$} at 3000 1300
  \put{$\scriptstyle C$} at 5100 1300
  \endmypicture
  \bigskip

\noindent Even under all of the hypotheses assumed so far, it is still
possible that either side of \lcite{\AssocEqtn} vanish, given the role
played by the sets $A,B$, and $C$.  Obviously, if both sides vanish
there is nothing to prove, so we shall assume without loss of
generality that the left-hand side is nonzero.  This entitles us to
assume the following
  \iaitem
  \aitem $u\in A\cup\{1\}$,
  \aitem $v\in B\cup\{1\}$.
  \medskip\noindent
  Since $A\c\D f$ and $B\c\D k$ we have that $fu$ and $kv$ are indeed
legal products, and   in addition
  \aitem $kv\its l$.
  \medskip\noindent
  As above let $\lcm(kv,l)=kvz=lw$ and, using \lcite{\PullBack}, write
$vz=xr$, and $w=yr$, for some $r\in\Mu$.  Using \lcite{\vpxq} we may
further write $z=ps$, and $r=qs$, with $s\in\Mu$.

We would very much like to be able to perform the multiplication
``$kvp$", but even though $kv$ and $vp$ are known to be legal
multiplications we cannot use the only clause of the associativity
property which might fail if $v=1$.  Briefly assuming that $v=1$,
notice that
  $$
  kvz = kz = k(ps),
  $$
  which implies that $kp$ is a legal multiplication, thus taking care
of our concern.  We next observe that
  $$
  (kv) p =
  k(vp) = k(xq) = (kx)q = (ly)q = l(yq),
  $$
  so both  $kv$ and $l$ divide $kvp$, and hence 
  $$
  kvz = \lcm(kv,l) \dil kvp.
  $$
  Since all elements are monic this implies that $z\dil p$, but we
have seen above that $p\dil z$, and hence $p=z$ by \lcite{\NoEquiv}.
This gives $s=1$, and hence $r=q$, and finally $w=yq$.
  Summarizing, 
  $$
  \lcm(kv,l)= kvp =l yq.
  $$
  Recall that we are assuming the non-vanishing of the left-hand side
of \lcite{\AssocEqtn}, which is given by
  $$
  \big(fu,u\inv(A)\cap v\inv(B),kv\big)\ \big(l,C,m\big)=
  \Big(fup,p\inv\big(u\inv(A)\cap v\inv(B)\big)\cap
(yq)\inv(C),myq\Big) \$=
  \Big(fup,
    p\inv\big(u\inv(A)\big)\cap 
    p\inv\big(v\inv(B)\big)\cap
    q\inv\big(y\inv(C)\big),
    myq\Big).
  \subeqmark FinalFormLHS
  $$
  Using \lcite{\RulesInverImag} we have that $p\in\D u\cup\{1\}$, so
$up$ is a legal multiplication.  Moreover, by \lcite{\InverseImage} we
have that $up\in A\cup\{1\}$.  Since $A\c\D g$ we are allowed to set
  $$
  t = gup.
  $$
  Speaking of the right-hand side of \lcite{\AssocEqtn}, we have
  $$
  (f,A,g) \Big((h,B,k) (l,C,m)\Big) =
  (f,A,g) \Big(hx,x\inv(B)\cap y\inv(C),my\Big),
  \subeqmark RightHandDevel
  $$
  and we claim that $t$, defined just above, is the least common
multiple of $g$ and $hx$.  It is clear that $g\dil t$ and 
  $$
  t = gup = (hv)p = h(vp) = h(xq) = (hx)q,
  \subeqmark gupIShxq
  $$
  so $hx\dil t$, as well.  Let $s$ be a common multiple of $g$ and
$hx$, and write $s=ga=hxb$, for suitable $a,b\in\Mu$.
  Using \lcite{\PullBack} there is $c\in\Mu$ such that
  $a=uc$, and $xb=vc$.
  Observing that 
  $$
  h(xb) = ga = g(uc) = (gu)c = (hv)c = h(vc),
  $$
  and that $h$ is monic, we have $vc=xb$, so we may write $c=pd$, and
$b=qd$, for some $d\in\Mu$.  Thus
  $$
  s = hxb= (hx)(qd) = ((hx)q)d \ \={(\gupIShxq)} \ td,
  $$
  proving that $t\dil s$.  This shows that $t = \lcm(g,hx)$, and by
\lcite{\gupIShxq} we have that \lcite{\RightHandDevel} equals
  $$
  \Big(fup,
  (up)\inv(A)\cap q\inv\big(x\inv(B)\cap y\inv(C)\big),myq\Big) \$= 
  \Big(fup,
    p\inv\big(u\inv(A)\big)\cap 
    q\inv\big(x\inv(B)\big)\cap
    q\inv\big(y\inv(C)\big),
    myq\Big) =
  $$
  which coincides with \lcite{\FinalFormLHS} because $pv=qx$.  This
concludes the proof of associativity, and it remains to prove that
$\SM$ is an inverse semigroup with the indicated adjoint operation.
The reader will find no difficulty in proving that
  $$
  (f,A,g) (g,A,f) (f,A,g)= (f,A,g),
  $$
  so what is really at stake is the uniqueness of the adjoint.  So,
suppose that we are given $s$ and $t$ in $\SM$ such that
  $$
  s t s= s
  \and
  t s t = t.
  $$
  If either $s$ of $t$ vanishes it is immediate that
$t=s^*$, so we will suppose that $s, t\neq0$.  This
also implies that all products involved are nonzero.  Write $s =
(f,A,g)$ and $t=(h,B,k)$ and, observing that $g\its h$ and $k\its
f$, write
  $$
  \lcm(g,h) = gu=hv
  \and
  \lcm(k,f)= kx = fy.
  $$
  We then have
  $$
  s =   s t s = 
  (f,A,g)(h,B,k)(f,A,g) = 
  \big(fu,u\inv(A)\cap v\inv(B),kv\big) \ (f,A,g) = \ldots
  $$
  Further writing $\lcm(kv,f)=kvz=fw$, the above equals
  $$
  \ldots = \Big(fuz,z\inv\big(u\inv(A)\cap v\inv(B)\big)\cap
w\inv(A),gw\Big).
  \subeqmark aaba
  $$
  Since this coincides with $s$ we have that $fuz=f$, and $gw=g$,
and hence $u=z=w=1$, because all elements are monic.  Therefore 
  $h\dil g$ and 
  $k\dil f$.  Applying the same reasoning to the equation $t
s t = t$
we deduce that 
  $g\dil h$ and 
  $f\dil k$, so $f=k$, and $g=h$, by \lcite{\NoEquiv}.  This also
implies that $v=1$, and turning to the middle coordinate of
\lcite{\aaba} we conclude that $A\cap B=A$, so $A\c B$.  By symmetry
we also have that $B\c A$, so in fact $A=B$, and this finally gives
$t=s^*$.
  \proofend

It is not hard to see that the idempotent semilattice $E(\SM)$ of
$\SM$ is formed by the elements $(f,A,g)\in\SM$, for which $f=g$.
Given the importance of the order relation in $E(\SM)$ we shall now
describe it in explicit terms: 

\state Proposition
  \label OrderInISGforGPD
  Let $(f,A,f)$ and $(g,B,g)$ be idempotents in $E(\SM)$, with $A\neq
\emptyset$.  Then
  \izitem 
  \zitem $(f,A,f)\leq(g,B,g)$, if and only if $g\dil f$ and, writing
$f=gh$, for $h\in\Mu$, one has that $A\c h\inv(B)$,
  \zitem if $f=g$, then $(f,A,f)\leq(f,B,f)$, if and only if $A\c B$,
  \zitem if $g\in\M$, and $B=\D g$, then $(f,A,f)\leq (g,\D g,g)$, if
and only if $g\dil f$.
  \zitem if $g=1$, and $f\in\M$, then $(f,A,f)\leq(1,B,1)$, if and
only if $f\in B$.

\proof Beginning with (i), supposing that $f=gh$, and that $A\c
h\inv(B)$, we have that $\lcm(f,g) = f = f1 = gh$, so
  $$
  (f,A,f)(g,B,g) = 
  (f,A\cap h\inv(B),gh) =
  (f,A,f),
  $$
  so $(f,A,f)\leq(g,B,g)$, as desired.  Conversely, assuming that
$(f,A,f)\leq(g,B,g)$ we have that 
  $$
  (f,A,f)(g,B,g) = (f,A,f) \neq 0,
  $$
  because of the assumption that $A\neq\emptyset$.  This implies that
$f\its g$, so we write $\lcm(f,g) = fk=gh$, with $k,h\in\Mu$, and then
  $$
  0 \neq (f,A,f) = (f,A,f)(g,B,g) = (fk,k\inv(A)\cap h\inv(B),gh).
  $$
  Notice that since the elements we are comparing above are nonzero,
there is no identification involved, meaning that equality only holds
when the correspondent components agree.
  This implies in particular that $f=fk=gh$, and hence $k=1$, by
\lcite{\DefineMonic}.  This also proves that $g\dil f$.  Another
conclusion to be drawn from the equation displayed above is that
  $$
  A = k\inv(A)\cap h\inv(B) = A\cap h\inv(B),
  $$
  which implies that $A \c h\inv(B)$.
  The other points follow easily from (i).
  \proofend

Referring to \lcite{\OrderInISGforGPD}, observe that if $A=\emptyset$,
then $(f,A,f)=0$, and hence $(f,A,f) \leq(g,B,g)$, regardless of any
other relationship between $f$, $g$, and $B$.

\state Proposition
  \label ConditionForIts
  Let $(f,A,f)$ and $(g,B,g)$ be idempotent elements in $E(\SM)$.
Then $(f,A,f)\its(g,B,g)$ if and only if there are $u,v\in\Mu$ such
that $fu=gv$, and $u\inv(A)\cap v\inv(B)$ is nonempty.

\proof
  Supposing that $(f,A,f)\its(g,B,g)$ we have that 
  $$
  0\neq (f,A,f)(g,B,g) = \big(fu,u\inv(A)\cap v\inv(B),gv\big),
  $$
  where $\lcm(f,g) = fu = gv$.  Obviously $u\inv(A)\cap v\inv(B)$ is
nonempty.

Conversely, suppose that $u$ and $v$ exist as in the statement.  Since
$u\inv(A)$ is nonempty we have by \lcite{\InverseImage} that $u\in
A\cup\{1\}\c\D f\cup\{1\}$, so $fu$ is meaningful, and so is $gv$.
Moreover it is clear that 
  $
  u\inv(A) \c \D{fu},
  $
  and
  $
  v\inv(B) \c \D{gv},
  $
  so we have that 
  $$
  u\inv(A) \cap v\inv(B) \c\D{fu}\cap \D{gv},
  $$
  proving that 
  $
  \big(fu,u\inv(A)\cap v\inv(B),gv\big),
  $
  is an element of $E(\SM)$, which is clearly nonzero.  Using
\lcite{\OrderInISGforGPD.i} we see that this element is smaller than
both   $(f,A,f)$ and $(g,B,g)$, and hence 
  $$
  0 \neq   \big(fu,u\inv(A)\cap v\inv(B),gv\big) \leq
  (f,A,f)  (g,B,g),
  $$
  proving that 
  $(f,A,f)\its (g,B,g)$.
  \proofend

The idempotents $(1,B,1)$ have an interesting property which is
described in our next result:

\state Proposition
  \label LessOrOrhogoNew
  Let $(f,A,f)$ be an idempotent such that $f\neq1$, and let $B\in\Q$.
  \izitem
  \zitem If $f\in B$, then $(f,A,f)\leq (1, B,1)$.
  \zitem If $f\notin B$, then $(f,A,f)\perp (1, B,1)$.

\proof
  We have
  $$
  (1, B,1) (f,A,f) =
  (f,f\inv( B)\cap A,f) = \ldots
  $$
  Assuming that $f\in B$ we have by \lcite{\InverseImage}, that
$f\inv( B)=\D f$.  In addition it is implicit that $A\c\D f$, so the
above equals
  $$
  \ldots = (f,\D f\cap A,f) = (f,A,f),
  $$
  proving (i).
  On the other hand, if $f$ is not in $ B$, we have that
$f\inv(B)=\emptyset$, and hence
  $(1, B,1) (f,A,f) =0$.
  \proofend


In view of the relevance of $E^*$-unitary inverse semigroups in the
characterization of the Hausdorff property for the groupoid of germs
given in \lcite{\HausdorffGPG} and \lcite{\ISGLattice}, it is
interesting to find sufficient conditions for $\SM$ to be
$E^*$-unitary.
  By analogy with \lcite{\DefineMonic} we will say that an element
$f\in\M$ is \stress{epic} if for every $g,h\in\Mu$ we have
  $$
  gf=hf \imply g=h.
  $$

\state  Proposition
  Let $\M$ be a semigroupoid satisfying \lcite{\StandingHyp}, and such
that all of its elements are epic.  Then $\SM$ is $E^*$-unitary
 
\proof Suppose that an element $(f,A,g)$ in $\SM$ dominates a nonzero
idempotent $(h,B,h)$.  Then
  $$
  0 \neq (h,B,h) = (f,A,g)(h,B,h),
  $$
  which implies that $g\its h$, so we may write $\lcm(g,h) = gu=hv$,
for suitable elements $u,v\in\Mu$, and 
  $$
  (h,B,h) = (fu,u\inv(A)\cap v\inv(B),hv).
  $$
  Since this is nonzero we conclude that $h=fu=hv$.  It follows that
  $$
  fu=hv=gu,
  $$
  and since $u$ is epic, we conclude that $f=g$, thus proving that
$(f,A,g)$ is an idempotent.
  \proofend


\section{Representations of semigroupoids}
  \label RepSGPDSect
  As before we fix a semigroupoid $\M$ satisfying
\lcite{\StandingHyp}.  
  We shall begin this section by introducing several important notions
inspired in \cite{\ExSemiGpdAlg}, most of which are homonyms of similar
notions introduced earlier in this work in the context of inverse
semigroups, such as \stress{representations}, \stress{covers}, and
\stress{tight representations}.  Once the appropriate context is clear
we believe the double meanings will cause no confusion.

  %
  %

\definition
  \label DefineRepSGPD
  A representation of the semigroupoid $\M$ in an inverse semigroup
with zero $\S$ is a map
  $
  \orep : \M \to \S,
  $
  such that for every $f,g\in\M$, one has that:
  \izitem 
  \zitem
  $
  \orep_f\orep_g  = \left\{\matrix{
    \orep_{fg}, &\hbox{ if } (f,g)\in\Mt, \cr
    \vrule height 15pt width 0pt
    0, & \hbox{ otherwise.}\hfill}\right.
  $
  \medskip\noindent
  Moreover the initial and final projections 
  $$
  \ini^\orep_f = \orep_f^*\orep_f \and \fin^\orep_g=\orep_g\orep_g^*, 
  \subeqmark DefineIniAndFin
  $$
  respectively, are required to satisfy
  \medskip
  \zitem $\fin^\orep_f\fin^\orep_g=0$, if $f\disj g$, 
  \zitem $\ini^\orep_f\fin^\orep_g = \fin^\orep_g$, if $(f,g)\in\Mt$.
  \medskip\noindent In case $\S$ is an inverse semigroup formed by
partial isometries on a Hilbert space $H$, and containing the zero
operator, we will say that $\orep$ is a \stress{representation of
$\M$ on $H$}.

\medskip We insist that, since the symbol ``$\disj$" is used in (ii)
in the context of semigroupoids, its meaning is to be taken from
\lcite{\DefineDisjSGPD}, and not from \lcite{\ItsAndPerpInSL}.

One might wonder what happens to the element appearing in the
left-hand side of the equation in \lcite{\DefineRepSGPD.iii}, in case
$(f,g)$ is not in $\Mt$.  The answer is provided by
\lcite{\DefineRepSGPD.i}, since in this case
  $$
  \ini^\orep_f\fin^\orep_g = \orep_f^*(\orep_f \orep_g) \orep_g^* =0.
  \eqmark DefineRepSGPDiv
  $$

  Should the context leave no room for confusion we will abbreviate
the notations $\ini^\orep_f$ and $\fin^\orep_g$, to $\ini_f$ and
$\fin_g$, respectively.

\definition
  \label DefineCoverSGPD
  Let $\Gamma$ be any subset of the semigroupoid $\M$.  A subset $H\c
\Gamma$ will be called a \stress{cover}
  for $\Gamma$ if for every $f\in \Gamma$ there exists $h\in H$ such
that $h\its f$.  If moreover the elements of $H$ are mutually disjoint
then $H$ will be called a \stress{partition} of $\Gamma$.


\definition
  \label DefTightRepSGPD
  Let $\S$ be an inverse semigroup of partial isometries on a Hilbert
space $H$, containing the identically zero operator.
  A representation $\orep$ of $\M$ on $\S$ is said to be
\stress{tight} if for every finite subsets $F,G\c\M$, and for
every finite covering $H$ of
  $$
  \M^{F,G}:= \bigcap_{f\in F} \D f\cap
  \bigcap_{g\in G} \M\setminus\D g,
  $$
  one has that
  $
  \ds \bigvee_{h\in H}\fin_h =
  \ini_{F,G},
  $
  where
  $$
  \ini_{F,G} :=
  \prod_{f\in F}\ini_f \prod_{g\in G}(1-\ini_g).
  $$
 
Observe that the Definition given in \scite{\ExSemiGpdAlg}{4.5} requires that
the above holds for every finite subsets $F$ and $G$ of $\Mu$, as
opposed to $\M$.  However notice that since $\D 1=\M$, and the
convention adopted there says that $\ini_1=1$, one has that
$\M^{F,G}=\emptyset$, whenever $1\in G$, and the above condition holds
vacuously.  If, on the other hand, $1$ is in $F$, then
$\M^{F,G}=\M^{F',G}$, where $F'=F\setminus\{1\}$, at the same time
that $\prod_{f\in F}\ini_f=\prod_{f\in F'}\ini_f$.  Therefore we see
that the above definition is equivalent to \scite{\ExSemiGpdAlg}{4.5},
regardless of our use of $\M$ in place of $\Mu$.

In the above definition the recipient inverse semigroup $\S$ needs to
be embedded in $B(H)$ or otherwise neither the supremum $\bigvee_{h\in
H}\fin_h$, nor the term $1-\ini_g$, would make sense.  This situation
may however be generalized by assuming that $E(\S)$ admits the
structure of a Boolean algebra which is compatible with the order of
$E(\S)$, in which case one might say that $\S$ is a \stress{Boolean
inverse semigroup}.  This and related results may be found in
\cite{\ExelAlgebra}.

  Tight representations have the following good behavior with respect
to least common multiples:

\state Proposition
  \label TightLCM
  Suppose that for every $f\in\M$, and every $h\in\D f$, there exists
a finite partition $H$ of $\D f$, such that $h\in H$.  If $f,g\in\M$
are such that $f\its g$, and $\orep$ is a tight representation of $\M$
one has that
  $$
  \fin_f\fin_g = \fin_m,
  $$
  where $m=\lcm(f,g)$.

\proof
  In case $f\dil g$, write $g=fh$, with $h\in\Mu$, and notice that
  $$
  p_f p_g  =
  \orep_f \orep_f^* \orep_f \orep_h \orep_h^*\orep_f^* =
  \orep_f \orep_h \orep_h^*\orep_f^* =
  \orep_g  \orep_g^* = p_g =  p_{m},
  $$
  since $m=g$.  The case $g\dil f$ may be treated similarly so
we next assume that there are $u$ and $v$ in $\M$ (as opposed to $\Mu$), such
that $fu=gv=m$.  By hypothesis let $H$ and $K$ be finite
partitions of $\D f$ and $\D g$, respectively, such that $u\in H$ and
$v\in K$.  Given that our representation is tight we have
  $$
  q_f = \bigvee_{h\in H} p_h = \sum_{h\in H} p_h,
  $$
  where the last equality follows from the fact that the $p_h$ are
pairwise orthogonal projections by \lcite{\DefineRepSGPD.ii}.
Therefore
  $$
  p_f =
  \orep_f \orep_f^* \orep_f \orep_f^* = 
  \orep_f q_f \orep_f^* =
  \sum_{h\in H} \orep_fp_{h} \orep_f^*=
  \sum_{h\in H} p_{fh},
  $$
  and similarly
  $
  p_g =
  \ds\sum_{k\in K} p_{gk}.
  $
  So
  $$
  p_f p_g =
  \sum_{h\in H} \sum_{k\in K} p_{fh}  p_{gk}.
  $$
  Among the pairs $(h,k)\in H\times K$ one clearly has the pair
$(u,v)$ for which
  $
  p_{fu} p_{gv} = p_{m},
  $
  and the proof will be complete once we show that
  $p_{fh} p_{gk} = 0$, for all other pairs $(h,k)$.  Thus assume that
$(h,k)\in H\times K$ is such that either $h\not=u$ or $k\not=v$.  We
will in fact prove that $fh\disj gk$, and hence the conclusion will
follow from \lcite{\DefineRepSGPD.ii}.  Arguing by contradiction
suppose that $fhx = gky$, where $x,y\in\Mu$.   By \lcite{\PullBack} we
have that $hx =ur$, and $ky=vr$, for some $r\in\Mu$, but this says
that  $h\its u$ and $k\its v$, a contradiction.
  \proofend

This result motivates the following:

\definition 
  \label DefinePreserveLCM
  A representation $\orep$ of $\M$ in an inverse semigroup $\S$ is
said to \stress{respect least common multiples} if for every
intersecting pair of elements $f$ and $g$ in $\M$, one has that
  $$
  p_fp_g = p_{\lcm(f,g)}.
  $$

The following is an important property of these representations.  It
is related to equation \scite{\KP}{3.1} in the context of finitely
aligned higher rank graphs.

\state Proposition
  \label TrocaEstrela
  Suppose that $\orep$ is a map from $\M$ into a *-semigroup%
  \fn{A *-semigroup is a semigroup equipped with an involution
$s\mapsto s^*$, which satisfies $(st)^*=t^*s^*$.  For reasons which
will soon become apparent we do not suppose that $\S$ is an inverse
semigroup here.}
  $\S$, satisfying \lcite{\DefineRepSGPD.i--iii} and the equation
displayed in \lcite{\DefinePreserveLCM}. Suppose moreover that 
  $$
  \orep_f\orep_f^*\orep_f=\orep_f
  \for f\in\M.
  $$
  (This is clearly the case if $\S$ is an inverse semigroup and
$\orep$ is a representation of $\M$ in $\S$ respecting least common
multiples).
  Given a pair of intersecting
elements $f,g\in\M$, with $f\neq g$, write $\lcm(f,g)=m=fh=gk$, with
$h,k\in\Mu$.  Then
  $$
  \orep_f^*\orep_g = \orep_h\orep_k^*.
  $$

\proof
  It is conceivable that $h$ or $k$ be equal to 1, in which case the
right hand side of the equation above needs clarification.  First
observe that since we are assuming that $f\neq g$, the situation in
which both $h$ and $k$ are equal to 1 will never arise.  If $h=1\neq
k$, then $\orep_h\orep_k^*$ is supposed to mean $\orep_k^*$, and vice
versa.  The best way to deal with this problem is to think that
$\orep_1=1$, where the last occurrence of $1$  is a \stress{multiplier}
of $\S$, meaning an element which may not belong to $\S$, but which is
allowed to multiply elements of $\S$ in such a way that
  $$
  1s=s1=s
  \for s\in\S.
  $$
  Also we set $1^*=1$.  

Observe that whenever $(f,g)\in\Mt$, we have by
  \lcite{\DefineRepSGPD.iii} that
  $$
  \orep_f^*\orep_f\orep_h =
  \orep_f^*\orep_f\orep_h\orep_h^*\orep_h =
  q_f p_h \orep_h =
  p_h \orep_h =
  \orep_h.
  \subeqmark QTimesRep
  $$
 
Let us now prove the statement under the special assumption that
$f\dil g$.  In this case we have $m=g$, and $k=1\neq h$.  Therefore
  $$
  \orep_f^*\orep_g = 
  \orep_f^*\orep_f\orep_h \={(\QTimesRep)}
  \orep_h = 
  \orep_h\orep_k^*.
  $$
  Assuming instead that $g\dil f$ one may give a similar proof, or
just use adjoints, so we next suppose that $f$ and $g$ do not divide
each other.  This implies that $h,k\neq1$, so $h\in\D f$ and $k\in \D
g$.
  We then have
  $$
  \orep_f^*\orep_g =
  \orep_f^*\orep_f\orep_f^*\orep_g\orep_g^*\orep_g =
  \orep_f^* p_f p_g\orep_g = 
  \orep_f^* p_m\orep_g =
  \orep_f^* \orep_m\orep_m^*\orep_g \$=
  \orep_f^* \orep_f\orep_h\orep_k^*\orep_g^*\orep_g =
  \orep_f^* \orep_f\orep_h(\orep_g^*\orep_g\orep_k)^* \={(\QTimesRep)}
  \orep_h\orep_k^*.
  \proofend
  $$

We will often deal with representations of $\M$ in inverse semigroups
of partial isometries on a Hilbert space and in the $\lcm$-preserving
case it is possible to omit any reference to that semigroup:

\state Proposition
  Let $H$ be a Hilbert space and let $\orep:\M\to B(H)$ be a map 
satisfying \lcite{\DefineRepSGPD.i--iii}.  Suppose moreover that, for
every $f,g\in\M$, one has that
  \izitem
  \zitem $\orep(f)$ is a partial isometry,
  \zitem $\ini_f$ and  $\ini_g$ commute,
  \zitem if $f\its g$, then $\fin_f\fin_g = \fin_{\lcm(f,g)}$.
  \medskip \noindent Then the smallest multiplicative subsemigroup of
$B(H)$ which is closed under adjoints and contains the range of
$\orep$ is an inverse semigroup and moreover $\orep$ is a
representation of $\M$ in it.

  \proof
  For each finite nonempty subset $F\c\M$, let 
  $$ 
  \ini_F = \prod_{f\in F}\ini_f.
  $$
  The order in which the above elements are  multiplied is irrelevant in
view of (ii).
  Extending $\orep$ to $\Mu$ by setting $\orep_1=1$, let 
  $$
  \calbox{S} = \{\orep_f \ini_F \orep_g^*: f,g\in\Mu,\ F\subseteq \M
\hbox{ is finite and nonempty}\} \cup \{0\}.
  $$
  Notice that for every $(f,g)\in\Mu\times\Mu\setminus\{(1,1)\}$, we
have that 
  $$
  \orep_f\orep_g^* = 
  \orep_f\orep_f^*\orep_f \orep_g^*\orep_g\orep_g^* =
  \orep_f\ini_{\{f,g\}}\orep_g^* \in 
  \calbox{S},
  $$
  so in particular $\calbox{S}$ contains the range of $\orep$.
  We next then claim that $\calbox{S}$ is a multiplicative
subsemigroup of $B(H)$.
  To prove it let us be given $f,g,h,k\in\Mu$, and finite nonempty
subsets $F,G\subseteq\M$. We will prove that
  $$
  \orep_f \ini_F \orep_g^* \  \orep_h \ini_G \orep_k^*
  \subeqmark RandomProduct
  $$
  either vanishes or equals
  $
  \orep_u \ini_H \orep_v^*,
  $
  for suitable $u,v\in\Mu$, and $H\subseteq\M$.  We divide the proof
in several cases, according to the values of $g$ and $h$:

\medskip\noindent {\tensc Case 1}: {$g=h=1$}.
   Take $u=f$, $H=F\cup G$, and $v=k$.

\medskip\noindent {\tensc Case 2}: {$g=1,\ h\not=1$}.
  Notice that 
  $
  \ini_F\orep_h = 
  \ini_F \fin_h \orep_h,
  $
  while 
  $$
  \ini_F\fin_h = 
  \left\{\matrix{
  \fin_h &, \hbox{ if } (f,h)\in\Mt,\ \forall f\in F, \cr\cr
  0 &, \hbox{ otherwise},\hfill
  }\right.
  $$
  by \lcite{\DefineRepSGPD.iii} and \lcite{\DefineRepSGPDiv}.
  Thus
  $\ini_F\orep_h$ either vanishes or agrees with $\orep_h$.  Therefore
\lcite{\RandomProduct} either vanishes or equals
  $$
  \orep_f \ini_F  \orep_h \ini_G \orep_k^* =
  \orep_f \orep_h \ini_G \orep_k^* =
  \orep_{fh} \ini_G \orep_k^*,
  $$
  where we have also assumed that $(f,h)\in\Mt$, or else
\lcite{\RandomProduct} again vanishes.

\medskip\noindent {\tensc Case 3}: {$g\not=1,\ h=1$}.  Follows from
case 2, and taking adjoints.

\medskip\noindent {\tensc Case 4}: {$g,h\in\M$, and $g\disj h$}.  Then
$\orep_g^*\orep_h = \orep_g^*\fin_g\fin_h\orep_h = 0$, by
\lcite{\DefineRepSGPD.ii}, and hence \lcite{\RandomProduct} vanishes.

\medskip\noindent {\tensc Case 5}: {$g=h\neq1$}.
   Take $u=f$, $H=F\cup\{g\}\cup G$, and $v=k$.

\medskip\noindent {\tensc Case 6}: {$g,h\in\M$, $g\neq h$, and $g\its
h$}.  Applying \lcite{\TrocaEstrela} to the multiplicative *-semigroup
of all bounded operators on $H$, we have that $\orep_g^*\orep_h =
\orep_u\orep_v^*$, with $u,v\in\Mu$.  Then \lcite{\RandomProduct}
equals
  $$
  \orep_f \ini_F \orep_u\orep_v^* \ini_G \orep_k^*,
  $$ 
  and the result follows as in case 2.  It is now clear that
$\calbox{S}$ is the *-subsemigroup of $B(H)$ generated by the range
of $\orep$.  We will now prove that $\calbox{S}$ consists of partial
isometries.  For this let $u\in\calbox{S}$ be a generic element and
write
  $
  u=\orep_f\ini_F \orep_g^*.
  $
  Observing that
  $$
  u=
  \orep_f\ini_f\ini_F \ini_g\orep_g^* =
  \orep_f\ini_{\{f\}\cup  F\cup \{g\}}\orep_g^*,
  $$
  we may assume that $f,g\in F$.  
  We then have
  $$
  uu^*u = 
  \orep_f\ini_F\orep_g^*\ \orep_g\ini_F \orep_f^*\ \orep_f\ini_F\orep_g^* =
  \orep_f\ini_F \ini_g \ini_F \ini_f\ini_F \orep_g^* =
  \orep_f\ini_F \orep_g^* = u,
  $$
  so $u$ is a partial isometry as claimed.  It is well known that any
subsemigroup of $B(H)$ consisting of partial isometries, and which is
closed under adjoints, is an inverse semigroup.  It is obvious that 
$\calbox{S}$ is closed under adjoints, so it is an inverse semigroup.
Obviously $\orep$ is then a representation of $\M$ in $\calbox{S}$.
  \proofend


Our next long term goal is to show a close relationship between tight
representations of $\M$ and tight representations of $\SM$ (which in
turn are related to representations of the groupoid of germs, by
\lcite{\EquivTightRepISGandGPD}).  An important ingredient in this
relationship is a representation of $\M$ in $\SM$ to be introduced next.
 
\state Proposition
  \label PropsOfTau
  The map $\tau:\M \to \SM$ defined by 
  $
  \tau_f = (f,\D f,1),
  $ 
  for all $f\in\M$,
  is a representation of $\M$ in $\SM$, which respects least common
multiples and moreover satisfies
  $$
  \ini^\tau_f = \tau_f^*\tau_f = (1,\D f,1)
  \and
  \fin^\tau_f=\tau_f\tau_f^* = (f,\D f,f),
  $$
  for all $f\in\M$.

\proof
  For $f\in\M$ one has that
  $$
  \tau_f^*\tau_f =
  (1,\D f,f) (f,\D f,1) =
  \big(1,1\inv(\D f) \cap 1\inv(\D f),1\big) =
  (1,\D f,1),
  $$
  and similarly one proves that 
  $
  \tau_f\tau_f^* =
  (f,\D f,f)
  $.
  If we are also given $g\in\M$, then
  $$
  \tau_f\tau_g =
  (f,\D f,1) (g,\D g,1) =
  \big(fg,g\inv(\D f)\cap\D g,1\big).
  $$
  If $(f,g)\in\Mt$ then $g\in\D f$ and hence $g\inv(\D f) = \D g$, by
\lcite{\InverseImage}, so
  $$
  \tau_f\tau_g =
  \big(fg,\D g,1\big) =
  \big(fg,\D {fg},1\big) =
  \tau_{fg}.
  $$
  If $(f,g)\notin\Mt$ then $g\notin\D f$ and using
\lcite{\InverseImage} again we have that $g\inv(\D f) = \emptyset$, so
  $$
  \tau_f\tau_g =
  (fg,\emptyset,1) = 0,
  $$
  regardless of the fact that $fg$ is meaningless.  With respect to
\lcite{\DefineRepSGPD.ii} assume that $f\disj g$.  Then
  $$
  \tau_f^*\tau_g = (1,\D f,f) (g,\D g, 1) = 0,
  $$
  by definition, from which one sees that $\fin_f\fin_g=0$.  If
$g\in\D f$ then
  $$
  \tau_f^* \tau_f \tau_g = 
  (1,\D f, 1) (g, \D g, 1) =
  (g, g\inv(\D f)\cap \D g, 1) =
  (g, \D g, 1) =
  \tau_g,
  $$
  and hence \lcite{\DefineRepSGPD.iii} follows.
  To conclude we must show that $\tau$ respects least common
multiples, so let $f,g\in\M$ be intersecting elements.  Write
$\lcm(f,g) = m = fu=gv$, for $u,v\in\Mu$, and notice that
  $$
  u\inv(\D f)=
  u\inv(f\inv(\M)) \={(\RulesInverImag.ii)}
  (fu)\inv(\M) =
  \D {fu} = \D m,
  $$
  and similarly
  $ v\inv(\D g)=\D m$.  So
  $$
  \fin^\tau_f \fin^\tau_g =
  (f,\D f,f) (g,\D g,g) = 
  (fu,u\inv(\D f)\cap v\inv(\D g),gv) =
  (m,\D m, m) = 
  \fin^\tau_m.
  \proofend
  $$

It is interesting to notice that if $f$ is a spring, that is, if $\D
f=\emptyset$, then $\tau_f=0$.  This is partly the reason why springs
are cumbersome elements to deal with.

\bigskip Given $A,B\in\Q$ it is immediate that
  $$
  (1,A,1)(1,B,1) = (1,A \cap B,1).
  $$
  So, given any $A\in\Q$, say
  $
  A = \inters_{h\in H}\D h,
  $
  where $H$ is a nonempty finite subset of $\M$, we have that
  $$
  (1,A,1) =
  \big(1,\inters_{h\in H}\D h,1\big) =
  \prod_{h\in H} (1,\D h,1) =
  \prod_{h\in H} \tau_h^*\tau_h.
  $$
  In addition, if $f,g\in\M$ are such that $\D f\cap\D g\supseteq A$,
we have
  $$
  (f,\D f, 1) (1,A,1)(1,\D g, g) =
  (f,\D f\cap A\cap \D g, g) =
  (f,A, g),
  $$
  so we have proved that:

\state Proposition
  \label LambdaGenerates
  Let $(f,A,g)\in\SM$, and write $A = {\inters}_{h\in H}\D h$, for
some nonempty finite subset $H\c\M$.  Then 
  $$
  (f,A,g) = 
  \tau_f \Big(\prod_{h\in H} \tau_h^*\tau_h\Big)\tau_g^*.
  $$

We therefore see that the range of $\tau$ generates $\SM$ as an
inverse semigroup.  Our next result uses $\tau$ to express the first
relationship between representations of $\M$ and representations of
$\SM$.

\state Proposition
  \label FromRepISGtoGPOD
  If $\irep$ is a representation of $\SM$ on a Hilbert space $H$
  \lcite{in the sense of \DefineRepSoDoCara} such that $\irep(0)=0$,
then the composition $\orep=\irep\circ\tau$ is a representation of
$\M$
  \lcite{in the sense of \DefineRepSGPD}, which respects least
common multiples. If moreover 
  $\M$ has no springs and $\irep$ is tight
  \lcite{in the sense of \DefTightRepISG}, then $\orep$ is tight 
  \lcite{in the sense of \DefTightRepSGPD}.

\proof
  That $\orep$ is a representation preserving least common multiples
follows immediately from \lcite{\PropsOfTau}, and the fact that
$\irep(0)=0$.
  Next, assuming that $\irep$ is tight, let us prove that the same
applies to $\orep$. So let $F,G\c\M$ be finite sets and let $H$ be a
cover for $\M^{F,G}$ in the sense of \lcite{\DefineCoverSGPD}.  We
must prove that
  $$
  \bigvee_{h\in H}\orep_h\orep_h^* =
  \prod_{f\in F}\orep_f^*\orep_f \prod_{g\in G}(1-\orep_g^*\orep_g).
  \subeqmark ThisIsWhatWeWant
  $$
  Letting
  $$
  \eqmatrix{
  X &=& \{(1,\D f,1): f\in F\}, \cr
  \pilar{12pt}
  Y &=& \{(1,\D g,1): g\in G\}, \cr
  \pilar{12pt}
  Z &=& \{(h,\D h,h): h\in H\},
  }
  $$
  we claim that $Z$ is a cover of $E(\SM)^{X,Y}$, in the sense of
\lcite{\DefineCoverInSLat}.
  In order to prove our claim let $(k,C,k)$ be a nonzero idempotent in
$E(\SM)^{X,Y}$.  Therefore $C$ is nonempty, so we pick some ${c}\in C$.
Given that $C\c\D k$, we may speak of $k{c}$, and it is easy to see,
based on the fact that $\D{k{c}}=\D {c}$, and
\lcite{\OrderInISGforGPD.i}, that
  $$
  (k{c},\D{k{c}},k{c})\leq (k,C,k).
  $$
  Since $(k,C,k)$ is in $E(\SM)^{X,Y}$, the same applies to 
$(k{c},\D{k{c}},k{c})$, and hence for every $f\in F$, and
$g\in G$, we have 
  $$
  (k{c},\D{k{c}},k{c})\leq (1,\D f,1)
  \and
  (k{c},\D{k{c}},k{c})\perp (1,\D g,1).
  $$
  Noticing that $(k{c},\D{k{c}},k{c})$ is nonzero because $\M$ has no
springs, and hence it cannot be simultaneously orthogonal and smaller
than any other element, we have by \lcite{\LessOrOrhogoNew} that
  $k{c}\in\D f$ and $k{c}\notin\D g$.  This says that
  $$
  k{c}\in
  \bigcap_{f\in F} \D f\cap
  \bigcap_{g\in G} \M\setminus\D g =
  \M^{F,G},
  $$
  so there exists some $h\in H$ such that $k{c}\its h$, and hence we may
write $k{c}x=hy$, for some $x,y\in\Mu$.  Using
\lcite{\OrderInISGforGPD.iii} one has that
  $(k{c}x,\D{k{c}x},k{c}x)$ is simultaneously smaller than $(h,\D h,h)$ and 
  $(k{c},\D{k{c}},k{c})$.  This implies that
  $$
  0\neq
  (k{c}x,\D{k{c}x},k{c}x) \leq
  (h,\D h,h)(k{c},\D{k{c}},k{c}) \leq
  (h,\D h,h)(k,C,k),
  $$
  proving that $(h,\D h,h)\its(k,C,k)$.  This shows that $Z$ is indeed
a cover for $E(\SM)^{X,Y}$, and because $\irep$ is assumed to be a
tight representation of $\SM$ we have
  $$
  \bigvee_{z\in Z} \irep_z =
  \prod_{x\in X} \irep_x 
  \prod_{y\in Y} (1- \irep_y).
  \subeqmark ThisIsWhatWeGot
  $$
  For every $z=(h,\D h,h)\in Z$, with $h\in H$, notice that
  $$
  \irep_z =
  \irep(h,\D h,h)=
  \irep(\tau_h\tau_h^*) =
  \orep_h\orep_h^*,
  $$
  and similarly 
  $$
  \irep_x = \orep_f^*\orep_f
  \and
  \irep_y = \orep_g^*\orep_g,
  $$
  for every $x=(1,\D f,1)\in X$, and $y=(1,\D g,1)\in Y$, so we see
that \lcite{\ThisIsWhatWeWant} follows from \lcite{\ThisIsWhatWeGot}.
This proves that $\orep$ is tight.
  \proofend


\section{The Boolean algebra of Domains}
  \label DomainSec
  In our pursuit of a bijective correspondence between representations
of the semigroupoid $\M$ and of its associated inverse semigroup $\SM$
we would like to show that any representation $\orep$ of $\M$ may be
\stress{extended} to $\SM$, meaning that there exists a representation
$\irep$ of $\SM$, such that $\orep = \irep\circ \tau$, thus obtaining
a converse to \lcite{\FromRepISGtoGPOD}.  In order to understand the
difficulties in doing so let us temporarily suppose that $\irep$ has
been found.  If $F$ is a nonempty finite subset of $\M$, then
  $$
  \prod_{f\in F} \orep_f^*\orep_f =
  \prod_{f\in F} \irep(\tau_f^*\tau_f) = 
  \irep\Big(\prod_{f\in F} \tau_f^*\tau_f\Big) =
  \irep\Big(\prod_{f\in F} (1,\D f, 1)\Big) \$=
  \irep\Big(1,\inters_{f\in F} \D f, 1\Big) =
  \irep\Big(1,Q^F, 1\Big),
  $$
  where our use of $Q^F$ is the same as in \lcite{\DefineMiddleSets}.
Implicit in the above calculation is the fact that 
  $\prod_{f\in F} \orep_f^*\orep_f$
  depends only on $Q^F$, and not on $F$. While for a general
representation $\orep$ this may fail, we will prove that this does
hold provided $\orep$ is tight.  Under that hypothesis we will not
only prove that $\irep$ exists, but also that it is tight.  Our
correspondence will therefore involve tight representations only.

A large part of the effort in accomplishing our goal will be spent on
studying the behavior of tight representations of $\M$ with respect to
a certain Boolean algebra of 
subsets of $\M$. 

 \cryout{Throughout this section we therefore fix a semigroupoid $\M$
satisfying \lcite{\StandingHyp} and a tight representation $\orep$ of
$\M$ on a Hilbert space $H$.}

\definition
  A subset $X\c\M$ will be called a \stress{domain}, provided it
belongs to the Boolean subalgebra $\Dom$ of $\P(\M)$ generated by
$\{\D f: f\in\M\}$.

If $F$ and $G$ are finite subsets of $\M$, then the set
  $$ 
  \M^{F,G} = \bigcap_{f\in F} \D f\cap
  \bigcap_{g\in G} \M\setminus\D g,
  $$
  already employed in \lcite{\DefTightRepSGPD}, is clearly a domain.
Moreover it is easy to see that any member of $\Dom$ may be written as
the union of a finite collection of sets each of which has the above
form.

If one is to decide whether or not some $h$ in $\M$ belongs to a given
domain $D\in\Dom$, that task will consist of a perhaps logically
complicated check depending on whether or not $h\in\D f$, for several
elements $f$ in $\M$.  It is therefore easy to see that for every
$k\in \D h$ one has
  $$
  h\in D \iff hk\in D.
  \eqmark PropOfDomain
  $$

In this section we will not be dealing with any representation
of $\M$ other  than $\orep$, so we will drop the superscripts 
in the $\ini^\orep_f$ and $\fin^\orep_g$ of 
\lcite{\DefineIniAndFin}.

We wish to define a map $\iniq:\Dom\to B(H)$ such that
  $\iniq(\D f)= \ini_f$, and which is a Boolean
algebra homomorphism in the sense that
  \bitem $\iniq(\emptyset) = 0$,
  \bitem $\iniq(C\cap D) = \iniq(C)\iniq(D)$, and 
  \bitem $\iniq(\tilde C) = 1-\iniq(C)$,
  \medskip\noindent
  for every $C,D\in\Dom$, where $\tilde C$ denotes the complement of
$C$ in $\M$.  Clearly we will have as a consequence that
  $$
  \iniq(C\cup D) = 1-\iniq(\tilde C\cap \tilde D) =
  1- \big(1-\iniq(C)\big)\big(1-\iniq(D)\big) =
  \iniq(C) + \iniq(D) - \iniq(C)\iniq(D),
  $$
  which is precisely the join, or supremum $\iniq(C)\vee\iniq(D)$, of
the commuting projections $\iniq(C)$ and $\iniq(D)$ in $B(H)$.
  If we are to succeed in obtaining $\iniq$ then for every finite subsets
$F,G\c\M$ we must have
  $$
  \iniq(\M^{F,G}) = 
  \prod_{f\in F} \iniq(\D f) \prod_{g\in G} \big(1-\iniq(\D g)\big) =
  \prod_{f\in F} \ini_f \prod_{g\in G} (1-\ini_g) =
  \ini_{F,G},
  $$
  where $\ini_{F,G}$ was already employed in
\lcite{\DefTightRepSGPD}.

In the next result we will take a first step in the direction of the
goal stated above by showing that $\ini_{F,G}$ does indeed depend only
on the set $\M^{F,G}$, and not on $F$ and $G$.

\state Proposition
  \label FirstIndependenceResultForQ
  Let $F$, $G$, $H$, and $K$ be finite subsets of $\M$.  
  \izitem 
  \zitem If 
  $\M^{F,G}\c\M^{H,K}$, then
  $\ini_{F,G}\leq\ini_{H,K}$.
  \zitem If 
  $\M^{F,G}=\M^{H,K}$, then 
  $\ini_{F,G}=\ini_{H,K}$.

  \proof
  Assume that $\M^{F,G}\c\M^{H,K}$. For each fixed $h\in H$ notice
that $\M^{F,G}\c\M^{H,K}\c\D h$, and hence
  $$
  \emptyset =
  \M^{F,G}\cap(\M\setminus\D h) = 
  \M^{F,G\cup\{h\}}.
  $$
  Since $\orep$ is tight we deduce that 
  $$
  0=
  \ini_{F,G\cup\{h\}} = 
  \ini_{F,G} (1-\ini_h),
  $$
  so that $\ini_{F,G} \leq \ini_h$.  On the other hand, for
every $k\in K$, we have that $\M^{F,G}\c\M\setminus \D k$, and hence
  $$ 
  \emptyset =
  \M^{F,G}\cap\D k =
  \M^{F\cup\{k\},G}.
  $$
  Since $\orep$ is tight we deduce that 
  \ $
  0=
  \ini_{F\cup\{k\},G} =
  \ini_{F,G} \,\ini_k,
  $ \
  so that \ $\ini_{F,G} \leq 1-\ini_k$. \  Therefore
  $$
  \ini_{F,G} \leq
  \prod_{h\in H}\ini_h \prod_{k\in K}(1-\ini_k) =
  \ini_{H,K},
  $$
  proving (i), and hence also (ii)
  \proofend

If a domain $D$ has the form $\M^{F,G}$, we may then define $\iniq(D) =
\ini_{F,G}$, without worrying about other possible descriptions of
$D$ in the form $\M^{H,K}$.  In the next result we shall consider the
possibility that some domains may be described in several ways as unions
of sets of the form $\M^{F,G}$.

\state Proposition
  \label OneAsUnionOfOthers
  Let $\{F_i\}_{i=0}^n$ and $\{G_i\}_{i=0}^n$ be two collections of
  finite subsets of $\M$, such that
  $$
  \M^{F_0,G_0} = 
  \M^{F_1,G_1} \cup
  \M^{F_2,G_2} \cup
  \ldots \cup
  \M^{F_n,G_n},
  $$
  Then
  \ $
  \ini_{F_0,G_0} = \bigvee_{i=1}^n \ini_{F_i,G_i} .
  $

\proof
  Let $H = F_0\cup G_0\cup F_1\cup G_1\cup \ldots\cup F_n\cup G_n$.
For each subset $X\c H$ we let
  $$
  E^X = \M^{X,H\setminus X}.
  $$
  It is then easy to see that the $E_X$ are pairwise disjoint and that
${\union}_{X\in\Psmall(H)}E^X = \M$.  Likewise, letting
  $$
  e_X = \ini_{X,H\setminus X},
  $$
  it is easy to see that the $e_X$ are pairwise orthogonal projections
such that $\sum_{X\in\Psmall(H)}e_X = 1$.  In order to prove the
statement it is therefore enough to show that for every $X\in\P(H)$
one has that
  $$
  \ini_{F_0,G_0}e_X = \bigvee_{i=1}^n \,\ini_{F_i,G_i}e_X.
  \subeqmark MainEqInFirstIndependenceResultForQ 
  $$ 
  Since $e_X=0$ whenever $E^X=\emptyset$, by
\lcite{\FirstIndependenceResultForQ}, we need only consider those $X$
for which $E^X$ is nonempty.  Let us thus fix $X\c H$, with
$E^X\neq\emptyset$.
  For each $i=0,\ldots,n$, observe that if 
  $F_i\c X$, and $G_i\c H\setminus X$, then $E^X\c\M^{F_i,G_i}$ and
$e_X\leq\ini_{F_i,G_i}$.  On the other hand if either $F_i\not\c X$,
or $G_i\not\c H\setminus X$, then necessarily
$E^X\cap\M^{F_i,G_i}=\emptyset$, and $e_X\perp \ini_{F_i,G_i}$.

\medskip\noindent{\tensc Case 1:} Assume that there exists some $i\geq
1$ such that $F_i\c X$, and $G_i\c H\setminus X$.  Then the right-hand
side of \lcite{\MainEqInFirstIndependenceResultForQ} equals $e_X$.
Moreover $E^X\c\M^{F_i,G_i}\c\M^{F_0,G_0}$, from where we deduce that
$F_0\c X$, and $G_0\c H\setminus X$, and hence
$e_X\leq\ini_{F_0,G_0}$, so  the left-hand side of
\lcite{\MainEqInFirstIndependenceResultForQ} also equals $e_X$.

\medskip\noindent{\tensc Case 2:} Assume that there is no $i\geq 1$
such that $F_i\c X$, and $G_i\c H\setminus X$.  Then $e_X\perp
\ini_{F_i,G_i}$, for all $i\geq1$, and hence the right-hand side of
\lcite{\MainEqInFirstIndependenceResultForQ} vanishes.  Moreover $E^X$
is disjoint from each $\M^{F_i,G_i}$, with $i\geq1$, and hence it is
also disjoint from $\M^{F_0,G_0}$.  Thus, it cannot be that $F_0\c X$,
and $G_0\c H\setminus X$, and hence $e_X\perp \ini_{F_0,G_0}$, proving
that the left-hand side of
\lcite{\MainEqInFirstIndependenceResultForQ} also vanishes.
  \proofend

The next result will finally allow us to define the map we are
seeking:

\state Proposition
  \label MapOnDomains
  For every $D\in\Dom$, write
  $
  D = {\union}_{j=1}^{n} \M^{F_j,G_j},
  $
  where the $F_j$ and $G_j$ are finite subsets of $\M$, and define
  $
  \iniq(D) = \bigvee_{j=1}^{n} \ini_{F_j,G_j}.
  $
  Then $\iniq:\Dom\to B(H)$ is a well defined map which moreover
satisfies
  \izitem 
  \zitem   $\iniq(\D f)= \ini_f$, 
  \zitem $\iniq(\emptyset) = 0$,
  \zitem $\iniq(C\cap D) = \iniq(C)\iniq(D)$,
  \zitem $\iniq(C\cup D) = \iniq(C)\vee\iniq(D)$,
  \zitem $\iniq(\tilde D) = 1-\iniq(D)$,
  \medskip\noindent
  for every $f\in\M$, and every $C,D\in\Dom$.

\proof
  To show well definedness   suppose that $D$ is a domain which may be written in two ways as
  $$
  D =
  \union_{j=1}^{n_1} \M^{F_j^1,G_j^1} =
  \union_{j=1}^{n_2} \M^{F_j^2,G_j^2},
  $$
  where the $F_j^i$ and $G_j^i$ are finite subsets of $\M$. Fix $k\leq
n_1$ and notice that
  $$
  \M^{F_k^1,G_k^1} =
  \M^{F_k^1,G_k^1} \cap D =
  \union_{j=1}^{n_2} \M^{F_k^1,G_k^1}\cap \M^{F_j^2,G_j^2} =
  \union_{j=1}^{n_2} \M^{F_k^1\cup F_j^2,G_k^1\cup G_j^2}.
  $$
  By \lcite{\OneAsUnionOfOthers} we conclude that 
  $$
  \ini_{F_k^1,G_k^1} =
  \bigvee_{j=1}^{n_2} \ini_{F_k^1\cup F_j^2,G_k^1\cup G_j^2} =
  \bigvee_{j=1}^{n_2} 
    \ini_{F_k^1,G_k^1} \
    \ini_{F_j^2,G_j^2} =
  \ini_{F_k^1,G_k^1} \Big(
  \bigvee_{j=1}^{n_2} 
    \ini_{F_j^2,G_j^2}\Big),
  $$
  showing that 
  $$
  \ini_{F_k^1,G_k^1} \leq \bigvee_{j=1}^{n_2} \ini_{F_j^2,G_j^2}.
  $$  
  Since $k$ is arbitrary we deduce that 
  $$
  \bigvee_{k=1}^{n_1} \ini_{F_k^1,G_k^1} \leq
  \bigvee_{j=1}^{n_2} \ini_{F_j^2,G_j^2},
  $$
  and by symmetry we obtain the reverse inequality, hence proving that
the two possibly different descriptions of $D$ lead to the same
proposed value of $\iniq(D)$.  

We leave it for the reader to prove (i--iii) and we will verify (v)
next.  Supposing initially that $D$ has the form $D=\M^{F,G}$, we have
  $$
  \tilde D =
  \union_{f\in F} \M \setminus \D f \cup \union_{g\in G} \D g,
  $$
  hence 
  $$
  \iniq(\tilde D) =
  \bigvee_{f\in F} (1-\ini_f) \vee \bigvee_{g\in G} \ini_g = 
  1- \prod_{f\in F} \ini_f \prod_{g\in G} (1- \ini_g) =
  1- \ini_{F,G} = 
  1 - \iniq(D).
  $$
  In the general case write $D = \union\nolimits_{j=1}^n D_j$, where
each $D_j$ is of the above form, then
  $$
  \iniq(\tilde D) = 
  \iniq\Big(\inters_{j=1}^n \tilde D_j\Big) =
  \prod_{j=1}^n \iniq(\tilde D_j) =
  \prod_{j=1}^n \big(1-\iniq(D_j)\big) =
  1-\bigvee_{j=1}^n \iniq(D_j) =
  1-\iniq(D).
  $$
  As already seen, (iv) follows from (iii) and (v).
  \proofend

For the record we notice the following:

\state Corollary
  \label iniqOfGenericQ
  If $H\c\M$ is a finite nonempty subset and $A = \inters_{h\in H}\D
h$ then $\iniq(A) = \prod\limits_{h\in H}\orep_h^*\orep_h$.

Recall that the condition for a representation of $\M$ to be tight is
that
  $$
  \bigvee_{h\in H}\fin_h =
  \prod_{f\in F}\ini_f \prod_{g\in G}(1-\ini_g),
  $$
  whenever $F,G\c\M$ are finite sets and $H$ is a finite cover for
$\M^{F,G}$.
  If we denote by $D$ the domain $D=\M^{F,G}$, then the above
condition may be expressed as
  $
  \bigvee_{h\in H}\fin_h =
  \iniq(D).
  $
  One may therefore ask if the same is true for every domain.  The
next result proves that this is in fact true.

\state Proposition
  \label ExtendedTight
  Let $D$ be a domain and let $H$ be a finite cover for $D$ in the
sense of \lcite{\DefineCoverSGPD}, then
  $$
  \bigvee_{h\in H}\fin_h = \iniq(D).
  $$

\proof
  Write 
  $
  D = \union\nolimits_{j=1}^{n} \M^{F_j,G_j},
  $
  where the $F_j$ and $G_j$ are finite subsets of $\M$.  By assumption
we have that $H\c D$, so if we put
  $$
  H_j = H\cap\M^{F_j,G_j},
  $$
  we will have that $H = \union\nolimits_{j=1}^{n} H_j$.  We claim
that $H_j$ is a cover for $\M^{F_j,G_j}$ for each $j\leq n$.  In fact,
given any $k\in\M^{F_j,G_j}$, we in particular have that $k\in D$.  By
hypothesis there exists some $h\in H$ such that $h\its k$, and we may
therefore choose $u,v\in\Mu$ such that $hu=kv$.  By
\lcite{\PropOfDomain} we have that $h\in\M^{F_j,G_j}$, so $h\in H_j$,
and the claim is proved.  Since we are assuming $\orep$ to be tight it
follows that
  $$
  \bigvee_{h\in H_j}\fin_h = \ini_{F_j,G_j},
  $$
  and therefore that
  $$
  \bigvee_{h\in H}\fin_h = 
  \bigvee_{j=1}^n \ \bigvee_{h\in H_j} \fin_h =
  \bigvee_{j=1}^n \ini_{F_j,G_j} =
  \iniq(D).
  \proofend
  $$


\section{Extending representations}
  \label ExtRepSec
  The sole aim of this section is to prove the following:

\state Theorem
  \label MainCorrTheorem
  Let $\M$ be a semigroupoid without springs in which every element is
monic, and such that every intersecting pair of elements admits a
least common multiple.  Given a tight representation $\orep$ of\/ $\M$
on a Hilbert space $H$, which respects least common multiples, there
exists a unique representation $\irep$ of the inverse semigroup $\SM$
  such that $\irep = \orep\circ \tau$.
Moreover $\irep$ is tight.

Observing that the representation of $\M$ given by
\lcite{\FromRepISGtoGPOD} necessarily respects least common multiples,
the above result cannot survive without assuming that $\orep$ also has
this property.

Given $(f,A,g)$ in $\SM$ write $A = \inters_{h\in H}\D h$, where
$H\c\M$ is a finite nonempty subset.  Recall from
\lcite{\LambdaGenerates} that
  $$
  (f,A,g) = \tau_f \Big(\prod_{h\in H} \tau_h^*\tau_h\Big) \tau_g^*,
  $$
  so if we want a representation $\irep$ of $\SM$ such that
$\irep\circ\tau=\orep$, we have no choice but to define
  $
  \irep(f,A,g) = \orep_f \Big(\prod_{h\in H} \orep_h^*\orep_h\Big)
\orep_g^*.
  $
  This immediately gives uniqueness and, in view of
\lcite{\iniqOfGenericQ}, it also suggests that we define
  $$
  \irep(f,A,g) = \orep_f \iniq(A) \orep_g^*,
  $$
  where $Q$ is given by \lcite{\MapOnDomains}.
  For $f\in\M$ we then have that
  $$
  \irep(\tau_f)=
  \irep(f,\D f,1) =
  \orep_f \iniq(\D f) =
  \orep_f \orep_f^* \orep_f =
  \orep_f,
  $$
  so $\orep = \irep\circ\tau$, as required.  
  It is also clear that $\irep$ preserves the star operation.

We next claim that if $A\in\Q$ and $f\in\M$ one has that
  $$
  \iniq(A)\orep_f = \left\{\matrix{
  \orep_f, & \hbox {if } f\in A,\hfill \cr
  0, & \hbox {otherwise.}}\right.
  \eqmark PiQRepf
  $$
  To prove it write $A=\inters\kern0pt_{h\in H}\D h$, for some finite
nonempty subset $H$ of $\M$.
  Assuming initially that $f\in A$, we then have that $f\in\D h$, for
all $h\in H$, and hence
  $
  \orep_h^* \orep_h \orep_f = \orep_f,
  $
  by \lcite{\DefineRepSGPD.iii}. So
  $$
  \iniq(A)\orep_f = \Big(\prod_{h\in H}\orep_h^*\orep_h\Big) \orep_f =
\orep_f.
  $$
  On the other hand, if $f\notin A$, then $f\notin\D h$ for some $h\in
H$, and hence 
  $\orep_h^* \orep_h \orep_f = 0$, by \lcite{\DefineRepSGPD.i}, so
our claim is proved.  Using \lcite{\InverseImage}, we may express
\lcite{\PiQRepf} alternatively as
  $$
  \iniq(A)\orep_f = \orep_f\iniq\big(f\inv(A)\big).
  \eqmark  PiQRepfBis
  $$
  The advantage of this over \lcite{\PiQRepf} is that it holds
inclusively for $f=1$, while the former does not.

We are now prepared to prove that $\irep$ is multiplicative.  Given
$(f,A,g)$ and $(h,B,k)$ in $\SM$ we then have to show that
  $$
  \irep(f,A,g)\irep(h,B,k) =
  \irep\big((f,A,g)(h,B,k)\big).
  \eqmark PiIsRep
  $$
  Observing that the left-hand side equals
  $$
  \orep_f \iniq(A)\orep_g^* \orep_h \iniq(B)\orep_k^*,
  $$
  we see that it vanishes whenever $g\disj h$, because
  $$
  \orep_g^* \orep_h =
  \orep_g^*\orep_g\orep_g^* \orep_h\orep_h^*\orep_h =
  \orep_g^*\fin_g \fin_h\orep_h =
  0,
  $$
  by \lcite{\DefineRepSGPD.ii}.  Still under the assumption that
$g\disj h$, we have that $(f,A,g)(h,B,k)=0$, by definition, and since
$\irep(0)=0$, the right-hand side of \lcite{\PiIsRep} also vanishes.
Thus \lcite{\PiIsRep} is true provided $g\disj h$, and we may then
suppose that $g\its h$, writing
  $\lcm(g,h) = m  = gu=hv$,
  with $u,v\in\Mu$.

Assuming, as we are, that $\orep$ respects least common multiples, we
wish to apply \lcite{\TrocaEstrela} to describe $\orep_g^* \orep_h$,
but for this we also need to assume we are in the special case in
which $g\neq h$.  Under this premise we have that $\orep_g^* \orep_h =
\orep_u\orep_v^*$, and hence the left-hand side of \lcite{\PiIsRep}
equals
  $$
  \orep_f \iniq(A)\orep_u\orep_v^* \iniq(B)\orep_k^* \={(\PiQRepfBis)}
  \orep_f
\orep_u\iniq\big(u\inv(A)\big)\iniq\big(v\inv(B)\big)\orep_v^*
\orep_k^* \$=
  \orep_f \orep_u\iniq\big(u\inv(A)\cap v\inv(B)\big)\orep_v^*
\orep_k^*.
  \eqmark LeftHandSideOfPiIsRep
  $$
  If $u\inv(A)\cap v\inv(B)=\emptyset$, then the above is zero, but so
is the right-hand side of \lcite{\PiIsRep}, and hence equality is
established.  In the event that $u\inv(A)\cap v\inv(B)$ is nonempty we
have by \lcite{\InverseImage} that 
  $$
  u\in A\cup\{1\}\c\D f\cup\{1\} 
  \and
  v\in B\cup\{1\}\c\D k\cup\{1\}.
  $$
  Thus $\orep_f\orep_u=\orep_{fu}$ and $\orep_k\orep_v=\orep_{kv}$, so
\lcite{\LeftHandSideOfPiIsRep} equals
  $$
  \orep_{fu}\iniq\big(u\inv(A)\cap v\inv(B)\big)\orep_{kv}^* = 
  \irep\big(fu,u\inv(A)\cap v\inv(B),kv\big) =
  \irep\big((f,A,g)(h,B,k)\big),
  $$
  proving \lcite{\PiIsRep} under the assumption that $g\neq h$, and the
only case left to be discussed is that in which $g=h$.  Under this
assumption observe that, since $A\c\D g$, we have by
\lcite{\MapOnDomains} that
  $$
  \iniq(A) \orep_g^*\orep_g =
  \iniq(A) \iniq(\D g) =
  \iniq(A \cap \D g) =
  \iniq(A),
  $$
  and hence the left-hand side of \lcite{\PiIsRep} equals
  $$
  \orep_f \iniq(A)\orep_g^* \orep_g \iniq(B)\orep_k^* =
  \orep_f \iniq(A)\iniq(B)\orep_k^* =
  \orep_f \iniq(A\cap B)\orep_k^* \$=
  \irep(f,A\cap B,k) =
  \irep\big((f,A,g)(g,B,k)\big),
  $$
  proving that $\irep$ is indeed a representation of $\SM$.
Summarizing our findings so far we have:

\state Lemma
  Under the hypotheses of \lcite{\MainCorrTheorem} the map $\irep:\SM
\to B(H)$, defined by
  $$
  \irep(f,A,g) = \orep_f \iniq(A) \orep_g^*
  \for (f,A,g)\in\SM,
  $$
  is a representation of $\SM$ satisfying $\irep\circ\tau=\orep$.

The remaining of this section will be dedicated to proving the last
sentence of \lcite{\MainCorrTheorem}, namely that $\irep$ is tight.
  The characterization of tightness given in \lcite{\AltLatTightRep}
will prove itself useful, but to employ it we must first check either
(i) or (ii) of \lcite{\ChaeckWithNonEmptX}.  We therefore suppose that
\lcite{\ChaeckWithNonEmptX.ii} fails, meaning that $E(\SM)$ admits a
finite cover, say $Z$.  Let us classify the elements $(f,A,f)$ of $Z$
according to whether $f=1$ or not by setting
  $$
  Z' = \{(f,A,f)\in Z: f\in\M\}
  \and
  Z'' = \{(f,A,f)\in Z: f=1\}.
  $$
  For each $(1,A,1)$ in $Z''$ write $A=\inters_{g\in G}\D g$, where
$G\c\M$ is finite and nonempty, and choose at random some $g_A\in G$.
Once this is done we have that $A\c\D{g_A}$ and hence
  $$
  (1,A,1) \leq (1,\D{g_A},1) =
  \tau_{g_A}^*\tau_{g_A}.
  $$
  With respect to the elements $(f,A,f)$ in $Z'$ notice that $A\c\D f$
and hence
  $$
  (f,A,f)\leq(f,\D f,f) =
  \tau_f\tau_f^*.
  $$
  Substituting each element of $Z$ appearing in the left-hand side of
the two inequalities displayed above by the respective right-hand side
we therefore obtain a set of the form
  $$
  W = \{\tau_g^*\tau_g: g\in G\}\cup\{\tau_f\tau_f^*: f\in F\},
  $$
  which is clearly also a cover for $E(\SM)$.
  We next claim that $F$ is a cover for 
  $$
  \M^{\emptyset,G} = \inters_{g\in G}\M\setminus \D g,
  $$
  in the sense of \lcite{\DefineCoverSGPD}.  To prove this let
$h\in\M^{\emptyset,G}$, and notice that $\tau_h\tau_h^*$ must
necessarily intersect some element of $W$.  If that element is of the
form $\tau_g^*\tau_g$, for some $g\in G$, then
  $$
  0\neq 
  \tau_h\tau_h^*\ \tau_g^*\tau_g =
  (h,\D h,h) (1,\D g,1) = 
  (h,\D h\cup h\inv(\D g),h),
  $$
  so $h\inv(\D g)$ is nonempty, and hence $h\in\D g$ by
\lcite{\InverseImage.iii}, but this contradicts the fact that
$h\in\M^{\emptyset,G}$.  The conclusion is that the element of $W$
which intersects $\tau_h\tau_h^*$ must be some $\tau_f\tau_f^*$, with
$f\in F$.  In this case
  $$
  0\neq
  \tau_h\tau_h^*\ \tau_f\tau_f^* =
  (h,\D h,h) (f,\D f,f)
  $$ which implies that $h\its f$, concluding the proof of our claim.
Since we are assuming that $\orep$ is tight we have that
  $$
  \bigvee_{f\in F}\orep_f\orep_f^* =
  \prod_{g\in G}(1- \orep_g^*\orep_g) = 
  1 - \bigvee_{g\in G}\orep_g^*\orep_g,
  $$
  and hence that 
  $$
  \Big(\bigvee_{f\in F}\orep_f\orep_f^*\Big) \vee
  \Big(\bigvee_{g\in G}\orep_g^*\orep_g\Big) =1.
  $$
  This implies that
  $$
  \bigvee_{w\in W} \irep(w) = 
  \Big(\bigvee_{f\in F}\irep(\tau_f\tau_f^*)\Big) \vee
  \Big(\bigvee_{g\in G}\irep(\tau_g^*\tau_g)\Big) =
  \Big(\bigvee_{f\in F}\orep_f\orep_f^*\Big) \vee
  \Big(\bigvee_{g\in G}\orep_g^*\orep_g\Big) =1.
  $$

We have therefore proven:

\state Lemma
  Either $E(\SM)$ does not admit any finite cover or there exists a
finite cover $W$ such that $\bigvee_{w\in W} \irep(w) =1$.

As already mentioned this result enables us to use
\lcite{\AltLatTightRep} to attempt a proof that $\irep$ is tight.
Therefore, given $(f,A,f)\in E(\SM)$ and a finite cover $Z$ for
$(f,A,f)$
  we need to prove that
  $$
  \bigvee_{z\in Z}\irep(z) \geq\irep(f,A,f).
  \eqmark BigGoalTight
  $$
  We will argue in two different ways according to whether $f=1$ or
not, so let us begin by assuming that $f=1$.
  As before write $Z = Z'\cup Z''$, where
  $$
  Z' = \{(h,C,h)\in Z: h\in\M\} = \big\{(h_i,C_i,h_i)\big\}_{i=1}^n,
  $$ and $$
  Z'' = \{(h,C,h)\in Z: h=1\} = \big\{(1,D_i,1)\big\}_{i=1}^m.
  $$
  Our proof will be by induction on $|Z'| =n$, so let us first treat
the case in which $n=0$.

Thus $Z''$ is a cover for $(1,A,1)$, and we claim that 
$A=\union\nolimits_{i=1}^m D_i$.  Since each $(1,D_i,1)\leq (1,A,1)$
we have that $D_i\c A$.  On the other hand, given $f\in A$, we have
that $(f,\D f,f)\leq (1,A,1)$.  Assuming that $\M$ has no springs 
  we see that $(f,\D f,f)$ is nonzero, so there is some $(1,D_i,1)\in
Z''$, such that $(f,\D f,f)\its(1,D_i,1)$, which means that $f\in D_i$, 
by \lcite{\LessOrOrhogoNew}.  This proves our claim, therefore
  $$
  \bigvee_{z\in Z}\irep(z) =
  \bigvee_{i=1}^m\irep(1,D_i,1) =
  \bigvee_{i=1}^m \iniq(D_i) \={(\MapOnDomains.iv)}
  \iniq\Big(\union_{i=1}^m D_i\Big) =
  \iniq(A) = 
  \irep(1,A,1),
  $$
  thus proving \lcite{\BigGoalTight} for $f=1$, and $n=0$.

Still assuming that $f=1$, but now that $n\geq1$, pick any $j\leq n$,
and let $(h,C,h)=(h_j,C_j,h_j)$.  Since $(h,C,h)\leq (1,A,1)$ we have
by \lcite{\OrderInISGforGPD.iv} that $h\in A$.  Incidentally notice
that \lcite{\OrderInISGforGPD.iv} requires that $C$ be nonempty, which
we may assume, since otherwise $(h,C,h)$ may be deleted from the
covering $Z$ without altering the left-hand side of
\lcite{\BigGoalTight}.
Given that $h\in A$, we have that
  $$
  (h,\D h,h)\leq (1,A,1).
  $$
  We next claim that 
  $$
  Z_h := \tau_h^*\, Z\, \tau_h
  $$
  is cover for $\tau_h^*\tau_h = (1,\D h,1)$.  To prove the claim let
$0\neq \gamma \leq\tau_h^*\tau_h$ and observe that
  $$
  \tau_h\gamma\tau_h^* \leq
  \tau_h\tau_h^*\tau_h\tau_h^* =
  \tau_h\tau_h^* = (h,\D h,h) \leq
  (1,A,1).
  $$
  Since $Z$ is a cover for $(1,A,1)$, and $\tau_h\gamma\tau_h^*$ is
nonzero (or else 
  $
  \gamma = \tau_h^*\tau_h\gamma\tau_h^*\tau_h =0
  $),
  there exists $z\in Z$ such that
$\tau_h\gamma\tau_h^* \its z$, so
  $$
  0\neq
  \tau_h\gamma\tau_h^*z =
  (\tau_h\tau_h^*)\tau_h\gamma\tau_h^*z =
  \tau_h\gamma\tau_h^*z(\tau_h\tau_h^*),
  $$
  which implies that
  $
  \gamma\tau_h^*z\tau_h \neq 0,
  $
  and hence that $\gamma\its \tau_h^*z\tau_h$, proving the claim.  

Let us now decompose $Z_h$ as the union $Z_h'\cup Z_h''$, where
  $$
  Z_h' = \{(g,B,g)\in Z_h: g\in\M\} \and
  Z_h'' = \{(g,B,g)\in Z_h: g=1\},
  $$
  in the same way we did with $Z$, because we are interested in the
number of elements of $Z_h'$, given that our proof is by induction on
this parameter.  Notice that for every $i\leq m$
  $$
  \tau_h^* (1,D_i,1)\tau_h =
  (1,\D h,h)(1,D_i,1)(h,\D h,1) =
  \big(1,\D h\cap h\inv(D_i),h\big)(h,\D h,1) \$=
  \big(1,\D h\cap h\inv(D_i)\cap\D h,1\big) \in Z_h'',
  $$
  which means that
  $$
  \tau_h^*\, Z'' \tau_h\c Z_h''.
  $$
  If the reader is expecting a similar inclusion with single primes
replacing double primes, he or she will be surprised to find that
there is an element of $Z'$ which migrates to $Z_h''$ when
\stress{conjugated} by $\tau_h$, namely
  $$
  \tau_h^* (h,C,h)\tau_h =
  (1,\D h,h)(h,C,h)(h,\D h,1) =
  (1,\D h\cap C\cap \D h,1) =
  (1,C,1) \in Z_h''.
  $$ 
  It follows that $Z_h'$ has at most $n-1$ elements and hence the
induction hypothesis applies to give
  $$
  \bigvee\limits_{z\in Z}\irep(\tau_h^* z \tau_h) \geq
\irep(\tau_h^*\tau_h),
  $$
  which translates into 
  $$
  \bigvee_{z\in Z}\orep_h^*\irep(z)\orep_h\geq
  \orep_h^*\orep_h.
  $$
  If this is left-multiplied by $\orep_h$, and right-multiplied by
$\orep_h^*$, we get
  $$
  \bigvee_{z\in Z}\orep_h\orep_h^*\irep(z)\orep_h \orep_h^*\geq
  \orep_h\orep_h^*\orep_h\orep_h^* =
  \orep_h\orep_h^*,
  $$
  which means that
  $$
  \orep_h\orep_h^* \leq \bigvee_{z\in Z}\irep(z).
  \eqmark GoingCrazyOne
  $$
  Leaving this aside for a moment consider the domain 
  $
  D = A\setminus \union\nolimits_{i=1}^m D_i,
  $
  and let $K$ be the set of all $h_i$'s belonging to $D$.  So far we
have been discussing several covers in the sense of semilattices
\lcite{\DefineCoverInSLat}, but now we claim that $K$ is a cover for
$D$, in the sense of semigroupoids \lcite{\DefineCoverSGPD}.  To see
this let $g\in D$, and notice that since $g\in A$, we have that $(g,\D
g,g)\leq (1,A,1)$.  It follows that there is some $z\in Z$ such that
$z\its(g,\D g,g)$, but notice that such a $z$ may not be in $Z''$,
since
  $$
  (g,\D g,g)\perp (1,D_i,1),
  $$
  by \lcite{\LessOrOrhogoNew}, because $g\notin D_i$.  Therefore
$(h_i,C_i,h_i)(g,\D g,g)\neq 0$, for some $i\leq n$.  In particular
this implies that $h_iu=gv$, for some $u,v\in \Mu$. Since $g\in D$, we
have by \lcite{\PropOfDomain} that $h_i\in D$, so in fact $h_i\in K$.
This proves our claim that $K$ is a cover for $D$, so
  $$
  \iniq(D) \={(\ExtendedTight)}
  \bigvee_{h\in K}\orep_h\orep_h^*
  \buildrel {(\GoingCrazyOne)}\over \leq
  \bigvee_{z\in Z}\irep(z).
  \eqmark GoingCrazyTwo
  $$
  Since $A\c D \cup \union\nolimits_{i=1}^m D_i$, we have
  $$
  \irep(1,A,1) = 
  \iniq(A) \leq 
  \iniq(D) \vee \bigvee_{i=1}^m \iniq(D_i) 
  \buildrel {(\GoingCrazyTwo)}\over \leq
  \bigvee_{z\in Z}\irep(z) \vee \bigvee_{i=1}^m \irep(1,D_i,1) =
  \bigvee_{z\in Z}\irep(z),
  $$
  proving \lcite{\BigGoalTight} for $f=1$, and arbitrary $n$.
Summarizing:

\state Lemma
  \label OKforCoversOneAOne
  If $A\in \Q$ and $Z$ is a cover for $(1,A,1)$ then
  $$
  \bigvee_{z\in Z}\irep(z) \geq\irep(1,A,1).
  $$

Let us now face \lcite{\BigGoalTight} in the most general situation,
so we assume that  $(f,A,f)$ is an arbitrary element of $\SM$ and that 
  $Z= \{(h_i,C_i,h_i)\}_{i=1}^n$
  is a cover for $(f,A,f)$.

Since $(h_i,C_i,h_i)\leq (f,A,f)$, for every $i$, we have by
\lcite{\OrderInISGforGPD} that $f\dil h_i$, so we may write
$h_i=fg_i$, with $g_i\in\Mu$, and in addition we have that $C_i\c
g_i\inv(A)$.  Observe that
  $$
  C_i\c\D {h_i} = \D{fg_i} \c \D {g_i},
  $$
  so $(g_i,C_i,g_i)\in \SM$.  Notice that
  $$
  (g_i,C_i,g_i) (1,A,1) = (g_i,C_i\cap g_i\inv(A),g_i) = (g_i,C_i,g_i),
  $$
  so $(g_i,C_i,g_i)\leq (1,A,1)$.  We then claim that
$\{(g_i,C_i,g_i)\}_{1=1}^n$ is a cover for $(1,A,1)$.  In order to
prove it let $(k,B,k)$ be a nonzero element with $(k,B,k)\leq
(1,A,1)$, so $B\c k\inv(A)$.  Given that $B$ is nonempty, the same is
true for $k\inv(A)$, so \lcite{\InverseImage} applies and gives $k\in
A\cup\{1\}\c\D f\cup\{1\}$, so $fk$ is a well defined element of
$\Mu$.
  Since $A\c\D f = f\inv(\M)$, we have
  $$
  B\c k\inv(A) \c k\inv(f\inv(\M)) \={(\RulesInverImag.ii)}
(fk)\inv(\M) = \D{fk},
  $$
  hence $(fk,B,fk)\in\SM$, and clearly $(fk,B,fk)\leq (f,A,f)$.
  Therefore there exists $i\leq n$ such that $(fk,B,fk)\its
(h_i,C_i,h_i)$, which may be interpreted via \lcite{\ConditionForIts},
by saying that there are $x,y\in\Mu$, such that 
  $$
  fkx = h_iy = fg_iy
  \and x\inv(B)\cap y\inv(C_i)\neq\emptyset.
  $$
  Since $f$ is monic we conclude that $kx = g_iy$ and, again by
\lcite{\ConditionForIts}, that $(k,B,k)\its (g_i,C_i,g_i)$, proving
our claim.  Therefore
  $$
  \bigvee_{i=1}^n \orep_{g_i} \iniq(C_i)\orep_{g_i}^* =
  \bigvee_{i=1}^n \irep(g_i,C_i,g_i) 
  \buildrel {(\OKforCoversOneAOne)} \over \geq
  \irep(1,A,1) = \iniq(A),
  $$
  which leads to 
  $$
  \irep(f,A,f) =
  \orep_f \iniq(A) \orep_f^* \leq
  \orep_f \Big(\bigvee_{i=1}^n \orep_{g_i} \iniq(C_i)\orep_{g_i}^*
\Big) \orep_f^* =
  \bigvee_{i=1}^n \orep_f\orep_{g_i} \iniq(C_i)\orep_{g_i}^*\orep_f^*
\$=
  \bigvee_{i=1}^n \orep_{h_i} \iniq(C_i)\orep_{h_i}^* =
  \bigvee_{i=1}^n \irep(h_i,C_i,h_i) =
  \bigvee_{z\in Z} \irep(z),
  $$
  and we are done!


  \def\OM{{\cal O}_\M}
  \def\OMlcm{{\cal O}^{{\rm lcm}}_\M}
  \def\nice{normal}%
  \def\orepu{\orep^u}
  \def\irepu{\irep^u}
  \def\ESMt{\widehat {E(\SM)}_\tightbox}

\section{The C*-algebra of a semigroupoid}
  In this section we fix a countable semigroupoid $\M$ without springs,
in which every element is monic, and such that every intersecting pair
of elements admits a least common multiple.  Our goal will be to study
the universal C*-algebra for representations of $\M$.

To single out the special kind of representations of $\M$ which we
will focus on we give the following:

\definition
   A representation $\orep$ of $\M$ on a Hilbert space $H$ will be
called
  \stress{\nice}, provided it is tight and respects least common
multiples.

  Recall from \cite{\ExSemiGpdAlg} that the C*-algebra of $\M$, denoted $\OM$,
is the C*-algebra generated by a universal tight representation of
$\M$.  By definition we therefore see that *-representations of $\OM$
correspond bijectively to tight representations of $\M$.  If $\M$
satisfies the hypothesis of \lcite{\TightLCM}, then it is automatic
that the tight representations we are talking about respect least
common multiples, and hence are {\nice} representations.

It is not clear to me whether or not one really needs the hypothesis
of \lcite{\TightLCM} to obtain that conclusion but, given the
dependence of our previous results on least common multiples, we simply cannot
live without it.  So much so that we are willing to impose it from the
outside:

\definition
  \label DefineOMlcm
  We will denote by $\OMlcm$ the C*-algebra generated by the range of
a universal {\nice} representation $\orepu$ of $\M$ (such as the
direct sum of all {\nice} representations of $\M$ on subspaces of
Hilbert's space $l_2$).

If $\M$ satisfies the hypothesis of \lcite{\TightLCM} it is then
obvious that $\OMlcm=\OM$, but in general all we can say is that
$\OMlcm$ is a quotient of $\OM$.

Observe that $\orepu$ is not necessarily injective by
\lcite{\ThouShallNotSeparate}.  The reader is referred to
\cite{\ExelAlgebra} for a thorough treatment of the injectivity
question.

Restricting our attention to $\OMlcm$ we therefore see that its
*-representations correspond bijectively to {\nice} representations of
$\M$, and in view of \lcite{\FromRepISGtoGPOD} and
\lcite{\MainCorrTheorem}, they also correspond bijectively to tight
representations of $\SM$.  Furthermore these correspond to
representations of the C*-algebra of the groupoid described in
\lcite{\EquivTightRepISGandGPD}.  

\definition
  We will denote by $\G_\M$ the
the groupoid of germs associated to the restriction of the action
$\act$ of \lcite{\MovingSemicharacters.iv} to the tight part of the
spectrum of $E(\SM)$.

  We thus arrive at one of the main results of this work:

\state Theorem
  \label CstarAlgOfGpd
  Let $\M$ be a countable semigroupoid with no springs, in which every
element is monic, and such that every intersecting pair of elements
admits a least common multiple.  Then $\OMlcm$ is naturally isomorphic
to $C^*(\G_\M)$.

\proof Let $X = \ESMt$, as defined in \lcite{\AltSpectrums}.
By \lcite{\RepXIsRep} we have that the map
  $$
  \irepu: s\in\SM \mapsto \iota(1^X_{\p s}\d_s )\in C^*(\G_\M),
  $$
  is a representation of $\SM$ in $C^*(\G_\M)$ (assumed to be an
algebra of operators via any faithful non-degenerated representation),
which is supported in $\ESMt$, and hence is tight by
\lcite{\CharactTightReps}.  The superscript ``$u$" in $\irepu$ is
justified by its universal property \lcite{\FromRepISGtoGroupoid}.

Employing \lcite{\FromRepISGtoGPOD} we deduce that the composition
  $$
  \orep : \M \labelarrow{\tau} \SM \labelarrow{\irepu} C^*(\G_\M)
  $$
  is a tight representation of $\M$ which respects least common
multiples, i.e., $\orep$ is a {\nice} representation.  Invoking the
universal property of $\OMlcm$ there exists a *-homomorphism
  $$
  \phi:\OMlcm \to C^*(\G_\M),
  $$
  such that the diagram
  
  \beginmypicture
  \setcoordinatesystem units <0.0010truecm, -0.0013truecm> point at 0
0
  \pouquinho = 600

  \put {$\SM$} at 200 -1500
  \morph {-1400}{0100}{0000}{-1500} \arwlabel{\tau}{-1000}{-0800}
  \morph {0400}{-1500}{1800}{-0100} \arwlabel{\irepu}{1400}{-1000}
  \morph {-1500}{0000}{1500}{0000} \arwlabel{\orep}{0000}{-200}
  \put {$\M$} at -1500 0000
  \put {$C^*(\G_\M)$} at 1800 0000
  \put {$\OMlcm$} at 0300 1500
  \morph {-1400}{0100}{0000}{1400} \arwlabel{\orepu}{-1000}{0900}
  \morph {0600}{1400}{1800}{0100} \arwlabel{\phi}{1400}{0900}
  \endmypicture
  \bigskip

\noindent commutes, where $\orepu$ was defined in \lcite{\DefineOMlcm}.
To define an inverse for $\phi$ recall that $\orepu$ is a {\nice}
representation of $\M$, and hence by \lcite{\MainCorrTheorem} there
exists a tight representation $\irep$ of $\SM$ such that $\orepu =
\irep\circ \tau$.
  The space of $\irep$ is evidently the same as the space $H^u$ of
$\orepu$. Since $\tau(\M)$ generates $\SM$, by
\lcite{\LambdaGenerates}, we deduce that the range of $\irep$ is
contained in the inverse semigroup of partial isometries on $H^u$
generated by the range of $\orepu$, which is obviously contained in
$\OMlcm$.  We may therefore regard $\irep$ as a map from $\SM$ to
$\OMlcm$.

   Being tight, $\irep$ is supported in $\ESMt$ by
\lcite{\CharactTightReps} and hence we may use
\lcite{\FromRepISGtoGroupoid} to conclude that there exists a
*-representation $\rho$ of $C^*(\G_\M)$ on $H^u$, such that
$\rho\circ\irepu = \irep$.  The range of $\irepu$ may be shown to
generate $C^*(\G_\M)$ as a C*-algebra, and hence we conclude as above
that $\rho$ may be regarded as a map from $C^*(\G_\M)$ to $\OMlcm$.

\bigskip
  \beginmypicture
  \setcoordinatesystem units <0.0013truecm, -0.0011truecm> point at 0
0
  \pouquinho = 600

  \put {$\SM$} at 200 -1500
  \morph {-1400}{0100}{0000}{-1500} \arwlabel{\tau}{-1000}{-0800}
  \morph {0400}{-1500}{1800}{-0100} \arwlabel{\irepu}{1400}{-1000}

 \morph {0200}{-1600}{0200}{1600} \arwlabel{\irep}{0400}{0000}

  \put {$\M$} at -1500 0000
  \put {$C^*(\G_\M)$} at 1800 0000
  \put {$\OMlcm$} at 0300 1500
  \morph {-1400}{0100}{0000}{1400} \arwlabel{\orepu}{-1000}{0900}
  \morph {1800}{0100}{0600}{1400} \arwlabel{\rho}{1400}{0900}
  \endmypicture
  \bigskip

We therefore have 
  $$
  \rho\circ\phi\circ\orepu =
  \rho\circ\orep =
  \rho\circ\irepu\circ \tau =
  \irep\circ \tau =
  \orepu,
  $$
  so $\rho\circ\phi$ coincides with the identity on the range of
$\orepu$, which is known to generate $\OMlcm$.  This proves that
$\rho\circ\phi$ is the identity map.  On the other hand
  $$
  \phi\circ\rho\circ\orep =
  \phi\circ\orepu =
  \orep,
  $$
  so 
  $\phi\circ\rho$ coincides with the identity on the range of $\orep$,
which again generates $C^*(\G_\M)$.  Therefore
  $\phi\circ\rho$ is the identity map, proving that $\phi$ and $\rho$
are each other's inverse, and hence isomorphisms.
  \proofend

It is interesting to notice that since $\tau(\M)$ generates $\SM$, the
groupoid $\G_\M$, which consists of germs for the action of $\SM$ on
$\ESMt$, is also in a sense generated by the action of $\M$, via
$\tau$.

A concrete understanding of this groupoid clearly depends on the
ability to describe $\ESMt$ in clear terms.  That is the purpose of
our next section.


\section{Categorical semigroupoids}
  \label CategoricalSec
  In this section we will give a concrete description for the space of
tight characters on the idempotent semilattice of $\SM$, where $\M$ is
a semigroupoid.  To reduce the technical difficulties to a minimum we
will assume that $\M$ possesses a crucial property well known to hold
on categories.

\definition
  \label DefineCatSGPD
  A semigroupoid $\M$ is said to be \stress{categorical}, if for every
$f,g\in\M$ one has that $\D f$ and $\D g$ are either equal or
disjoint.

With this notion we wish to capture the essential characteristic of
categories which is relevant to our work.  Obviously any small
category is a categorical semigroupoid.

In order to apply our results to a categorical semigroupoid $\M$ we
must assume that it satisfies our crucial working hypotheses, namely
the conditions listed in \lcite{\StandingHyp}, often adding the
absence of springs.
    With respect to the requirement that every element is monic we
should stress that, although the term we use is inspired in the Theory
of Categories, our use of it is strictly different.  In particular,
requiring an element $f$ to be monic impedes the existence of a right
unit to $f$, namely an element $u$ such that $fu=f$.  According to
Definition \lcite{\DefineMonic}, the only element $u$ which is allowed
to satisfy such an equation is the added unit $1$, as in
$\Mu=\M\,\dot\cup\,\{1\}$.

If one is to apply our theory to a classical small category, one
should therefore first remove all of its identities, and then hope
that the products of the remaining elements never come out to being an
identity.  See below for a discussion of this issue in the context of
higher rank graphs.

\state Proposition
  \label SGPDFromCatego
  Let $\Catego$ be a small category such that no morphism is
right-invertible, except for the identities.  Then the set $\M$ of all
non-identity morphisms admits the structure of a categorical
semigroupoid.  In addition:
  \izitem
  \zitem If every morphism in $\Catego$ is a monomorphism (in the
usual sense of the word), then every element of $\M$ is monic (in the
sense of Definition \lcite{\DefineMonic}).
  \zitem If for every object $v$ in $\Catego$ there exists a morphism
$f\neq id_v$, such that $\r(f)=v$, then $\M$ has no springs.
  \zitem Suppose that whenever $fu=gv$ in $\Catego$, there exists a
pull-back for the pair $(f,g)$.  Then $\M$ admits least common
multiples.

\proof
  Given $f,g\in\M$ suppose that $fg$ is an identity morphism,
necessarily the identity on $v:=\r(f)$.  Then $f$ is right-invertible
and hence by hypothesis, $f=id_v\notin\M$, a contradiction.  So
whenever $f,g\in\M$, and $fg$ is defined in $\Catego$, one has that
$fg\in\M$.  We may then put
  $$
  \Mt = \{(f,g)\in\M\times\M: \s(f)=\r(g)\},
  $$
  and it is clear that $\M$ is a categorical semigroupoid with
composition as multiplication.

Under hypothesis (i) suppose that $fg=fh$, for $f\in\M$, and
$g,h\in\Mu$.  If $g,h\in\M$, we have that $g=h$ because $f$ is a
monomorphism, by hypothesis.  If $g\in\M$ and $h=1$, then
  $$
  fg = f = f\,id_{\s(f)}.
  $$
  Using again that $f$ is monic we deduce that 
  $g =id_{\s(f)}\notin\M$, a contradiction.  This shows that every
element is monic in the sense of \lcite{\DefineMonic}.

Point (ii) is elementary.  With respect to (iii) let $f,g\in\M$ be
such that $f\its g$.  If $g\dil f$ then it is obvious that
$f=\lcm(f,g)$, and similarly $g=\lcm(f,g)$, if $f\dil g$.
  Otherwise, assuming that neither $g\dil f$, nor $f\dil g$, there are
$u$ and $v$ in $\M$ (as opposed to $\Mu$) such that $fu=gv$.
  So there let $(p,q)$ be a pull back for $(f,g)$, which in particular
entails $fp=gq$.
  Since $f$ and $g$ do not divide each other we have that neither $p$
nor $q$ are identity morphisms.
  Setting $m=fp$, notice that $m$ is not an identity either because
$f$ is not right-invertible, and so $m\in\M$. It is then clear that
$m$ is a common multiple of $f$ and $g$ (relative to $\M$).  If
$n\in\M$ is another common multiple of $f$ and $g$, then $n=fx=gy$,
for some $x,y\in\Mu$.  But, since $f$ and $g$ do not divide each other
we see that $x,y\in\M$, and hence the equation $fx=gy$ makes sense in
$\Catego$.  By definition of pull-backs, there is a morphism $r$ such
that $x=pr$, and $y=qr$.  Therefore $n=fx = fpr = mr$, and hence
$m\dil n$ in $\M$, regardless of whether or not $r\in\M$.
  \proofend

\cryout{From now on we assume that $\M$ is a fixed categorical semigroupoid
satisfying \lcite{\StandingHyp}, and having no springs.}

  The greatest
simplification brought about by restricting one's attention to such
semigroupoids is in the structure of the semilattice $\Q$ of
elementary domains defined in \lcite{\DefineMiddleSets}, which is
easily seen to be just
  $$
  \Q=\{\D f: f\in \M\} \cup\{\emptyset\}.
  $$

If $f$ lies in $\D {g_1}$ and $\D {g_2}$, for two elements
$g_1,g_2\in\M$, then evidently $\D {g_1}\cap\D {g_2}$ is nonempty, and
hence by hypothesis $\D {g_1}=\D {g_2}$.  We may then define the
\stress{range} of $f$, denoted
  $$
  \r(f),
  $$
  to be the only element $A\in \Q$ for which $f\in A$.  It is possible
that some $f\in\M$ is not in any $\D g$, in which case $\r(f)$ will
not be defined. We then conclude that
  $$
  (f,g)\in\Mt \iff \D f = \r(g)
  \for f,g\in\M,
  $$
  where we consider the expression in the right-hand side to be false
if $\r(g)$ is not defined.  The reader is invited to compare this with
the criteria for two morphisms in a category to be composable.

The above simple form of $\Q$ leads to a simplified $\SM$, which then
consists of the disjoint union of the following sets:
  $$
  \big\{(f,A,g): f,g\in\M,\ \D f = \D g = A\big\}, $$ $$
  \big\{(f,\D f,1): f\in \M\big\}, $$ $$
  \big\{(1,\D g,g): g\in \M\big\}, $$ $$
  \hbox to 19 pt {\hfill} \big\{(1,\D f,1): f\in \M\big\}, \hbox to 19
pt { and} $$ $$
  \{0\}.
  $$
  The all important semilattice $E\big(\SM\big)$ is then simply the disjoint
union of the   sets%
  \fn{%
  In case $\M$ is obtained from a category $\Catego$, as in
\lcite{\SGPDFromCatego}, then $\Q$ is in one-to-one correspondence
with the objects in $\Catego$, or at least those which are the
co-domain of a non-identity morphism.  It is therefore curious that,
after the identities have been put to sleep in
\lcite{\SGPDFromCatego}, they were suddenly awakened by this
expression for $E(\SM)$.}
  $$
  E(\SM) =
  \big\{(f,\D f,f): f\in \M\big\} \cup
  \big\{(1,A,1): A\in\Q\}.
  $$

\sysstate{Notations}{\rm}{From now on we shall adopt the following
shorthand notations:
  \izitem
  \zitem $\E=E\big(\SM\big)$,
  \zitem $\E_p=\big\{(f,\D f,f): f\in \M\big\}$,
  \zitem $\E_q=\big\{(1,A,1): A\in\Q\}$,
  \zitem $\fin_f = \fin^\tau_f = (f,\D f,f)$, for all $f\in\M$,
  \zitem $\ini_A = (1,A,1)$, for all $A\in\Q$.
  }

This in turn evokes the notations
  $\spec$ from \lcite{\BasicSpectrum},
  $\specz$ from \lcite{\DefineCharacter}, 
  in addition to $\speci$ and $\spect$ from \lcite{\AltSpectrums}.
     It is our purpose here to describe the most important of these,
namely $\spect$.  Recall from \lcite{\ClosureOfUltrafilters} that
$\spect$ is the closure of $\speci$ in $\specz$.  Being left out of
this equation, $\spec$ will not matter much to us.

The following is a compilation of properties relating to the order
relation on $\E$, some of which we have already encountered in
\lcite{\OrderInISGforGPD} and \lcite{\LessOrOrhogoNew}, and which
completely describes the structure of $\E$, as a semilattice.

\state Proposition
  \label BasicPFQA
  If $f,g\in\M$, and $A,B\in \Q$. Then
  \izitem
  \zitem $\fin_f\fin_g=\fin_{\lcm(f,g)}$, if $f\its g$,
  \zitem $\fin_f\leq\fin_g$, if and only if $g\dil f$,
  \zitem $\fin_f\perp\fin_g$, if $f\disj g$,
  \zitem $\fin_f\leq \ini_A$, if $f\in A$,
  \zitem $\fin_f\perp\ini_A$, if $f\notin A$,
  \zitem $\ini_A\perp\ini_B$, if $A\neq B$.

\definition
  Let $\xi$ be a filter in $E$.  We will say that $\xi$ is of 
  \izitem
  \zitem $p$-type, if $\xi\c\E_p$,
  \zitem $q$-type, if $\xi\c\E_q$,
  \zitem $pq$-type, if $\xi\cap\E_p$, and $\xi\cap\E_q$ are nonempty.

If $\xi$ is a filter of $q$-type then all of its elements are of the
form $\ini_A$, for some nonempty $A\in \Q$.  But since any two of
these are disjoint by \lcite{\BasicPFQA.vi}, only one such element is
allowed.  We thus see that $\xi=\{\ini_A\}$, for a single nonempty
$A\in \Q$.  On the other hand, given any $A\in\Q$, with
$A\neq\emptyset$, it is easy to see that the singleton
  $$
  \xi_A=\{\ini_A\}
  \eqmark FilterXiA
  $$
  is a filter of $q$-type.

 The next concept is borrowed from \scite{\infinoa}{5.5}.

\definition
  Let $\xi$ be a filter in $\E$. We will say that the \stress{stem} of
$\xi$ is the set
  $$
  \stem_\xi = \{f\in\M: \fin_f\in\xi\}.
  $$

It is clear that a filter is of $q$-type if and only if its stem is
empty.
  The following elementary result describes all filters according to
their type:

\state Proposition
  \label FiltersAndTypes
  Let $\xi$  be  a filter in $\E$.
  \izitem
  \zitem If $\,\xi$ is of $q$-type, then 
  $\xi= \xi_A :=\{\ini_A\}$, for some $A\in\Q$, with $A\neq\emptyset$.
  \zitem If $\,\xi$ is of $p$-type, then
  $\xi = \{\fin_f: f\in\stem_\xi\}$, and   moreover $\r(f)$ is not defined
for any $f\in\stem_\xi$.
  \zitem If $\,\xi$ is of $pq$-type, then there is some
$A\in\Q$, such that $\stem_\xi\c A$.  In addition
  $\xi = \{\fin_f: f\in\stem_\xi\}\cup\{\ini_A\}$.

\proof
  Point (i) was already discussed above.  Under the hypothesis of  (ii),
suppose that $f$ is an element of $\stem_\xi$ such that $f\in A$, for
some $A\in\Q$.  Then $\fin_f\leq\ini_A$, by \lcite{\BasicPFQA.iv} and
hence $\ini_A\in\xi$, by \lcite{\DefineFilter.ii}.  This contradicts
the fact that $\xi$ is of $p$-type, and hence $f$ does not belong to
any $A$, which means that $\r(f)$ is not defined.  The first sentence
of (ii) is obvious.

If $\xi$ is of $pq$-type, then by assumption $\xi$ contains some
$\ini_A$, for $A\in\Q$.  As already argued, only one such element is
allowed and hence $\xi\cap\E_q$ must be a singleton $\{\ini_A\}$.  The
other elements of $\xi$ must be of the form $\fin_f$, for $f\in\M$,
and hence $\xi = \{\fin_f: f\in\stem_\xi\}\cup\{\ini_A\}$.  Given any
$f\in\stem_\xi$, we have that both $\fin_f$ and $\ini_A$ lie in $\xi$,
and hence $\fin_f\its\ini_A$, by \lcite{\DefineFilter.iii}.  It then
follows from \lcite{\BasicPFQA.iv} that $f\in A$.
  \proofend

\state Proposition
  \label PropertiesOfStem
  Given a filter $\xi$ on $\E$ one has that:
  \izitem
  \zitem if $f\in\stem_\xi$ and $g\in\M$ is such that $g\dil f$, then
$g\in\stem_\xi$,
  \zitem for every $f,g\in \stem_\xi$, one has that $f\its g$, and moreover
$\lcm(f,g)\in\stem_\xi$.

\proof
  If \ $g\dil f\in \stem_\xi$, \ we have that $\xi\ni\fin_f\leq
\fin_g$, so $\fin_g\in\xi$, and hence $g\in\stem_\xi$.
  In order to prove (ii) let us be given $f,g\in \stem_\xi$ and
suppose by contradiction that $f\disj g$.  Then 
  $$
  0=\fin_f\fin_g\in\xi,
  $$
  which is impossible. This proves that $f\its g$.  Moreover,
  $$
  \xi \ni \fin_f\fin_g = \fin_{\lcm(f,g)},
  $$
  so $\lcm(f,g)\in\stem_\xi$.
  \proofend

Based on the findings of the above result we introduce the following
generalization of the notion of paths in a graph:

\definition
  \label DefinePath
  A \stress{path} in a semigroupoid $\M$ is a subset $\stem\in\M$
such that, 
  \izitem
  \zitem if $f\in\stem$, and $g\in\M$ is such that $g\dil f$, then
$g\in\stem$,
  \zitem for every $f,g\in \stem$, one has that $f\its g$, and
moreover $\lcm(f,g)\in\stem$.
  \medskip\noindent
  An \stress{ultra-path} is a path which is not properly contained in
any other path.

It is therefore obvious that $\stem_\xi$ is a path for every filter
$\xi$, possibly the empty path if $\xi$ is of $q$-type.

  Given a nonempty path $\stem$, suppose that $f,g\in\stem$ and that
$f\in A$, for some $A\in\Q$.  By \lcite{\DefinePath.ii} we may write
$fu=gv$, for suitable $u,v\in\Mu$, and hence $g\in A$, by
\lcite{\DfgDgDois}.  This means that, if $\r(f)$ is defined for some
$f\in\stem$, then $\r(g)=\r(f)$ for every $g\in\stem$.  In this case
we say that $A$ is the range of $\stem$, in symbols
  $$
  \r(\stem) = A.
  $$
  Otherwise $\r(\stem)$ is not defined.
  
  Notice that if we are given some  $f\in\M$ then
  $$
  \stem^f:= \{g\in\M: g\dil f\}
  $$
  is clearly a path.  By Zorn's Lemma there exists an ultra-path
containing $\stem^f$, and hence any element of $\M$ belongs to some
ultra-path.  Another consequence of this is that even if the definition
allows for paths to be empty, the empty path is never an ultra-path
(unless $\M=\emptyset$).

A filter of the form $\xi_A$, as defined in
\lcite{\FiltersAndTypes.i}, is never an ultra-filter because if $f\in
A$, then the set of all elements in $\E$ which are bigger than or
equal to $\fin_f$ forms a filter properly containing $\xi_A$.  For
that reason the filters of $q$-type are left out of the following
characterization of ultra-filters.

\state Proposition
  The correspondence \ $\xi\mapsto \stem_\xi$ \ is a bijection from
the set 
  of all filters in $\E$, bar the $q$-types,
  and the set of all nonempty paths.  This also gives a one-to-one
correspondence from the set of all ultra-filters to the set of all
ultra-paths.

\proof Given a nonempty path $\stem$ consider the subset $\xi_\stem$
of $\E$ defined by
  $$
  \xi_\stem  = \left\{\matrix{
    \quad \{\fin_f : f\in\stem\} \cup\{\ini_{\r(\stem)}\},& \hbox{if
$\r(\stem)$ is defined,}\cr\cr
    \{\fin_f : f\in\stem\}, & \hbox{otherwise.}\hfill
  }\right.
  $$
  It is then easy to see that $\xi_\stem$ is a filter and that the
resulting map \ $\stem\to\xi_\stem$ \ gives the inverse of the
correspondence in the statement.  Since the two correspondences
referred to preserve inclusion, it is clear that ultra-filters
correspond to ultra-paths.
  \proofend

From now on we will use \lcite{\XiPhi} and \lcite{\PhiXi} to identify
characters with filters, without further warnings.  Therefore $\specz$
will be seen as the set of all filters in $\E$.  This said, $\speci$
corresponds to ultra-filters, and the filters corresponding to the
elements of $\spect$ will be referred to as \stress{tight-filters}.

Were we only interested in $\speci$, it would be sensible to use the
above result to replace $\speci$ by the set of all ultra-paths.
However our primary interest is in tight filters, and unfortunately
paths fail to capture the topological complexity of filters.

\state Proposition
  \label CharacTightInCatego
  Let $\xi$ be a filter in $\E$.
  \iaitem
  \aitem Suppose that $\xi$ is of $q$-type, and write $\xi=\xi_A$, as
in \lcite{\FilterXiA}.  Then $\xi$ is tight if and only if $A$ admits
no finite cover (in the sense of \lcite{\DefineCoverSGPD}),
  \aitem Suppose that $\xi$ is of $p$-type. Then $\xi$ is tight if and
only if for every $f\in\stem_\xi$, and every finite cover $H$ for $\D
f$ (again in the sense of \lcite{\DefineCoverSGPD}), there is some
$h\in H$ such that $fh\in\stem_\xi$.
  \aitem Suppose that $\xi$ is of $pq$-type.  Then $\xi$ is tight if
and only if the condition in (b) is satisfied and moreover for every
finite cover (ditto) $H$ of $\r(\stem_\xi)$, one has that
$h\in\stem_\xi$, for some $h\in H$.

\proof 
  Before we begin it is convenient to notice the following auxiliary
result: if $A\in\Q$ is nonempty then the covers of $\ini_A$ (in the
sense of \lcite{\DefineCoverInSLat}) which do not contain $\ini_A$
itself, correspond to the covers of $A$ (in the sense of
\lcite{\DefineCoverSGPD}) in the following way: given a cover $H$ of
$A$, the set $\{\fin_h: h\in A\}$ is a cover for $\ini_A$.  On the
other hand, given a cover $Z$ of $\ini_A$ which does not contain
$\ini_A$, the set $\{h\in\M: \fin_h\in Z\}$ is a cover for $A$.

To prove it let $H$ be a cover for $A$.  Then every nonzero element
$z\in\E$, which is smaller than $\ini_A$ is either $\ini_A$ itself, in
which case $z$ intercepts every $\fin_h$, or $z = \fin_g$, for some
$g\in A$ by \lcite{\BasicPFQA}.  In the latter case $g\its h$, for
some $h\in H$, and hence
  $$
  z\fin_h = \fin_g\fin_h = \fin_{\lcm(g,h)} \neq 0,
  $$
  so $z\its \fin_h$.  Conversely, assuming that $Z$ is a cover for
$\ini_A$ not containing $\ini_A$, we have by \lcite{\BasicPFQA} that
$Z$ must have the form
  $$
  Z = \{\fin_{h}: h\in H\},
  $$
  where $H$ is a subset of $A$.  To prove that $H$ is a cover for $A$,
let $f\in A$.  Then $\fin_f\leq\ini_A$ by \lcite{\BasicPFQA.iv}, so
$\fin_f\fin_h$ is nonzero for some $h\in H$, which means that $f\its
h$.

Addressing (i) suppose that $\xi_A$ is a tight filter. Arguing by
contradiction let $H$ be a finite cover for $A$, so that
$\{\fin_h:h\in H\}$ is a cover for $\ini_A$.  If $\phi$ is the tight
character associated to $\xi$ according to \lcite{\PhiXi}, then
  $$
  1  = \phi(\ini_A) =  \bigvee_{h\in H} \phi(\fin_h),
  $$
  and hence $\fin_h\in\xi$, for some $h\in H$.  This would seem to
indicate that $h\in\stem_\xi$, which is a contradiction.  Thus no
cover for $A$ may exist.

Conversely suppose that $A$ admits no finite cover.
  Again denoting by $\phi$ the associated character, as in
\lcite{\PhiXi}, notice that $\phi(\ini_A)=1$, and hence condition
\lcite{\ChaeckWithNonEmptX.i} is satisfied so we may use
\lcite{\AltLatTightRep} in order to prove that $\phi$ is tight.  So
let $x\in\E$ and let $Z$ be a finite cover for $x$.  We must then
prove that
  $$
  \bigvee_{z\in Z}\phi(z) \geq \phi(x).
  \subeqmark CharacTightInCategoGoal
  $$
  Observe that, except for $x=\ini_A$, one has that $\phi(x)=0$, in
which case the above inequality holds trivially.  We may then restrict
our attention to the case in which $x=\ini_A$.  Excluding the trivial
case in which $\ini_A$ itself belongs to $Z$, we have that
  $
  Z = \{\fin_{h}: h\in H\},
  $
  where $H$ is a finite cover for $A$, but since $A$ admits no finite
cover by hypothesis, there is nothing to be proven.

Suppose now that $\xi$ is a tight filter of $p$-type or $pq$-type,
which implies that $\stem_\xi$ is nonempty.
  Let $f\in\stem_\xi$ and let us be given a finite cover $H$ for $\D
f$.  This time we claim that $\{\fin_{fh}:h\in H\}$ is a cover for
$\fin_f$.  In fact, any element of $\E$ which is smaller than $\fin_f$
is necessarily of the form $\fin_g$, for some $g\in\M$ which is a
multiple of $f$, by \lcite{\BasicPFQA}.  So write $g=fk$, with
$k\in\Mu$.  If $k=1$ then obviously $\fin_g\its\fin_{fh}$, for all
$h\in H$.  Otherwise $k\in\D f$, and hence $h\its k$, for some $h\in
H$.  This implies that $fh\its fk$, or equivalently that $fh\its g$,
whence
  $$
  \fin_{fh}\fin_g = \fin_{\lcm(fh,g)} \neq 0,
  $$
  proving the claim.  Because the character $\phi$ associated to $\xi$
is tight we deduce that
  $$
  1 = \phi(\fin_f) = \bigvee_{h\in H} \phi(\fin_{fh}),
  $$
  so there exists some $h\in H$, such that $\fin_{fh}\in\xi$, and
hence $fh\in\stem_\xi$.
  In the special case in which $\xi$ is a tight filter of $pq$-type we
must still address the last assertion in (c).  Let $A=\r(\stem_\xi)$,
and picking any $f\in\stem_\xi$ we have that $f\in A$, so
  $$
  \xi\ni \fin_f \leq\ini_A,
  $$
  whence $\ini_A\in\xi$, which is to say that the associated character
$\phi$ satisfies $\phi(\ini_A)=1$.
  Let $H$ be a finite cover for $A$.  By the auxiliary result proved
above one has that $\{\fin_h:h\in H\}$ is a cover for $\ini_A$ and
hence
  $$
  1  = \phi(\ini_A) =  \bigvee_{h\in H} \phi(\fin_h),
  $$
  from where we deduce that $\phi(\fin_h)=1$, for some $h\in H$,
meaning that $h\in\stem_\xi$.

Let us now address the converse implications in (b) and (c)
simultaneously.  So let $\xi$ be a filter with nonempty stem
satisfying the condition in (b).  In case $\xi$ is of $pq$-type, we
assume in addition that it also satisfies the condition in (c).

  Since the associated character $\phi$ is nonzero, there must exist
some $x\in\E$, such that $\phi(x)=1$, and hence we may again use
\lcite{\AltLatTightRep} in order to prove that $\phi$ is tight.  So
let $x\in\E$, and let $Z$ be a finite cover for $x$. We must prove
\lcite{\CharacTightInCategoGoal}. Excluding the trivial case in which
$\phi(x)=0$, we suppose that $x\in\xi$.

  The proof will be broken up in two cases, the first one
corresponding to $x=\fin_f$, for some $f\in\stem_\xi$.
  Since $Z$ is a cover for $\fin_f$ we have by \lcite{\BasicPFQA.ii}
that,
  $$
  Z = \{\fin_g: g\in G\},
  $$
  where $G$ is a finite subset of $\M$ consisting of multiples of
$f$.  We may therefore rewrite $Z$ as 
  $$
  Z  = \{\fin_{fh}: h\in H\},
  $$
  where $H$ is a finite subset of $\D f\cup\{1\}$.  If $\fin_f$ itself
belongs to $Z$ then obviously the right-hand side of
\lcite{\CharacTightInCategoGoal} is 1, and the proof is finished.  So
assume that $H\c\D f$.  We then claim that $H$ is a cover for $\D f$.
To prove it let $k\in\D f$.  Then $\fin_{fk}\leq\fin_f$, so that
$\fin_{fk}\its\fin_{fh}$, for some $h\in H$, meaning that $fkx=fhy$,
for suitable $x,y\in\Mu$.  Canceling out $f$ we deduce that $kx=hy$,
and hence that $k\its h$, proving our claim.  The hypothesis therefore
applies and we have that $fh\in\stem_\xi$, for some $h\in H$, which
may be rephrased by saying that $\phi(\fin_{fh})=1$, proving that the
left-hand side of \lcite{\CharacTightInCategoGoal} is 1.

Assume next that $x=\ini_A$, for some $A\in\Q$.  As we are supposing
that $\phi(x)=1$, and hence that $\ini_A\in\xi$, this can only happen
if $\xi$ is of $pq$-type, in which case we moreover have that
$A=\r(\stem_\xi)$.  The fact that $Z$ is a cover for $\ini_A$ implies
that either $\ini_A\in Z$, when \lcite{\CharacTightInCategoGoal} is
readily proved, or
  $Z = \{\fin_h : h\in H\}$, where $H$ is a cover for
$A=\r(\stem_\xi)$.
  This may then be combined with our hypothesis to give
$h\in\stem_\xi$, for some $h\in H$.  Then $\fin_h\in\xi$ and hence
$\phi(\fin_h)=1$.  Since $\fin_h$ is in $Z$, we have that the
left-hand side of \lcite{\CharacTightInCategoGoal} is 1, concluding
the proof.
  \proofend


  \def\bp{\bar p}
  \def\bq{\bar q}
  \def\F{\bar f}
  \def\Cat{\Lambda}
  \def\M{{\underline{\Lambda}}}
  \def\rank{{\partial}}

\section{Higher rank graphs}
  In this section we wish to apply our theory to higher rank graphs.
The reader should consult the references listed in the introduction
for more information on this subject.

  From now on we assume that $k\geq 1$ is an integer and $\Cat$ is a
$k$-graph, with \stress{rank} map given by
  $$
  \rank:\Cat \to\N^k.
  $$
  The well known \stress{unique factorization property} states that for every
morphism $f$ in $\Cat$, and for every $n,m\in\N^k$ such that
  $\rank(f) = n+m$,
  there exists a unique pair of morphisms $(g,h)$ such that $f=gh$,
$\rank(g)=n$, and $\rank(h)=m$.

  As usual we will say that $f$ is an \stress{edge} if $\rank(f)$ is an
element of the canonical basis $\{e_i\}_{i=1}^k$ of $\N^k$.  For an
edge $f$, one sometimes refer to $\rank(f)$ as the \stress{color} of $f$.
While one does not really have to attach ``colors" to the $e_i$, it
does make sense to say that two edges have, or do not have the same
color.

The possibility of studying $\Cat$ with our tools naturally hinges on
whether or not we may verify our working hypotheses, namely
\lcite{\StandingHyp}, and the absence of springs.  

We will soon specialize to a situation in which we may apply all of
the points in \lcite{\SGPDFromCatego}, hence obtaining our working
hypotheses.  We do so mainly to avoid technical complications, but we
nevertheless believe that our methods, and Theorem
\lcite{\EquivTightRepISGandGPD} in special, may be applied to the
inverse semigroup constructed in \cite{\Muhly} in the most general
case, obtaining the same description of the C*-algebra of $\Cat$ as a
groupoid C*-algebra.

  Notice that the identities in $\Cat$ are precisely the morphisms
with rank zero.  Moreover if $f,g\in\Cat$ are such that $fg$ is an
identity, then
  $$
  0 = \rank(fg) = \rank(f) + \rank(g),
  $$
  which implies that $\rank(f) = \rank(g)=0$,
  so $f$ and $g$ are both identities.  This says that no morphism
other than the identities may be right-invertible, and hence we have
by \lcite{\SGPDFromCatego} that the set $\M$ of all non-identity
morphisms is a categorical semigroupoid.

With respect to \lcite{\SGPDFromCatego.i}, if $f,g,h$ are morphisms in
$\Cat$ such that $fg=fh$, then $\rank(g)=\rank(h)$, and the uniqueness of
the factorization implies that $g=h$.  This says that every morphism
in $\Cat$ is a monomorphism, so we may apply
\lcite{\SGPDFromCatego.i} to collect another of our working
hypotheses.

  For each vertex (object) $v$ in $\Cat$ and each $n\in\N^k$ one
usually denotes by $\Cat_n^v$ the set of all morphisms $f$ in $\Cat$
with $\r(f) = v$ and $\rank(f)=n$.

Recall that $\Cat$ is said to be \stress{row-finite} if $\Cat_n^v$ is
finite for every $v$ and $n$.  If $\Cat_n^v$ is never empty then one
says that $\Cat$ \stress{has no sources}.  Notice that in order for
the associated semigroupoid $\M$ to have no springs one does not
necessarily need to rule out all sources of $\Cat$.  It is clearly
enough to suppose that
  $$
  \D v = \union_{n\neq 0}\Cat_n^v \neq\emptyset,
  $$
  for every object $v$ in $\Cat$.

The last requirement we will impose on $\Cat$ is designed to allow for
the use of \lcite{\SGPDFromCatego.iii}, and it is related to the
question of finite alignment.  Recall that $\Cat$ is said to be
\stress{finitely aligned}, if for every $f,g\in\Cat$ one has that
  $$
  \Cat^{\rm min}(f,g) :=
  \big\{(p,q)\in\Cat\times\Cat: fp=gq, \hbox{ and } \rank(fp) = \rank(f)\vee
\rank(g)\big\}
  $$
  is finite.

Observe that for any pair $(p,q)$ in 
  $\Cat^{\rm min}(f,g)$, one has that $m:=fp$ is a common multiple of
$f$ and $g$.  If $(p',q')$ is another pair in $\Cat^{\rm min}(f,g)$,
then $m'=fp'$ is another common multiple but neither $m\dil m'$, nor
$m'\dil m$, because $\rank(m) = \rank(m')$.  Thus, unless
  $\Cat^{\rm min}(f,g)$ has at most one element, $\M$ will not admit
least common multiples.

\definition
  We shall say that $\Cat$ is \stress{singly aligned}, if $\Cat^{\rm
min}(f,g)$ has at most one element for every $f$ and $g$ in $\Cat$.

We would like to reach the conclusion that $\Cat$ is singly aligned,
and also that $\M$ admits least common multiples, starting with the
following apparently weaker concept:

\definition
  \label DefineLittlePullBack
  We shall say that $\Cat$ satisfies the \stress{little pull-back
property} if, given two commuting squares 
  $$
  \matrix{
  &\bullet\cr
  ^{p_1}\nearrow && \searrow^{f_1}\cr
  \bullet\ \ &&\ \ \bullet\cr
  _{q_1}\searrow && \nearrow_{g_1}\cr
  &\bullet
  }
  \qquad  \qquad {\rm and}  \qquad \qquad
  \matrix{
  &\bullet\cr
  ^{p_2}\nearrow && \searrow^{f_2}\cr
  \bullet\ \ &&\ \ \bullet\cr
  _{q_2}\searrow && \nearrow_{g_2}\cr
  &\bullet
  }
  $$
  such that 
  \izitem 
  \zitem all arrows involved are edges,
  \zitem all northeast edges are of the same  color,
  \zitem all southeast edges are of the same color, but not the same
as the northeast ones,
  \medskip\noindent then
  $$
  f_1=f_2  \and g_1=g_2 \quad \Longrightarrow \quad p_1=p_2 \and q_1=q_2.
  $$

It is obvious that a $k$-graph which does not satisfy the little
pull-back property cannot be singly aligned.

Speaking of either one of the diagrams above, say the one on the
left-hand side, one sometimes think of the two-dimensional figure
formed by it as a geometrical representation of the element $f_1p_1$
of $\Cat$.  The algebraic structure of $\Cat$ is based on the idea
that this \stress{square} is determined by the \stress{sides} $p_1$
and $f_1$.  In particular, there cannot be two different squares
sharing these two sides.  A similar observation clearly holds for the
sides $q_1$ and $g_1$.  The little pull-back property goes very much
in this direction by stating that there cannot be two different
squares sharing the sides $f_1$ and $g_1$.  A similar property, which
could be called the \stress{little push-out property}, would say that
two different squares cannot share  the sides $p_1$ and $q_1$.
That property may be shown to imply the existence of push-outs in
$\Cat$.

  \state Proposition
  \label UniqueLongSquares
  Suppose that $\Cat$ satisfies the little pull-back property, and let
$f_i$, $g_i$, $p_i$ and $q_i$ be morphism (rather than edges) such
that $f_ip_i=g_iq_i$, for $i=1,2.$ Suppose also that
  $\rank(f_i)\inf \rank(g_i)=0$, and
  $\rank(p_i)\inf \rank(q_i)=0$,
  for $i=1,2$.  Then the implication at the end of
\lcite{\DefineLittlePullBack} holds true.

\proof  
First observe that
  $$
  \rank(f_i)-\rank(g_i) = \rank(q_i)-\rank(p_i),
  $$
  so one necessarily has
  $$
  \rank(f_i)=\rank(q_i)
  \and 
  \rank(g_i)=\rank(p_i),
  $$
  as a consequence of orthogonality.

Letting $f=f_1=f_2$, and $g=g_1=g_2$, observe that it is enough to show
that $p_1=p_2$, since this would imply that
  $$
  gq_1 = fp_1 = fp_2 = gq_2,
  $$
  and the uniqueness of the factorization would give $q_1=q_2$.

  Suppose first that 
$\rank(f)=0$. Then $\rank(q_1)=\rank(q_2)=0$, so $f$, $q_1$, and $q_2$ are the
identity morphisms on their respective domains.
Therefore
  $$
  p_1 =
  fp_1 =
  gq_1=
  g =
  gq_2=
  fp_2 =
  p_2.
  $$
  A similar argument proves the result if $\rank(g)=0$.  We therefore
suppose, from now on, that $\rank(f)$ and $\rank(g)$ are both nonzero.

We will now proceed by induction on $|\rank(f)|+|\rank(g)|$, observing that
when $|\rank(f)|+|\rank(g)|\leq 1$, the conclusion follows from the above
arguments.   

If $|\rank(f)|+|\rank(g)|=2$, since $|\rank(f)|,|\rank(g)|>0$, we must have that
$|\rank(f)|=|\rank(g)|=1$, and hence $f$ and $g$ are edges. By hypothesis the
$p_i$ and $q_i$ are also edges, so the conclusion follows from the
little pull-back property.

We thus assume that $n\geq2$, and $|\rank(f)|+|\rank(g)|=n+1$.  Hence either
$|\rank(f)|\geq2$, or $|\rank(g)|\geq2$.  Without loss of generality we assume
that $|\rank(f)|\geq2$.  So by the factorization property there are morphisms $f'$
and $f''$ such that $f=f'f''$, and $|\rank(f')|,|\rank(f'')|<|\rank(f)|$.  Since, for $i=1,2$,
  $$
  \rank(q_i) = \rank(f) = \rank(f')+\rank(f''),
  $$
  we may write $q_i=q_i'q_i''$, with $\rank(q_i')=\rank(f')$, and
$\rank(q_i'')=\rank(f'')$.
  We furthermore observe that
  $$
  \rank(fp_i) = \rank(f')+\rank(f'')+\rank(p_i) =
\rank(f')+\rank(p_i)+\rank(q_i''),
  $$
 and hence we factorize $fp_i=\phi_i h_i \psi_i$, with
  $$
  \rank(\phi_i)=\rank(f'),\quad
  \rank(h_i)=\rank(p_i)\and
  \rank(\psi_i)=\rank(q_i'').
  $$
  Notice that
  $$
  f'f''p_i =
  fp_i =
  \phi_i h_i \psi_i.
  $$
  By the uniqueness of the factorization we conclude that
  $$
  f'= \phi_i
  \and
  f''p_i =
  h_i \psi_i.
  \eqmark FirstGreatMess
  $$
  On the other hand, notice that 
  $$
  f' h_i \psi_i =
  \phi_i h_i \psi_i =
  fp_i =
  gq_i = 
  gq_i'q_i''.
  $$
  Again by the uniqueness of the factorization we conclude that
  $$
  f'h_i = gq_i'
  \and 
  \psi_i = q_i''.
  $$

  Observe that $\rank(f')\inf \rank(g) \leq \rank(f)\inf \rank(g) =0$, that
  $\rank(f')=\rank(q_i')$, and that $\rank(h_i)=\rank(g)$.  By the induction
hypothesis we have  that $h_1=h_2$, and $q'_1=q'_2$.
  Let us thus use the simplified notation $h=h_1=h_2$, and $q'=q'_1=q'_2$.
  By  \lcite{\FirstGreatMess} we then deduce that
  $$
  f''p_i =
  h \psi_i.
  $$
  Again we have
  $
  \rank(f'')\inf \rank(h) \leq \rank(f)\inf \rank(g)=0,
  $
  $
  \rank(f'')=\rank(\psi_i), 
  $
  and
  $
  \rank(p_i)=\rank(h).
  $
  By induction we conclude that $p_1=p_2$, and $\psi_1=\psi_2$,
finishing  the proof.
  \proofend

While the result above deals with uniqueness, the next result will
provide existence:

\state Lemma
  \label PullLemma
  Let $f_1,f_2,p_1,p_2$ be morphisms such that $f_1p_1=f_2p_2$.  Then
there are morphisms $r,\bp_1,\bp_2$, such that, for every $i=1,2$, one has
  \izitem
  \zitem $f_1\bp_1=f_2\bp_2$,
  \zitem $p_i=\bp_i r$,
  \zitem $\rank(\bp_1)\inf \rank(\bp_2)=0$,
  \zitem $\rank(f_i\bp_i)=\rank(f_1)\vee \rank(f_2)$.

\proof Since $\rank(f_1),\rank(f_2)\leq \rank(f_1p_1)=\rank(f_2p_2)$, we have that $\rank(f_1)\vee
\rank(f_2)\leq \rank(f_ip_i)$, and hence there are morphisms $s$ and $r$ such that
$sr=f_ip_i$, and $\rank(s) =\rank(f_1)\vee \rank(f_2)$.
  Notice that 
  $$
  \rank(s)+\rank(r) =
  \rank(f_ip_i)  =
  \rank(f_i) + \rank(p_i) \leq 
  \rank(f_1)\vee \rank(f_2) +\rank(p_i) =
  \rank(s) +  \rank(p_i),
  $$
  and hence $\rank(r)\leq \rank(p_i)$.
  By the factorization property we may factor $p_i=\bp_i r_i$, with
$\rank(r_i)=\rank(r)$.  Notice that
  $$
  f_i\bp_ir_i =
  f_ip_i =
  sr.
  $$
  By the uniqueness of the factorization we conclude that
$f_i\bp_i=s$, and $r_i=r$, hence proving (i) and (ii).
  In addition we have 
  $$
  \rank(f_i\bp_i) =
  \rank(s) = 
  \rank(f_1)\vee \rank(f_2),
  $$
  taking care of (iv).  In order to show (iii) suppose that  $n\in\N^k$ is such
that $n \leq \rank(\bp_i)$, for all $i$, then 
  $$
  \rank(f_i) \leq 
  \rank(f_i)+  \rank(\bp_i) - n =
  \rank(f_1)\vee \rank(f_2) -n,
  $$
  whence $\rank(f_1)\vee \rank(f_2) \leq \rank(f_1)\vee \rank(f_2) -n$, so that $n=0$.
  \proofend

The following result can be proved by applying \lcite{\PullLemma} to
the opposite category  $\Cat\,^{op}$.

\state Lemma
  \label PushLemma
  Let $q_1,q_2,g_1,g_2$ be morphisms such that $q_1g_1=q_2g_2$.  Then
there are morphisms $s,\bq_1,\bq_2$, such that, for every $i=1,2$, one has
  \izitem
  \zitem $\bq_1g_1=\bq_2g_2$,
  \zitem $q_i=s\bq_i$,
  \zitem $\rank(\bq_1)\inf \rank(\bq_2)=0$,
  \zitem $\rank(\bq_ig_i)=\rank(g_1)\vee \rank(g_2)$.

\state Proposition
  \label UniqueLongRectangleOne
  Assume that $\Cat$ satisfies the little pull-back property and 
  let $f_1$, $f_2$, $p_1$, $p_2$, $p'_1$ and $p'_2$ be morphisms such
that for all $i=1,2$,
  \izitem
  \zitem $f_1p_1=f_2p_2$, and $f_1p'_1=f_2p'_2$, 
  \zitem $\rank(p_1)\inf \rank(p_2)=0$, and  $\rank(p'_1)\inf \rank(p'_2)=0$.
  \medskip\noindent
  Then $p_i=p'_i$, for $i=1,2$.

\proof
  First observe that 
  $$
  \rank(p_1)-\rank(p_2) = \rank(f_2)-\rank(f_1) =   \rank(p_1')-\rank(p_2'),
  $$
  so $\rank(p_i)=\rank(p_i')$, by (ii).
  Using \lcite{\PushLemma} with $g_i=p_i$, and $q_i=f_i$, let $\F_i$
and $s$ be such that $\F_1p_1=\F_2p_2$, $f_i=s\F_i$, and
$\rank(\F_1)\inf \rank(\F_2)=0$.  Since
  $$
  \rank(\F_1) - \rank(\F_2) = \rank(p_2) - \rank(p_1),
  $$
  we have $\rank(\F_1)=\rank(p_2)$, and $\rank(\F_2)=\rank(p_1)$.

Replacing $p_i$ by $p_i'$ in our application of \lcite{\PushLemma}
just above we would get $\F'_i$ and $s'$ such that
$\F_1'p_1'=\F_2'p_2'$, $f_i=s'\F_i'$, $\rank(\F_1')\inf \rank(\F_2')=0$.  As
above we may also prove that $\rank(\F_1')=\rank(p_2')$, and
$\rank(\F_2')=\rank(p_1')$.  Therefore
  $$
  \rank(\F_1) = \rank(p_2) = \rank(p_2') = \rank(\F_1'),
  $$ and $$
  \rank(\F_2) = \rank(p_1) = \rank(p_1') = \rank(\F_2').
  $$
  Furthermore
  $$
  \rank(s) = \rank(f_i)-\rank(\F_i) = \rank(f_i)-\rank(\F_i') = \rank(s'),
  $$
  and hence the identity $s\F_i =s'\F_i'$, together with the
uniqueness of the factorization gives $s=s'$ and $\F_i=\F_i'$. The
two identities
  $$
  \F_1p_1=\F_2p_2
  \and
  \F_1p_1'=\F_2p_2'
  $$
  and \lcite{\UniqueLongSquares} thus give the conclusion.
  \proofend

So here is the result we were looking for:

  \state Theorem A $k$-graph $\Cat$ satisfying the little pull-back property
is singly aligned and the associated semigroupoid $\M$ admits least
common multiples.

  \proof
  Let $f_1,f_2,p_1,p_2$ be morphisms such that $f_1p_1=f_2p_2$.  Pick
  $r,\bp_1,\bp_2$ as in \lcite{\PullLemma}.  We claim that
$(\bp_1,\bp_2)$ is a pull-back for $(f_1,f_2)$.  

In order to prove
this let $q_1,q_2$ be morphisms such that $f_1q_1=f_2q_2$. Again pick
  $s,\bq_1,\bq_2$ as in \lcite{\PullLemma}, so that
  $f_1\bq_1=f_2\bq_2$,
  $q_i=\bq_i s$, and
  $\rank(\bq_1)\inf \rank(\bq_2)=0$.
  By \lcite{\UniqueLongRectangleOne} we deduce that $\bp_i=\bq_i$, and
hence $q_i=\bp_i s$, as desired.  It is also clear that $s$ is unique
by the factorization property.  If then immediately follows that
$\Cat$ is singly aligned.  That $\M$ admits least common multiples is
then a consequence of \lcite{\SGPDFromCatego.iii}.
  \proofend

The little pull-back property is the last restriction we need to
impose on $\Cat$ in order to be able to apply all of the conclusions
of \lcite{\SGPDFromCatego}.

In view of \scite{\Muhly}{3.8.(3)} it is reasonable to restricts one's
attention to representations of $\M$ which respects least common
multiples.  The following is the main result of this section.

\state Theorem 
  \label HHGrap
  Let $\Cat$ be a countable $k$-graph satisfying the little pull-back
property and such that for every vertex $v$ there is some morphism $f$,
other than the identity on $v$, with $\r(f)=v$.  Then, removing the
identities from $\Cat$ we obtain a semigroupoid $\M$ which has no
springs, contains only monic elements and in which every intersecting
pair of elements admits a least common multiple.  Moreover the the
C*-algebra generated by the range of a universal tight representation
of\/ $\M$, respecting least common multiples, is naturally isomorphic to
the C*-algebra of the groupoid $\G_\M$ of germs for the standard action
of\/ $\SM$ on the tight part of the spectrum of its idempotent
semilattice.

\proof
  Follows from \lcite{\SGPDFromCatego} and \lcite{\CstarAlgOfGpd}.
  \proofend

\cryout{From now we fix a $k$-graph $\Cat$ satisfying the hypothesis of
\lcite{\HHGrap}.}  

To conclude this section we will give a description 
of $\spect$, where $\E$ is the idempotent
semilattice of $\SM$.  Given a path $\stem$ on $\M$, let $f,g\in\stem$
with $\rank(f)=\rank(g)$.  Since $f\its g$, by \lcite{\DefinePath.ii}
we may write $fu=gv$.  Extending $\rank$ to $\Mu$ by defining
$\rank(1)=0$, we then have that $\rank(u)=\rank(v)$, and then $f=g$ by
the unique factorization property.  This says that $\stem$ may contain
at most one element $f$ with $\rank(f)=n$, for each $n\in\N^k$.

\state Proposition
  Given a nonempty path $\stem$ on $\M$, let $D$ be the image of $\stem$
under the rank function $\rank$, and for each $n\in D$, let
$\mu(n)$ be the unique element $f$ in $\stem$ with $\rank(f)=n$.
Then
  \izitem
  \zitem $D\cup\{0\}$ is a hereditary subset of\/ $\N^k$,
  \zitem if $n,m\in D$, then $n\vee m\in D$,
  \zitem if $n,m\in D$, and $n\leq m$, then $\mu(n)\dil
\mu(m)$,
  \zitem $\stem =\{\mu(n): n\in D\}$.

\proof
  Let $n,m\in N^k$, with $m\in D$, and $0\neq n\leq m$.  Set
$f=\mu(m)$, so that $\rank(f)=m$.  Writing $m = n + (m-n)$, the unique
factorization property implies that $f=gh$, with $\rank(g)= n$, and
$\rank(h)=m-n$.  Since $g\dil f$ we conclude that $g\in\stem$, and
hence $n\in D$, proving (i).
  It is also clear that $g=\mu(n)$, so (iii) is also proved.  To prove
(ii) let $n,m\in D$, so $f:= \lcm(\mu(n),\mu(m))\in\stem$, and hence
$\rank(f)\in D$.  It may be proved that $\rank(f)=m\vee n$, but it
suffices to notice that, since $f$ is a common multiple of $\mu(n)$
and $\mu(m)$, one has that
  $n,m\leq \rank(f)$, and consequently 
  $n\vee m\leq \rank(f)$.  Thus (ii) follows from (i).  The last point
is trivial.
  \proofend

The following is a converse to the above:

\state Proposition
  \label GeneralPath
  Let $D$ be a subset of\/ $\N^k$ not containing 0, but such that
$D\cup\{0\}$ is a hereditary subset of\/ $\N^k$.  Assume that $D$ is
closed under ``$\,\vee$" and let
  $\stem:D \to \M$ be any map such that for every $n,m\in D$,
  \izitem 
  \zitem $\rank(\mu(n))=n$,
  \zitem $\mu(n)\dil \mu(m)$, if $n\leq m$,
  \medskip\noindent
  Then the set $\stem = \{\mu(n): n\in D\}$ is a path in $\M$.

\proof
  If $f\in\M$, and $f\dil \mu(m)$, for some $m\in D$, write
$\mu(m)=fu$, for some $u\in\Mu$.  This clearly implies that
  $$
  n:=\rank(f)\leq\rank(\mu(m)) = m,
  $$
  so $n\in D$, and $\mu(m) = \mu(n)v$, for some $v\in\Mu$.  By the
unique factorization we have that $f=\mu(n)\in\stem$.

To prove \lcite{\DefinePath.ii} suppose that $n,m\in D$.  Then
$\mu(n\vee m)$ is a common multiple of $\mu(n)$ and $\mu(m)$ and,
recalling that under our assumptions $\M$ admits least common
multiples, we have
  $$
  \lcm(\mu(n),\mu(m)) \dil \mu(n\vee m).
  $$
  However it is easy to see that
$\rank\big(\lcm\big(\mu(n),\mu(m)\big)\big) \geq n\vee m$, so
  $$
  \lcm(\mu(n),\mu(m)) = \mu(n\vee m) \in\stem.
  \proofend
  $$

Notice that for any set $D$ as above one may define the supremum of
$D$ as an element 
  $$
  m \in(\N\cup\{\infty\})^k, 
  $$
  and   hence
  $D = \Omega_{k,m} := \{n\in\N^k: n\leq m\}$, as Defined in
\scite{\Muhly}{3.2}.

It therefore follows that paths in $\M$ correspond to maps $\mu$, as
in \lcite{\GeneralPath}, and hence also to the usual notion of paths
in higher rank graphs \scite{\Muhly}{5.1}.  One may then use
\lcite{\CharacTightInCatego} to relate elements of $\spect$ to the
boundary paths of \cite{Muhly}.

  \endgroup

  \begingroup
  \bigskip\bigskip 
  \font \sc = cmcsc8 \sc
  \parskip = -1pt

  Departamento de Matem\'atica 

  Universidade Federal de Santa Catarina

  88040-900 -- Florian\'opolis -- Brasil

  \eightrm r@exel.com.br

  \endgroup
  \end
\references


  \bibitem{\AbadieGroupoid}
  {F. Abadie}
  {On partial actions and groupoids}
  {\it Proc. Amer.  Math. Soc., \bf 132 \rm (2004), 1037--1047
(electronic)}

  \bibitem{\BHRS}
  {T. Bates, J. Hong, I. Raeburn and W. Szyma{\'n}ski}
  {The ideal structure of the C*-algebras of infinite graphs}
  {\it Illinois J. Math., \bf 46 \rm (2002), 1159--1176}

  \bibitem{\BPRS}
  {T. Bates, D. Pask, I. Raeburn and W. Szyma{\'n}ski}
  {The C*-algebras of row-finite graphs}
  {\it New York J. Math., \bf 6 \rm (2000), 307--324 (electronic)}

  \bibitem{\Connes}
  {A. Connes}
  {A survey of foliations and operator algebras}
  {Operator algebras and applications, Part I (Kingston, Ont.,1980),
\it Proc. Sympos. Pure Math., \bf 38 \rm (1982), 521--628}

  \bibitem{\Cuntz}
  {J. Cuntz}
  {Simple $C^*$-algebras generated by isometries}
  {\sl Comm. Math. Phys. \bf 57 \rm (1977), 173--185}

  \bibitem{\CK}
  {J. Cuntz and W. Krieger}
  {A Class of C*-algebras and Topological Markov Chains}
  {\sl Inventiones Math., \bf 56 \rm (1980), 251--268}

  \bibitem{\Watatani}
  {M. Enomoto, Y. Watatani}
  {A graph theory for C*-algebras}
  {\it Math. Japon., \bf 25 \rm (1980), 435--442}

  \bibitem{\newpim}
  {R. Exel}
  {Circle actions on C*-algebras, partial automorphisms and a
generalized Pimsner--Voiculescu exact sequence}
  {\it J. Funct. Analysis, \bf 122 \rm (1994), 361--401}

  \bibitem{\tpa}
  {R. Exel}
  {Twisted partial actions, a classification of regular C*-algebraic
bundles}
  {\it Proc. London Math. Soc., \bf 74 \rm (1997), 417--443}

  \bibitem{\ExSemiGpdAlg} 
  {R. Exel}
  {Semigroupoid C*-Algebras}
  {preprint, Universidade Federal de Santa Catarina, 2006,
[arXiv:math.OA/0611929]}

  \bibitem{\ExelAlgebra}
  {R. Exel}
  {Tight representations of semilattices and inverse semigroups}
  {in preparation}

  \bibitem{\infinoa}
  {R. Exel and M. Laca}
  {Cuntz--Krieger algebras for infinite matrices}
  {\it J. reine angew. Math. \bf 512 \rm (1999), 119--172}

  \bibitem{\Muhly}
  {C. Farthing, P.  Muhly, and T. Yeend}
  {Higher-rank graph C*-algebras: an inverse semigroup and groupoid
approach}
  {\it Semigroup Forum, \bf 71 \rm (2005), 159--187}

  \bibitem{\FLR}
  {N. Fowler, M. Laca, and I. Raeburn}
  {The C*-algebras of infinite graphs}
  {\it Proc. Amer. Math. Soc., \bf 128 \rm (2000), 2319--2327}

  \bibitem{\KPActions}
  {A. Kumjian and D. Pask}
  {$C\sp *$-algebras of directed graphs and group actions}
  {\it Ergodic Theory Dynam. Systems, \bf 19 \rm (1999), 1503--1519}

  \bibitem{\KP}
  {A. Kumjian and D. Pask}
  {Higher-rank graph C*-algebras}
  {\it New York J. Math., \bf 6 \rm (2000), 1--20 (electronic)}

  \bibitem{\KPZActions}
  {A. Kumjian and D. Pask}
  {Actions of $Z\sp k$ associated to higher rank graphs}
  {\it Ergodic Theory Dynam. Systems, \bf 23 \rm (2003), 1153--1172}

  \bibitem{\KPR}
  {A. Kumjian, D. Pask and I. Raeburn}
  {Cuntz-{K}rieger algebras of directed graphs}
  {\it Pacific J. Math., \bf 184 \rm (1998), 161--174}

  \bibitem{\KPPR}
  {A. Kumjian, D. Pask, I. Raeburn and J.  Renault}
  {Graphs, groupoids, and Cuntz-Krieger algebras}
  {\it J. Funct. Anal., \bf 144 \rm (1997), 505--541}

  \bibitem{\Lawson}
  {M. V. Lawson}
  {Inverse semigroups, the theory of partial symmetries}
  {World Scientific, 1998}

  \bibitem{\McClanahan}
  {K. McClanahan}
  {$K$-theory for partial crossed products by discrete groups}
  {\it J. Funct. Anal., \bf 130 \rm (1995), 77--117}

  \bibitem{\PQR}
  {D. Pask, J. Quigg and I. Raeburn}
  {Fundamental groupoids of {$k$}-graphs}
  {\it New York J. Math., \bf 10 \rm (2004), 195--207 (electronic)}

  \bibitem{\PRRS}
  {D. Pask, I. Raeburn, M. Rordam and A. Sims}
  {Rank-two graphs whose C*-algebras are direct limits of circle algebras}
  {\it J. Funct. Anal., \bf 239 \rm (2006), 137--178}

  \bibitem{\PatBook} 
  {A. L. T. Paterson}
  {Groupoids, inverse semigroups, and their operator algebras}
  {Birkh\"auser, 1999}

  \bibitem{\PatGraph}
  {A. L. T. Paterson}
  {Graph inverse semigroups, groupoids and their C*-algebras}
  {\it J. Operator Theory, \bf 48 \rm (2002), 645--662}


  \bibitem{\QuigSieb}
  {J. Quigg and N. Sieben}
  {C*-actions of {$r$}-discrete groupoids and inverse semigroups}
  {\it J. Austral. Math. Soc. Ser. A, \bf 66 \rm (1999), 143--167}

  \bibitem{\RaeBook}
  {I. Raeburn}
  {Graph algebras}
  {CBMS Regional Conference Series in Mathematics, \bf 103 \rm (2005),
pp. vi+113}

  \bibitem{\RSY}
  {I. Raeburn, A. Sims and T. Yeend}
  {The C*-algebras of finitely aligned higher-rank graphs}
  {\it J. Funct. Anal., \bf 213 \rm (2004), 206--240}

  \bibitem{\RaeSzy}
  {I. Raeburn and W. Szyma\'nski}
  {Cuntz-Krieger algebras of infinite graphs and matrices}
  {\it Trans. Amer. Math. Soc., \bf 356 \rm (2004), 39--59
(electronic)}

  \bibitem{\PTW}
  {I. Raeburn, M. Tomforde and D. P. Williams}
  {Classification theorems for the C*-algebras of graphs with sinks}
  {\it Bull. Austral. Math. Soc., \bf 70 \rm (2004), 143--161}

  \bibitem{\RenaultThesis}  
  {J. Renault}
  {A groupoid approach to $C^*$-algebras}
  {Lecture Notes in Mathematics vol.~793, Springer, 1980}

  \bibitem{\RenaultInfinoa}
  {J.  Renault}
  {Cuntz-like algebras}
  {Operator theoretical methods (Timisoara, 1998), 371--386, Theta
Found., Bucharest, 2000}


  \bibitem{\Sieben}
  {N. Sieben}
  {C*-crossed products by partial actions and actions of inverse
semigroups}
  {\it J. Austral. Math. Soc. Ser. A, \bf 63 \rm (1997), 32--46}

  \bibitem{\Szendrei}
  {M. B. Szendrei}
  {A generalization of McAlister's P-theorem for E-unitary regular
semigroups}
  {\it Acta Sci. Math., \bf 51 \rm (1987), 229--249}

  \bibitem{\Tomforde}
  {M. Tomforde}
  {A unified approach to Exel-Laca algebras and C*-algebras
associated to graphs}
  {\it J. Operator Theory, \bf 50 \rm (2003), 345--368}
